\documentclass{cedram-aif}

\usepackage[english,frenchb]{babel}

\usepackage[latin1]{inputenc}
 \usepackage[T1]{fontenc}

\usepackage{amsmath,amssymb}
\usepackage{cite}
% ...

\newcommand\SL{\mathord{\mathsf{SL}}}

\newcommand\rank{{\mathord{\mathrm{rk}}}}
% ...

\title
[Kac-Moody groups, hovels and Littelmann paths]
{Kac-Moody groups, hovels and Littelmann paths}

\alttitle{Kac-Moody groups, hovels and Littelmann paths}

\author{\firstname{Stéphane} \lastname{Gaussent}}
\author{\firstname{Guy} \lastname{Rousseau}}
\address{Institut Elie Cartan\\ 
Unité Mixte de Recherche 7502\\
Nancy-Université, CNRS, INRIA\\ 
Boulevard des Aiguillettes\\ 
B.P. 239\\ 
F-54506 Vand\oe uvre-lès-Nancy\\}

\email{gaussent@iecn.u-nancy.fr, rousseau@iecn.u-nancy.fr}

  

\begin{document}
\begin{abstract}

Nous d\'efinissons une sorte d'immeuble $\mathcal I$ associ\'e \`a un groupe de Kac-Moody sym\'etrisable sur un corps $K$ muni d'une valuation discr\`ete avec un corps r\'esiduel contenant $\mathbb C$. Nous l'appelons masure (hovel) \`a cause de l'absence d'une propri\'et\'e importante des immeubles. Cependant, de bonnes propri\'et\'es restent, par exemple l'existence de retractions  de centre un germe de quartier. Cela nous permet de g\'en\'eraliser plusieurs r\'esultats prouv\'es par S. Gaussent et P. Littelmann dans le cas semi-simple. En particulier, si
$K=\mathbb C(\!(t)\!)$, les segments g\'eod\'siques dans $\mathcal I$, d'extr\'emit\'e un sommet sp\'ecial et se r\'etractant sur un chemin donn\'e $\pi$, sont param\'etr\'es par un ouvert de Zariski $P$ de $\mathbb C^N$. Cette dimension $N$ est maximale quand $\pi$ est un chemin LS et alors $P$ est fortement associ\'e \`a un cycle de Mirkovi\'c-Vilonen.
\end{abstract}

\begin{altabstract}
  We give the definition of a kind of building $\mathcal I$ for a 
symmetrizable Kac-Moody group over a field $K$ endowed with a discrete
valuation and with a residue field containing $\mathbb C$. Due to the lack of some important property of buildings,
we call it a hovel. Nevertheless, some good ones remain, for example, the
existence of retractions with center a sector-germ. This enables us to
generalize many results proved in the semisimple case by S. Gaussent and P.
Littelmann. In particular, if
$K=\mathbb C(\!(t)\!)$, the geodesic segments in $\mathcal I$, ending in a special vertex and 
retracting onto a given path $\pi$, are parametrized by a Zariski open subset
$P$ of $\mathbb C^N$. This dimension $N$ is maximal when $\pi$ is a LS path and then
$P$ is closely related to some Mirkovi\'c-Vilonen cycle.
\end{altabstract}

\maketitle

\section{Introduction}
\label{se:Intro}

Let $\mathfrak g^\vee$ be a complex symmetrizable Kac-Moody algebra. To capture 
the combinatorial essence of the representation theory of $\mathfrak g^\vee$, P.
{Littelmann} \cite{Littelmann94, Littelmann95} introduced the path model. Particularly, this model gives
a method to compute the multiplicity of a weight $\mu$ in an irreducible
representation of highest weight
$\lambda$ (a dominant weight), by counting some ``Lakshmibai-Seshadri'' (or LS) paths  of
shape $\lambda$ starting from $0$ and ending in $\mu$. When $\mathfrak g^\vee$ is semi-simple and
$G$ is an algebraic group with Lie algebra the Langlands dual $\mathfrak g$ of $\mathfrak g^\vee$, I.
Mirkovi\'c and K. Vilonen \cite{MirkovicVilonen00} gave a new interpretation of this multiplicity: it is
the number of irreducible components (the MV cycles) in some subvariety $X^\mu_\lambda$ of
the affine grassmannian $\mathcal G=G(\mathbb C(\!(t)\!)\,)/G(\mathbb C[[t]]\,)$.

 S. Gaussent and P. Littelmann \cite{GaussentLittelmann05} gave a link between these two theories  (when
$G$ is semi-simple). Actually, the LS paths are drawn in a vector space~$V$ which is an
apartment $A$  of the Bruhat-Tits building of $G$ (over any non archimedean valuated
field $K$, in particular $K=\mathbb C(\!(t)\!)$). They replaced  the LS paths of shape $\lambda$ from
$0$ to $\mu$ by ``LS galleries'' of type $\lambda$ from $0$ to $\mu$. This gives a new
``gallery model'' for the representations of~$\mathfrak g^\vee$. Moreover, let $\rho$ be the retraction of the Bruhat-Tits building $\mathcal I$  of
$G$ over $K=\mathbb C(\!(t)\!)$ onto $A$ with center some sector-germ $\mathfrak S_{-\infty}$ in $A$.
Then, the image under $\rho$ of a minimal gallery $\Gamma$ of type $\lambda$ starting from $0$
in $\mathcal I$ is a gallery $\gamma$ in $A$ of type~$\lambda$ which looks much like a LS gallery: it
is ``positively folded''. Conversely, any positively folded gallery $\gamma$ in $A$ of type
$\lambda$ from $0$ to
$\mu$ is the image under $\rho$ of many minimal galleries $\Gamma$ in $\mathcal I$. These 
galleries are parametrized by a complex variety $X_\gamma$ combinatorially defined from~$\gamma$. Moreover, $\gamma$ is a LS gallery if, and only if, $\dim(X_\gamma)$ is maximal, and then~
$X_\gamma$ is isomorphic to an open subset of a MV cycle in $X^\mu_\lambda$.

 It was natural (and suggested to us by P. Littelmann) to try to generalize this  
when $G$ is a Kac-Moody group. Actually, G. Rousseau \cite{Rousseau06} had constructed some building for a Kac-Moody group over a discretely valuated field $K$. But, the apartments of this microaffine
building are not appropriate to define LS paths. So, we construct a new set $\mathcal I$ (see \ref{de:Hovel})
associated to the Kac-Moody group $G$ over $K=\mathbb C(\!(t)\!)$ (or, more generally, over any discretely valuated field $K$ with residue field containing $\mathbb C$). By definition, the group $G(K)$ acts on $\mathcal I$, LS paths can be drawn on its apartments and the action of $G(K)$ on them is transitive.
Unfortunately, any two points in $\mathcal I$ are not always in a same apartment: this was
already noticed (in a different language and in the affine case) by H. Garland \cite{Garland95}, who remarked that Cartan decomposition is true only after some twist (cf. Remark~\ref{re:Preorder}). Because of this pathological behaviour, $\mathcal I$ is called a {\it hovel}. Moreover, the system of walls in an apartment of $\mathcal I$ is not discrete, so the notion of chamber in $\mathcal I$ is unusual (but follows an idea of F. Bruhat and J. Tits \cite{BruhatTits72}) and $\mathcal I$ is not gallery-connected (see Section~\ref{sse:AffineApart}). Therefore, we have to come back to the path model.

 Nevertheless, we get good generalizations of Gaussent-Littelmann's results. 
First, any sector-germ and any point in $\mathcal I$ are always in a same apartment (this is
equivalent to Iwasawa decomposition, Proposition~\ref{pr:Iwasawa}). Next, we fix a maximal torus in $G$ and a system of positive roots, this gives us an apartment $A$ in $\mathcal I$ and a sector-germ $\mathfrak S_{-\infty}$ in $A$. So, we get a retraction $\rho$ of $\mathcal I$ onto $A$ with center $\mathfrak S_{-\infty}$ (\ref{sse:Retraction}). Following a definition of M. Kapovich and J. Millson \cite{KapovichMillson05}, we say that a
Hecke path is a piecewise linear path $\pi_1$ in~$A$ which is positively folded along
true walls (see Definition~\ref{de:Hecke}). Now, an analogue of some results due,
in the semisimple case, to Kapovich-Millson or Gaussent-Littelmann may be proven:

\begin{theo}[see \ref{th:RetractingSeg}]\label{theo11}
If $\pi$ is a geodesic segment (of shape $\lambda$) in  $\mathcal I$, then $\rho\pi$ is a Hecke path (of shape $\lambda$) in $A$.
\end{theo}

 Conversely, any Hecke path $\pi_1$ in $A$ is the image under $\rho$ of a geodesic
segment $\pi$ in $\mathcal I$. But, if we want a finite dimensional variety of parameters, we 
can no longer look at segments with a given starting point but rather at
segments with a given end. We get (when $K=\mathbb C(\!(t)\!)$):

\begin{theo}[see \ref{th:RetractingOnHecke}]\label{theo12}
Let $\pi_1$ be a Hecke path of shape $\lambda$ in $A$ with endpoint a special vertex $y$. Then, there exist geodesic segments $\pi$ in
$\mathcal I$ with endpoint $y$ such that $\rho\pi=\pi_1$ and they are parametrized by a Zariski
open subset $P(\pi_1,y)$ of $\mathbb C^N$, stable under the natural action of $(\mathbb C^*)^N$.
\end{theo}

Here, $N$ is the so-called {\it dual dimension} of $\pi_1$ (\ref{de:DdimCodim}) and it is maximal (among 
Hecke paths of shape $\lambda$ with the same starting and ending points) if and only if
$\pi_1$ is a LS path.

 This result enables us to state that $P(\pi_1,y)$ is isomorphic to a dense open 
subset of some Mirkovi\'c-Vilonen cycle. In the semi-simple case, this MV cycle is, up
to isomorphism, the classical one associated to the reverse path of $\pi_1$.

 The paper is organized as follows. In Section~\ref{se:KMgroupsAndApart}, we recall some results 
on Kac-Moody groups and their affine apartments. Actually, in the literature one finds
many kinds of Kac-Moody groups. We choose the  minimal one, the most algebraic. But (in Section~\ref{sse:MaxKM}), we will also have to use the maximal one which appears to be a
formal completion of the minimal one and has better commutation relations. 

The construction of the hovel $\mathcal I$ is explained completely  in Section~\ref{se:Hovel}
and the first properties are developed in Section~\ref{se:HovelProperties}. The proofs are rather
involved but we get all what is needed: in particular, the Iwasawa
decomposition (\ref{pr:Iwasawa}), the retraction with respect to a sector-germ (\ref{sse:Retraction}), the twin building structure of the residue of $\mathcal I$ at some point $x$ (\ref{sse:Residue}) and, for some
``good'' subsets $\Omega$ in apartments, the structure of their fixator ({\it i.e.} pointwise
stabilizer) $G_\Omega$ and the transitivity of the action of 
$G_\Omega$ on apartments containing $\Omega$ (\ref{sse:Good} to \ref{sse:Applications}).

In Section~\ref{se:Littelmann}, we give the definitions of LS paths, Hecke paths, dual dimension and codimension. We prove the characterization of LS paths as Hecke paths with maximal dual dimension
(resp. minimal codimension).

 We get in Section \ref{se:Segments} the results on paths explained above, in particular 
Theorems \ref{theo11} and \ref{theo12}. The last theorem (Theorem~\ref{th:Preorder}) asserts that there is, on $\mathcal I$, a
preorder relation which induces, on each apartment, the preorder given by the Tits cone.

 We thank Peter Littelmann for his suggestion to look at these problems and
for some interesting discussions. We also thank Michel Brion for his careful
reading of a previous version of the present paper and his comments.

\section{Kac-Moody groups and the apartment}
\label{se:KMgroupsAndApart}

 We recall here the main results on Kac-Moody groups (in \ref{sse:KMgroups}). A good 
reference is \cite{Remy02}, see also \cite{Rousseau06}. We introduce the model apartment of
our hovel and call it, in analogy with the classical case, the affine apartment (see~\ref{sse:AffineApart}). Because the set of walls is not locally finite anymore, the definition of faces
needs the notion of filter (see Sections \ref{ssse:Faces} to \ref{ssse:Sectors}).

\subsection{Kac-Moody groups} 
\label{sse:KMgroups}

\subsubsection{Kac-Moody algebras}
\label{ssse:KMalg}

A {\it Kac-Moody matrix} (or generalized Cartan matrix)
 is a square matrix $A = (a_{i,j})_{i,j\in I}$, with integer coefficients, 
indexed by a finite set $I$ and satisfying:
\begin{itemize}
\item[(i)] $a_{i,i} = 2 \quad \forall i \in I$,
 \item[(ii)] $a_{i,j} \leq 0 \quad \forall i \not = j $,
\item[(iii)] $a_{i,j} = 0 \iff a_{j,i} = 0 $.
\end{itemize}
A {\it root generating system} \cite{Bardy96} is a 5-tuple 
${\mathcal S}=(A,X,Y,(\alpha_i)_{i\in I}, (\alpha_i^\vee)_{i\in I})$ made of a 
Kac-Moody matrix $A$ indexed by $I$, of two dual free $\mathbb Z-$modules  $X$ (of 
characters) and
 $Y$ (of cocharacters) of finite rank $\rank (X)$, a family $(\alpha_i)_{i\in I}$ 
(of simple roots) in $X$ and a family  $(\alpha_i^\vee)_{i\in I}$ (of simple 
coroots) in $Y$. They have to satisfy the 
following compatibility condition: $\;$ $a_{i,j} = \alpha_j(\alpha_i^\vee) $.

 The {\it  Langlands dual} of $\mathcal S$ is  ${\mathcal S}^\vee=(^tA,Y,X,(\alpha_i^\vee)_{i\in I} ,(\alpha_i)_{i\in I})$, where $^tA$ is the transposed matrix of $A$.

The {\it Kac-Moody algebra} $\mathfrak g=\mathfrak g_{\mathcal S}$ is a complex Lie algebra 
generated by the standard Cartan subalgebra $\mathfrak h=Y\otimes\mathbb C$ and the Chevalley 
generators $(e_i)_{i\in I}$, $(f_i)_{i\in I}$; we shall not explain here the relations, 
see for instance \cite{Kac90}.

 The adjoint action of $\mathfrak h$ on $\mathfrak g$ gives a grading on $\mathfrak g$: 
$\mathfrak g= \mathfrak h\oplus(\oplus_{\alpha\in\Delta}\;\mathfrak g_\alpha)$, where 
$\Delta\subset X\setminus\{0\}\subset\mathfrak h^*$ is the set of {\it roots} of $\mathfrak g$ 
(with respect to $\mathfrak h$). For all $i\in I$, $\mathfrak g_{\alpha_i}=\mathbb C e_i$ and 
$\mathfrak g_{-\alpha_i}=\mathbb C f_i$. If $Q^+ = \sum_i\; \mathbb N\alpha_i$, $\Delta^+ = \Delta\cap Q^+$ 
and $\Delta^- = - \Delta^+$, one has $\Delta = \Delta^+ \bigsqcup \Delta^-$.

\begin{enonce*}[remark]{N.B} For simplicity we shall assume throughout the paper the 
following condition: 
$$F(\mathcal S)\; \hbox{The family}\ (\alpha_i)_{i\in I}\ \hbox{is free in}\ X\ \hbox{and the
family}\  (\alpha_i^\vee)_{i\in I}\  \hbox{is free in}\  Y.$$
See Section \ref{ssse:Freedom}.
\end{enonce*}

Starting from Section \ref{sse:Parahoric}, we shall also assume that the Lie algebra $\mathfrak g$ is symmetrizable, that is, endowed with a nondegenerate invariant $\mathbb C-$valued symmetric bilinear form.

\subsubsection{Weyl group and real roots} 
\label{ssse:WeylGroup}

 Let $V = Y\otimes\mathbb R\subset \mathfrak h$; every element in $X$ defines a linear form
 on this $\mathbb R-$vector space. For $i \in I$, the formula 
$r_i(v) = v - \alpha_i(v)\alpha_i^\vee $ defines an involution in $V$ (or $\mathfrak h$), 
more precisely a reflection of hyperplane Ker$(\alpha_i)$.

 The (vectorial) Weyl group $W^v$ is the subgroup of $GL(V)$ generated by  the set
$\{r_i\}_{i\in I}$. One knows that $W^v$ is a Coxeter group; it stabilizes the lattice
 $Y$ of $V$. It also acts on $X$ and stabilizes $\Delta$.

 One denotes by $\Phi=\Delta_{\rm re}$ the set of real roots, those which 
can be written as $\alpha = w(\alpha_i)$ with $w \in W^v$ and $ i \in I$. This set 
$\Phi$  is infinite except in the classical case, where $A$ is a Cartan matrix and 
$\mathfrak g$ is finite-dimensional (reductive). 
If $\alpha \in \Phi$, then $r_{\alpha} = w.r_i.w^{-1} $ is well determined 
by $\alpha$, independently of the choice of $w$ and of $i$ such that $\alpha = 
w(\alpha_i)$.
 For $v\in V$ one has $r_{\alpha}(v) = v - \alpha(v)\alpha^\vee $, where the coroot
$\alpha^\vee \in Y $ associated to $\alpha$ satisfies $\alpha(\alpha^\vee) = 2 $.
 Hence $r_{\alpha}$ is the 
reflection with respect to the hyperplane $M(\alpha) = $ Ker$(\alpha)$ which is  called
the  {\it wall} of $\alpha$. The {\it half-apartment} associated to $\alpha$ is 
$D(\alpha) = \{ v \in V \mid \alpha(v) \geq 0 \}$.

 The set $\Phi$ is a system of (real) roots in the sense of  
\cite{MoodyPianzola95}. The set $\Delta$ is a system of roots in the sense of  
\cite{Bardy96}. The imaginary roots (those in $\Delta_{\rm im}=\Delta\setminus\Phi$) 
will not be very much used here. We define $\Phi^\pm =\Phi\cap\Delta^\pm $.

 A subset $\Psi$ of $\Phi$ (or $\Delta$) is said to be {\it closed} in $\Phi$ (or 
$\Delta$)  if: $\alpha , \beta \in \Psi$, $\alpha + \beta \in \Phi$ (or $\Delta)
\Rightarrow \alpha + \beta \in \Psi$. The subset $\Psi$ is said to be {\it prenilpotent} if there
 exist $w, w' \in W^v$ such that $w\Psi \subset \Delta^+$ and $w'\Psi \subset
 \Delta^-$. Then $\Psi$ is finite and contained in the subset
$w^{-1}(\Phi^+) \cap (w')^{-1}(\Phi^-)$ which is {\it nilpotent} ({\it i.e.} 
prenilpotent and closed).

 One denotes by $Q^\vee$ (resp. $P^\vee, Q$ ) the {\it 
``coroot-lattice''}  (resp. {\it ``coweight-lattice''}, {\it ``root-lattice''}), {\it i.e.} the
subgroup of $Y$ generated by the 
$\alpha_i^\vee$ (resp. $P^\vee = \{ y\in Y\otimes\mathbb Q \, \mid \, \alpha_i(y) 
\in \mathbb Z , \forall i \in I\, \}$, $Q = \sum_i\; \mathbb Z\alpha_i$); one has  $Q^\vee
\subset  Y \subset P^\vee$. Actually, $Q^\vee$, $P^\vee$ or $Q$ is a lattice in $V$ or
$V^*$ if and only
 if the $\alpha_i^\vee$ generate $V$ {\it i.e.} the ${\alpha}{_i}$ generate $V^*$ {\it i.e.}
$\vert I\vert = \rank(X) = \dim(V)$. We define the set of {\it dominant weights} 
$X^+=\{\chi\in X \mid \chi(\alpha_i^\vee)\geq 0 , \forall i\in I\}$ and $X^-=-X^+$. 
Dually,  the set of {\it dominant coweights} is
$Y^+=\{\lambda\in Y \mid \lambda(\alpha{_i})\geq 0 , \forall i\in I\}$ and $Y^-=-Y^+$.

\subsubsection{The Tits cone} 
\label{ssse:Titscone}

 The {\it positive fundamental chamber} $C^v_f = \{ u\in V\mid\alpha_i(u) > 0
 \quad\forall i \in I \}$ is a nonempty open convex cone. Its closure 
$\overline{C^v_f}$ is the 
disjoint union of the {\it faces} $ F^v(J) =$  $\{\; u\in V\;\mid \;\alpha_i(u) = 0 
\; \forall i \in J \,; \; \alpha_i(u) > 0 \;	\forall i \notin J \} $ for $J\subset I$;
 one has $C^v_f = F^v(\emptyset)$. We define $V_0=F^v(I)$, it is a vector subspace. These  faces
are called vectorial because they are convex  cones with base point $0$. One says that
the face $F^v(J)$ or the set
$J$ is {\it spherical} (or {\it of finite type}) if the matrix  $A(J) =
(a_{i,j})_{i,j\in J}$ is a  Cartan  matrix (in the classical sense), i.e. if $W^v(J) =
\langle r_i  \mid i \in J \rangle$ is finite.  This holds for the chamber $C^v_f$  or
its {\it panels} $F^v(\{i\})$, $\forall i \in I$.

 The {\it Tits cone} is the union $\mathcal T$ of the positive {\it closed-chambers}
 $w.\overline{C^v_f}$ for $w\in W^v$. Its interior is the {\it open Tits cone}
 ${\mathcal T}^o$, disjoint union of the  (positive) {\it spherical faces}
 $w.F^v(J)$ for $w$ in $W^v$ and $J$ spherical. Both $\mathcal T$, ${\mathcal T}^o$ and
their closure $\overline{\mathcal T}$ are  
convex cones, stable under $W^v$. They may be defined as: 
\begin{align*}
{\mathcal T}&=\{v\in V\mid \alpha(v)<0\ \hbox{only for a finite number of }\alpha\in\Delta^+  (\hbox{or }
\Phi^+)\;\},\\
{\mathcal T}^o&=\{v\in V\mid \alpha(v)\leq 0\ \hbox{only for a finite number of } \alpha\in\Delta^+ (\hbox{or }
\Phi^+)\;\},\\
\overline{\mathcal T}&=\{v\in V\mid \alpha(v)\geq 0\quad\forall\alpha\in\Delta^+_{\rm im}\;\}.
\end{align*}

 The action of $W^v$ on the positive chambers is simply transitive. The fixator 
(pointwise stabilizer)  or the stabilizer of $F^v(J)$ is $W^v(J)$.

 We shall also consider the negative Tits cones $-\mathcal T$, $-{\mathcal T}^o$, 
 $-\overline{\mathcal T}$ and all negative faces, 
chambers... which are obtained by change of sign.

 Actually, ${\mathcal T}^o\cap-{\mathcal T}^o=\emptyset$ except in the classical case 
(where ${\mathcal T}^o=-{\mathcal T}^o=V$) and ${\mathcal T}\cap-{\mathcal T}=\{v\in V\;\mid\;
\alpha(v)=0$ for almost all $\alpha\in\Phi$ (or $\Delta)\}$ is reduced to 
$V_0=\bigcap_{\alpha\in\Delta}\;$Ker$(\alpha)$ if no connected component of $I$ is spherical.

\subsubsection{The Kac-Moody groups} 
\label{ssse:KMgroups}

 One considers the (split, complex) Kac-Moody group $G=G_{\mathcal S}$ associated 
to the above root generating system as defined by Tits \cite{Tits87}, see also 
\cite[Chapitre 8]{Remy02}. It is actually an affine ind-algebraic-group \cite[7.4.14]{Kumar02}.

For any field $K$ containing $\mathbb C$, the group $G(K)$ of the points of $G$ in $K$ 
is generated by the following subgroups:
\begin{itemize}
\item the fundamental torus $T(K)$ where $T = {\rm Spec}(\mathbb Z[X]) $, hence $T(K)$ is 
isomorphic to the group $(K^*)^n=(K^*)\otimes_{\mathbb Z} Y$ and the character (resp.
cocharacter)  group of $T$ is $X$ (resp. $Y$).
\item  the root subgroups $U_{\alpha}(K)$ for $\alpha \in \Phi$, each isomorphic to
the additive group $(K,+)$ by an isomorphism (of algebraic groups) $x_{\alpha}$. 
\end{itemize}

 Actually, we consider an isomorphism $x_\alpha: K\simeq \mathfrak g_\alpha\otimes K\rightarrow
U_{\alpha}(K)$ where the additive group $\mathfrak g_\alpha\otimes K$ is identified with $K$ by the
choice of a Chevalley generator $e_\alpha$ of the $1-$dimensional complex space $\mathfrak g_\alpha$.

 Let $M$ be an $\mathfrak h-$diagonalizable $\mathfrak g-$module with weights in $X$, and where 
the action of each $\mathfrak g_\alpha$ is locally nilpotent ({\it e.g.} $\mathfrak g$ itself). Then $G(K)$  acts
on $M\otimes K$: the torus $T(K)$ acts via the character $\lambda$ on 
$M_\lambda\otimes K$ and the action of $x_\alpha(a)$ for $a\in\mathfrak g_\alpha\otimes K$ is the 
exponential of the action of $a$.

\subsubsection{About the freedom condition $F(\mathcal S)$}
\label{ssse:Freedom}

This condition is used in \ref{ssse:Titscone} to show that $C^v_f$ is nonempty, and in \ref{sse:NewLS} below for the existence of $\rho_{\Phi^+}$. If it fails, then $W^v$ as defined in \ref{ssse:WeylGroup} could be smaller than wanted (finite), and the roots of $\mathfrak g_{\mathcal S}$ could not
be defined by the adjoint action  of $\mathfrak h$. But, actually, $F(\mathcal S)$ is not necessary 
to define~$\mathfrak g_{\mathcal S}$ or $G_{\mathcal S}$ \cite{Remy02}.

 In \cite{Kac90}, \cite{Kumar02} and \cite{Littelmann94}, the condition $F(\mathcal S)$ and a 
minimality condition for the rank of $X$:\quad $\rank(X)=\;\vert I\vert+$corank$(A)$ are
assumed. Kumar requires moreover a  simple-connectedness condition:
$$(SC)\qquad \sum_{i\in I}\;\mathbb Z\alpha_i^\vee \hbox{ is cotorsion-free in\ } Y.$$
For a root generating system $\mathcal S$, define  $\mathcal S_{sc}=(A,X_{sc},Y_{sc},(\alpha_i)_{i\in I},({\alpha}{_i^*})_{i\in I})$ by
$X_{sc}=X\oplus\mathbb Z^I$, $Y_{sc}$ is the dual of $X_{sc}$ and 
${\alpha}{_i^*}(x+(n_j)_{j\in I})=\alpha_i^\vee(x)+n_i$. The group $G_{\mathcal S}$ is
a quotient of $G_{{\mathcal S}_{sc}}$ by a subtorus of the torus $T_{{\mathcal S}_{sc}}$,
central in $G_{{\mathcal S}_{sc}}$.  

 Starting from $\mathcal S$ satisfying $F(\mathcal S)$, then $\mathcal S_{sc}$ satisfies $F(\mathcal S_{sc})$ 
and $(SC)$. The group $G_{{\mathcal S}_{sc}}$ is the direct product of a  torus and a 
group with the properties assumed by Kumar. So, there is no trouble in using the
results of \cite{Kumar02} for the groups we define.

\subsubsection{Some commuting relations}
\label{ssse:CommRela}

We present here some relations in the group $G$, for more details, see \cite[1.5 and 1.6]{Rousseau06}. 
 
 If $\{ \alpha , \beta \}$ is a prenilpotent pair of roots in $\Phi$ 
(hence $\alpha\neq-\beta$), one denotes by $]\alpha,\beta[$ the finite set of the roots
$\gamma = p\alpha + q\beta \in \Phi$ with $p$ and $q$ strictly positive integers and
 $[\alpha,\beta] = ]\alpha,\beta[ \cup \{\alpha,\beta\}$;
we choose any total order on this set. Then, the product map:  $\prod_{\gamma\in[\alpha,\beta]}
U_\gamma(K) 
\rightarrow G(K)$ is a bijection onto the group $U_{[\alpha,\beta]}(K)$ generated by 
these $U_\gamma(K)$; it is actually an isomorphism of algebraic varieties. The  commutator
group $[U_\alpha(K),U_\beta(K)]$ is contained in $U_{]\alpha,\beta[}(K)$. More precisely, for $u
, v \in K$, one has: \quad$[x_\alpha(u),x_\beta(v)] =
 \prod x_\gamma( C_{p,q}u^pv^q)$\quad where the product runs over the 
$\gamma = p\alpha + q\beta \in ]\alpha,\beta[$  (in the fixed order) and the $C_{p,q}$ are integers.

 The group $T(K)$ normalizes $U_\alpha(K)$: if $t\in T(K)$ and $u\in K$ one has 
$tx_\alpha(u)t^{-1} = x_\alpha(\alpha(t)u)$. The subgroup $G^{(\alpha)}$ of $G$ generated by
$U_\alpha$, $U_{-\alpha}$ and $T$ is,  up to its center, isomorphic to $PGL_2$. In
particular, for $\alpha\in\Phi$ and $u\in U_\alpha(K)$, the set $U_{-\alpha}(K)uU_{-\alpha}(K)$
contains  a unique element  $m(u)$  conjugating $U_\beta(K)$ to $U_{r_\alpha(\beta)}(K)$ for
each 
$\beta\in\Phi$. Moreover, for $v , v' \in K^*$:
$$
\begin{array}{l}
m(x_\alpha(v)) = m(x_{-\alpha}(v^{-1})) =
x_{-\alpha}(-v^{-1})x_\alpha(v)x_{-\alpha}(-v^{-1}) =\cdots\\
\hfill\cdots=\,\alpha^\vee{\, }(v)m(x_\alpha(1)) = 
m(x_\alpha(1))\alpha^\vee{\, }(v^{-1})\\
m(x_\alpha(v))x_\alpha(v')m(x_\alpha(v))^{-1} = x_{-\alpha}(v^{-2}v')\\
m(x_\alpha(v))^2 = \alpha^\vee{\, }(-1).\\
\end{array}
$$

If $N$ is the normalizer of $T$ in $G$, then $N(K)$ is the group generated by
$T(K)$ and the $m(u)$ for all $\alpha \in \Phi$ and all $u \in U_{\alpha}(K)$.

\begin{lemm}
\label{le:actionN}
(\cite{Rousseau06} or \cite{Remy02}) 
There exists an algebraic 
homomorphism $\nu^{v}$ from $N$ onto $W^{v}$ such that $\nu^{v}(m(u)) = r_{i}$ for
 $u \in U_{\pm\alpha_{i}}(K)$ and ${\rm Ker}(\nu^{v}) = T$. 
As $K$ is infinite, $N(K)$ is the normalizer of $T(K)$ in $G(K)$ and all
 maximal split subtori of $G(K)$ are conjugate of $T(K)$. 
\end{lemm}

 The conjugacy action of $N$ on $T$ is given by
 $\nu^v$ where $W^v$ acts on $T$ through its action on $X$ or $Y$.

\subsubsection{Borel subgroups}
\label{ssse:Borels}

 The subgroup $U^+(K)$ of $G(K)$ is generated by the groups $U_\alpha(K)$ for 
$\alpha\in\Phi^+$;  it is normalized by $T(K)$. We define the same way $U^-(K)$ and
$U(\Psi)(K)$ for any  subset $\Psi$ of $\Phi^+$ or $\Phi^-$.

 The groups $B(K)=B^+(K)=T(K).U^+(K)$, $B^-(K)=T(K).U^-(K)$ are the 
{\it standard (positive, negative) Borel subgroups} of $G(K)$.

 One has $U^+(K)\cap B^-(K)=U^-(K)\cap B^+(K)=\{1\}$; more generally, one has the following decompositions.

\textit{Bruhat decompositions}:
$$G(K)=U^+(K)N(K)U^+(K)=U^-(K)N(K)U^-(K)$$
 \noindent Moreover, the maps from $N(K)$ onto $U^\pm (K)\backslash G(K)/U^\pm (K)$  are
one to one.

\textit{Birkhoff decompositions}: 
$$G(K)=U^+(K)N(K)U^-(K)=U^-(K)N(K)U^+(K)$$
 \noindent Moreover, the maps from $N(K)$ onto $U^\pm (K)\backslash G(K)/U^\mp(K)$ are  one
to one.

\subsection{The affine apartment}
\label{sse:AffineApart}

The affine apartment $\mathbb A$ is $V$ considered as an affine space.

\subsubsection{Affine Weyl group and preorder relation}
\label{ssse:AffineWeyl}

The group $W^v$ acts $\mathbb Z-$linearly on $Y$, hence it acts 
$\mathbb R-$linearly on $\mathbb A=V$. One has also an action of $V$ by translations. 
Finally, one obtains an affine action on $\mathbb A$ of the semi-direct product 
 $W_{\mathbb R} = W^v\ltimes V$.

 For $\alpha\in\Phi$ and $k\in\mathbb Z$,  $M(\alpha,k) =
 \{ v\in \mathbb A \, \mid \, \alpha(v) + k = 0 \}$
 is the {\it wall} associated to $(\alpha,k)$, it is closed in 
$\mathbb A$. One has $M(\alpha,k) = M(-\alpha,-k)$.

 For $\alpha\in X\setminus\{0\}$ and $k\in\mathbb Z$,  we define
\quad $D(\alpha,k) = \{ v\in \mathbb A\,\mid \, \alpha(v)+k \geq 0 \}$,
\noindent  it is closed in $\mathbb A$. When $\alpha\in\Phi$, we call $D(\alpha,k)$ 
the {\it half-apartment} associated to $(\alpha,k)$ and the set $D^{\circ}(\alpha,k) =
D(\alpha,k)\setminus M(\alpha,k) =  V\setminus D(-\alpha,-k)$ is the {\it open-half-apartment}
associated to $(\alpha,k)$.
 
 The {\it  reflection} associated to the wall $M(\alpha,k)$ is 
$\; r_{\alpha,k} : \mathbb A \rightarrow \mathbb A$ given by the formula: 
$$r_{\alpha,k}(y) = r_\alpha(y) - k\alpha^\vee.$$
The group generated by the $r_{\alpha,k}$ is
 $W = W^v\ltimes Q^\vee \subset W_{\mathbb R}$.

 The subgroup of $W_{\mathbb R}$ of all elements stabilizing the set 
of walls is $W_P = W^v\ltimes P^\vee$. One defines also
$W_Y = W^v\ltimes Y$ and one has:
$$W \subset W_Y \subset W_P \subset W_{\mathbb R}=W^v\ltimes V.$$

\begin{defi} 
\label{de:apart}
The affine space $\mathbb A$ together with its walls 
and its Tits cone is the  
{\it affine apartment} of $G$ associated to $T$, its {\it (affine) Weyl group} 
$W=W(\mathbb A)$  is generated by the reflections with respect to the walls.
\end{defi}

 As the Tits cone $\mathcal T$ is convex, we can define a {\it preorder-relation}  on
$\mathbb A$ given by $x\leq y\;\Leftrightarrow\; y-x\in\mathcal T$. This is a genuine 
(= antisymmetric) order relation only when $\Phi$ generates $V$ and the Kac-Moody 
matrix $A$ has no factor of finite type.

For $x\in\mathbb A$, the set $\Delta_x$ of all roots $\alpha$ such that  $\alpha
(x)\in\mathbb Z$ is a closed  subsystem of roots of $\Delta$ in the sense of \cite{Bardy96} Section 5.1.
The associated Weyl group 
$W_x^{\rm min}$ is the subgroup of $W$ generated by all the reflections associated  to
the walls  containing $x$. It is isomorphic to its image $W_x^v$ in $W^v$ and is a
Coxeter group,  as shown in [loc. cit.; 5.1.12]. The canonical generators of $W_x^v$
are the $r_\alpha$  for $\alpha$ simple in $\Phi_x^+=\Delta_x\cap\Phi^+$; their number may be infinite.

 The point $x$ is {\it special} when $\Phi_x=\Delta_x\cap\Phi$ 
is equal to $\Phi$, {\it i.e.} when $W_x^v=~W^v$.

\subsubsection{Faces}
\label{ssse:Faces}

 The faces in $\mathbb A$ are associated to the above systems of walls 
and half-apartments. As in \cite{BruhatTits72}, they 
are no longer subsets of $\mathbb A$, but filters of subsets of $\mathbb A$.

\begin{defi}[\cite{BruhatTits72}, \cite{Rousseau04} or \cite{Rousseau06}]
\label{de:Filters}
A {\it filter} in a set $E$ is a nonempty set $F$ of nonempty subsets of $E$, 
such that, if $S,S'\in F$ then $S\cap S'\in F$ and, if $S'\supset S\in F$  then
$S'\in F$.  If $Z$ is a nonempty subset of $E$, the set $F(Z)$ of subsets of $E$
containing 
$Z$ is a filter (usually identified with $Z$). If $E\subset E'$, to any filter 
$F$ in $E$ is associated the filter $F_{E'}$ in $E'$ consisting of all subsets of 
$E'$ containing some $S$ in $F$; one usually makes no difference between $F$ and 
$F_{E'}$.
\end{defi}

 A filter $F$ is said to be contained in another filter $F'$: $F\subset F'$ 
(resp.  in a subset $Z$ in $E$: $F\subset Z$) if and only if any set in $F'$ 
(resp. if $Z$) is in $F$. The union of a family of filters in $E$ is the filter 
consisting of subsets which are in all the filters. Note that these definitions 
are opposite the natural ones.

 A group $\Gamma$ acting on $E$ fixes pointwise (resp. stabilizes) a filter $F$, 
if and only if  every $\gamma$ in $\Gamma$ fixes pointwise some $S\in F$ (resp. for all
$\gamma$ in $\Gamma$ and all~$S$ in $F$, $\gamma S\in F$).

 If $E$ is a topological space, the closure of a filter $F$ in $E$ is 
the filter $\overline F$ consisting of all subsets of $E$ containing the 
closure of a set in $F$. 

If $F$ is a subset of $E$ containing an element $x$ in its closure, 
the {\it germ} of $F$ in $x$ is the filter ${\rm germ}_x(F)$ consisting of all subsets 
of $E$ which are intersections of $F$ and neighbourhoods of $x$.

 If $E$ is a real affine space and $x\neq y\in E$, then the {\it segment-germ} 
$[x,y)$ is the germ of the segment $[x,y]$ in $x$.

 All the above definitions for filters are compatible with the corresponding 
definitions for subsets and the identification of a subset $Z$ with the filter~$F(Z)$.

 We say that a family $\mathcal F$ of filters {\it generates} 
a filter $\Omega$ if: a set $S$ is in $\Omega$ if and only if it is in some filter  $F\in\mathcal F$.  If $\mathcal B$ is a basis of the filter $\Omega$, then the (filters canonically
associated  to the) sets in $\mathcal B$ generate $\Omega$.

 The {\it enclosure} $cl(F)$ of a filter $F$ of subsets of $\mathbb A$  is the filter 
made of  the subsets of $\mathbb A$ containing any intersection of half-spaces
$D(\alpha,k)$  (for $\alpha\in\Delta$ and $k\in\mathbb Z$), which is in $F$. With this
definition, the enclosure of a subset $\Omega$ is the closed subset intersection of all $D(\alpha,k)$  (for $\alpha\in\Delta$ and $k\in\mathbb Z$) containing $\Omega$. For $\mathcal P$ a non empty subset of
$X\setminus\{0\}$,  we define also the {\it $\mathcal P-$enclosure} $cl_{\mathcal P}(F)$ by 
the same definition, just replacing $\Delta$ by $\mathcal P$.

\begin{defi} 
\label{de:Faces}
A {\it face} $F$ in the apartment $\mathbb A$ is associated
 to a point $x\in \mathbb A$ and a  vectorial face $F^v$ in $V$; 
it is called spherical according to the nature of $F^v$. More 
precisely,  a subset $S$ of $\mathbb A$ is an element of the face $F(x,F^v)$ if and only if
it contains  an intersection of half-spaces $D(\alpha,k)$ or $D^\circ(\alpha,k)$ (for
$\alpha\in\Delta$ and $k\in\mathbb Z$) which contains
 $\Omega\cap(x+F^v)$, where $\Omega$ is an open neighborhood of $x$ in $\mathbb A$. The 
enclosure of a face $F=F(x,F^v)$ is its closure: the closed-face $\overline F$; it 
is the  enclosure of the {\it local-face in} $x$, ${\rm germ}_x(x+F^v)$.
\end{defi}

 Actually, in the classical case where $\Phi$ is finite, this definition 
is still valid: $F(x,F^v)$ is a subset $Z$ of $\mathbb A$ (more precisely: is the 
filter of subsets containing a subset $Z$ of $\mathbb A$), and this subset $Z$ is 
a face in the sense of \cite[\S1]{BruhatTits72} or \cite[6.1]{Brown89}.


 Note that the union of the faces $F(x,F^v)$ is not always the filter of
neighborhoods  of $x$; it is contained in $(x+\mathcal T)\cup(x-\mathcal T)$ if $x$ is special.

\subsubsection{Chambers, panels...} 
\label{ssse:Chambers}

 There is an order on the faces: the assertions ``$F$ is a face of $F'$ '', 
``$F'$ covers $F$ '' and ``$F\leq F'$ '' are by definition equivalent to
$F\subset\overline{F'}$.

 Any point $x\in \mathbb A$ is contained in a unique face $F(x,V_0)$ which is minimal 
(but
 seldom spherical); $x$ is a {\it vertex} if and only if $F(x,V_0)=\{x\}$. 
When $\Phi$ generates V ({\it i.e.} ${\rm rk}(X)=\vert I\vert$ ), a 
special point is a vertex, but the converse is not true.

 The {\it dimension} of a face $F$ is the smallest dimension of an affine 
space generated by some $S\in F$. The (unique) such affine space $E$ of minimal 
dimension is the {\it support} of $F$. Any $S\in F$ contains a non empty open 
subset of $E$.

 A {\it chamber} (or alcove) is a maximal face, or, equivalently, a face  such
that all  its elements contain a nonempty open subset of $\mathbb A$ or a face of dimension 
$\rank(X)=\dim(\mathbb A)$.

 A {\it panel} is a spherical face maximal among faces which are not chambers, or, 
equivalently, a spherical face of dimension $n-1$. Its support is a wall.

 So, the set  of spherical faces of $\mathbb A$ completely determines the set 
$\mathcal H$ of walls.

  A {\it wall of a chamber} $C$ is the support $M$ of a panel $F$ covered by 
$C$. Two chambers are called {\it adjacent} (along $F$ or $M$) if they cover a
common panel ($F$ of support $M$). But there may exist a chamber covering no panel,
and hence having no wall. So,
$\mathbb A$ is far from being ``gallery-connected''.

\subsubsection{Sectors}
\label{ssse:Sectors}

 A {\it sector} in $\mathbb A$ is a $V-$translate $\mathfrak s=x+C^v$ of a vectorial chamber 
$C^v=\pm w.C^v_f$ ($w \in W^v$), $x$ is its {\it base point} and $C^v$ its  {\it
direction}.  Two sectors have the same direction if and only if they are conjugate
by $V-$translation,
 and if and only if their intersection contains another sector. 

 The {\it sector-germ} of a sector $\mathfrak s=x+C^v$ in $\mathbb A$ is the filter $\mathfrak S$ of 
subsets of~$\mathbb A$ consisting of the sets containing a $V-$translate of $\mathfrak s$, it is well 
determined by the direction $C^v$. So the set of 
translation classes of sectors in $\mathbb A$, the set of vectorial chambers in $V$ and 
 the set of sector-germs in $\mathbb A$ are in canonical bijection.

 The sector-germ associated to the positive (resp. negative) fundamental 
chamber $C^v_f$ (resp. $-C^v_f$) is called the positive (resp. negative) fundamental 
sector-germ and is denoted by $\mathfrak S_{+\infty}$ (resp. $\mathfrak S_{-\infty}$).

 A {\it sector-face} in $\mathbb A$ is a $V-$translate $\mathfrak f=x+F^v$ of a vectorial face
$F^v=\pm wF^v(J)$. The {\it sector-face-germ} of $\mathfrak f$ is the filter $\mathfrak F$ of 
subsets containing a translate $\mathfrak f'$ of $\mathfrak f$ by an element of $F^v$ ({\it i.e.} $\mathfrak
f'\subset \mathfrak f$). If $F^v$ is spherical, then $\mathfrak f$ and $\mathfrak F$ are also called
spherical. The sign of $\mathfrak f$ and $\mathfrak F$ is the sign of $F^v$.

\section{The hovel, definition}
\label{se:Hovel}

 To define something like an affine building associated to the Kac-Moody 
group $G$ and the apartment $\mathbb A$  of $T(K)$, we have to define the action of
$N(K)$ on $\mathbb A$ and the fixator $\hat{P}_x$ in $G(K)$ of a point $x$ in 
$\mathbb A$, {\it i.e.} the associated parahoric subgroup. This fixator $\hat{P}_x$ should 
contain the fixator $\hat{N}_x$ of $x$ in $N(K)$ and the groups $U_{\alpha,k}$
for $x\in D(\alpha,k)$.  When $x$ is $0$ (``origin'' of $\mathbb A$), the group 
$\hat{P}_x$ should be $G(\mathcal O)$, so that the orbit of $0$ in the  ``building''
is the affine grassmannian $\mathcal G=G(K)/G(\mathcal O)$ (see Example \ref{ex:G(O)} below).

 But, we have also to define and study parahoric subgroups associated to 
more  general points or faces in $\mathbb A$ and this will lead to difficulties.
Moreover,  the expected Bruhat decomposition for parahoric subgroups is
actually false in our case (see Remark \ref{re:Preorder}).  So, the ``building'' we can construct has bad
properties, therefore, we call it an {\it hovel}  (in french ``masure'').

Here, we give an overview of the present section. First, we describe the action of $N(K)$ on the apartment $\mathbb A$ (\ref{sse:ActionNK}). Then, given a filter of 
subsets $\Omega$ in $\mathbb A$, we define a subgroup $P^{\rm min}_\Omega$ of $G(K)$
(\ref{sse:FirstObjects}), in the same way as in \cite{BruhatTits72}. But, due to some bad
commutation relations in $G(K)$, we have to work in larger groups, the formal
completions of $G$. There are two ways of doing it which lead to two groups
$P^{pm}_\Omega$ and $P^{nm}_\Omega$ (still in $G(K)$) both containing $P^{\rm min}_\Omega$
(see \ref{sse:3Groups}). All these groups are defined by generators. The ideal situation is
when they coincide, but in general, we need another group to compare them. So,
we define a fourth group $\widetilde P_\Omega$ as the stabilizer of some
subalgebra and some submodule for the action of $G(K)$ on highest weight
representations (in \ref{sse:Tilde}); it contains all the previous ones. However,
$\widetilde P_\Omega$ is a bit too big to be the fixator of $\Omega$ for the action
of $G(K)$ on the hovel. We get the ``parahoric group'' $P_\Omega$ by assuming that
$\mathfrak g$ is symmetrizable and by taking the fixator of some subalgebra for the
action of $\widetilde P_\Omega$. Finally, the ``right candidate'', as a fixator of
$\Omega$, is the group $\hat{P}_\Omega$ obtained by adding to $P_\Omega$ the fixator
$\hat{N}_\Omega$ of $\Omega$ in $N(K)$ (for its action on $\mathbb A$), see \ref{sse:Parahoric}. 

\subsection{Action of $N(K)$ on $\mathbb A$}
\label{sse:ActionNK}

We suppose now the field $K$ endowed with a discrete 
valuation $\omega$, assumed normalized: $\omega(K^*)=\mathbb Z$. The ring of integers  is $\mathcal O$; we choose  a uniformizing parameter $\varpi$, so $\omega(\varpi)=1$, $\mathcal O^*=\mathcal
O\setminus\varpi\mathcal O$ and the residue field is $\kappa=\mathcal O/\varpi\mathcal O$. Moreover,
we assume that $\kappa$ contains $\mathbb C$ (so, if $K$ is complete for $\omega$, then 
$K=\kappa(\!(\varpi)\!)$ and $\mathcal O=\kappa[[\varpi]])$. For the definition of the hovel 
for a Kac-Moody group $G$ over any valuated field, one needs more knowledge about $G$,
it should appear in \cite{Rousseau08b}.

 For $\alpha \in \Phi , u \in U_\alpha(K)$, $u\not = 1$, we define:  $\varphi_\alpha(u)
= \omega(t)$,
 if $u = x_\alpha(t)$ with $t \in K$. For all $k \in \mathbb R\cup\{+\infty\}$, 
the set $U_{\alpha,k} = \varphi_\alpha^{-1}([k,+\infty])$ is a subgroup of 
$U_\alpha(K)$ and $U_{\alpha,\infty} = \{1\}$. See also \cite[2.2]{Rousseau06}.  

\bigskip
The group $T(K)$ acts on $\mathbb A$ by translations:
\quad if $t \in T(K)\ ,\ \nu(t)$ is the element in $V$ such that 
$\chi(\nu(t)) = -\omega(\chi(t)) \; , \; \forall \chi \in X$. This action is $W^v-$equivariant.

The following lemma is a trivial consequence of the corresponding result 2.9.2 in 
\cite{Rousseau06}.

\begin{lemm} 
\label{le:ActionNK}
There exists an action $\nu$ of $N(K)$ on $\mathbb A$ 
which induces the preceding one on $T(K)$ and such that for $n\in N(K)$, 
$\nu(n)$ is an affine map with associated linear map $\nu^v(n)$.
\end{lemm}

\begin{remas}
\label{re:ActionNK}
\begin{enumerate}
\item[1)] The image of $N(K)$ in $Aut(\mathbb A)$ is $\nu(N) = W_Y$. The kernel 
$H={\rm Ker}(\nu)\subset T(K)$ is $H=\mathcal O^*\otimes Y=T(\mathcal O)$. 
\item[2)] By construction $\nu(N(\mathbb C))$ fixes $0$, the origin of $\mathbb A$, so, 
$\nu(m(x_\alpha(1)))$ is the reflection
 $r_\alpha = r_{\alpha,0}$ with respect to the wall $M(\alpha,0)$. Moreover, 
$m(x_\alpha(u)) = 
\alpha^\vee(u)m(x_\alpha(1))$, hence the image  $\nu(m(x_\alpha(u)))$ is the
reflection $r_{\alpha,\omega(u)}$ with respect to the wall $M(\alpha,\omega(u))$, as  by
definition
 one has: $\alpha(\nu(\alpha^\vee(u))) = - 
\omega(\alpha(\alpha^\vee(u)))= - \omega(u^2) = -2\omega(u)$.
\end{enumerate}
\end{remas}

\subsection{First objects associated to $\Omega$ and the group $P^{\rm min}_\Omega$}
\label{sse:FirstObjects}

Let $\Omega$ be a filter of subsets in $\mathbb A$. For $\alpha\in\Delta$, let 
$f_\Omega(\alpha)= $ inf$\{k\in\mathbb Z \mid \Omega\subset D(\alpha,k)\}=$ 
inf$\{k\in\mathbb Z\mid \alpha(\Omega)+k\subset[0,+\infty)\}\in\mathbb Z\cup\{+\infty\}$; 
by this second equality, $f_\Omega$ is defined on $X$.
 The function $f_\Omega$ is concave\cite{BruhatTits72}: $\forall \alpha, \beta \in X $, 
$f_\Omega(\alpha+\beta)\leq f_\Omega(\alpha)+ f_\Omega(\beta)$ and $f_\Omega(0)=0$; in particular
$f_\Omega(\alpha)+f_\Omega(-\alpha)\geq 0$. We say that $\Omega$
 is {\it narrow} (resp. {\it almost open}) if and only if 
$f_\Omega(\alpha)+f_\Omega(-\alpha)\in\{0,1\}$
 (resp. $\neq 0$), $\forall\alpha\in\Phi$. The filter $\Omega$ is almost open if  and
only if it is  not contained in any wall, this is true for a chamber. A point
or a face is narrow.  Actually, in the classical case, $\Omega$ is narrow if and
only if it is included in  the closure of a chamber.

We define $U_\Omega$ as the subgroup of $G(K)$ generated by the groups 
$U_{\alpha,\Omega}=U_{\alpha,f_{\Omega}(\alpha)}$ for $\alpha\in\Phi$, and  $U^\pm _\Omega=U_\Omega\cap
U^\pm (K)$. For $\alpha\in\Phi$, $U^{(\alpha)}_\Omega$  ($\subset U_\Omega$) is generated
by 
$U_{\alpha,\Omega}$ and $U_{-\alpha,\Omega}$; $N^{(\alpha)}_\Omega=N(K)\cap U^{(\alpha)}_\Omega$. The 
group $N^{u}_\Omega$ ($\subset N(K)\cap U_\Omega$) is generated by all 
$N^{(\alpha)}_\Omega$ for $\alpha\in\Phi$.

All these groups are normalized by $H$. In particular, 
one can define the groups $N^{\rm min}_\Omega=H.N^{u}_\Omega$ and $P^{\rm min}_\Omega=H.U_\Omega$. 
These groups depend only on the enclosure of $\Omega$ (not on $\Omega$ itself).

\begin{lemm} 
\label{le:SL2}
Let $\Omega$ be a filter of subsets in $\mathbb A$ and $\alpha\in\Phi$ a root.
\begin{enumerate}
\item[1)] $U^{(\alpha)}_\Omega=U_{\alpha,\Omega}.U_{-\alpha,\Omega}.N^{(\alpha)}_\Omega
=U_{-\alpha,\Omega}.U_{\alpha,\Omega}.N^{(\alpha)}_\Omega$.

\item[2)] If $f_\Omega(\alpha)+f_\Omega(-\alpha)>0$, then $N^{(\alpha)}_\Omega\subset H$.
 If $f_\Omega(\alpha)=-f_\Omega(-\alpha)=k$, then $\nu(N^{(\alpha)}_\Omega)=r_{\alpha,k}$.

\item[3)] $N^{(\alpha)}_\Omega$ fixes $\Omega$ {\it i.e.} $\forall n\in N^{(\alpha)}_\Omega$,
 $\exists S\in \Omega$ pointwise fixed by  $\nu(n)$.
\end{enumerate}
\end{lemm}

\begin{enonce*}[remark]{Consequence} 
The group $W^{\rm min}_\Omega=N^{\rm min}_\Omega/H$ is isomorphic
 to its image $W^v_\Omega$ in $W^v$, it is generated by the reflections 
$r_{\alpha,k}$ for which $\Omega\subset M(\alpha,k)$ ($\alpha\in\Phi$, $k\in\mathbb Z$). The group
$N^{\rm min}_\Omega$ is included in the group $\hat{N}_\Omega$, fixator in $N(K)$ of 
$\Omega$ which normalizes $H$, $U_\Omega$ and $P_\Omega^{\rm min}$. The group 
$\hat{W}_\Omega=\hat{N}_\Omega/H$ is also  isomorphic to a subgroup of $W^v$.
\end{enonce*}

\begin{proof}
Parts 1) and 2) are proved by an easy computation in 
$SL_2$ or $PGL_2$; one can also refer to \cite[6.4.7]{BruhatTits72} where a more 
complicated result (non split case) is proved. Clearly, 3) is a consequence  of
2).
\end{proof}

We define  $\mathfrak g_\Omega=\mathfrak h_{\mathcal O}\bigoplus(\oplus_{\alpha\in\Delta}\;\mathfrak g_{\alpha,\Omega})$,  where $\mathfrak h_{\mathcal O}=\mathfrak h\otimes_{\mathbb C}\mathcal O$, $\mathfrak g_{\alpha,\Omega}=\mathfrak g_{\alpha,f_{\Omega}(\alpha)}$ and  (in general) $\mathfrak g_{\alpha,k}=\mathfrak g_{\alpha}\otimes_{\mathbb C}\{t\in K\mid\omega(t)\geq k\}$. This is a  sub-$\mathcal O$-Lie-algebra of $\mathfrak g_K=\mathfrak g\otimes_{\mathbb C}K$.

The Lie algebra $\mathfrak g_\Omega$ depends only on the enclosure of $\Omega$ (not 
on $\Omega$ itself).  This is also true for the algebras and groups defined above in the Consequence of Lemma~\ref{le:SL2} (except for $\hat{N}_\Omega$ and $\hat{W}_\Omega$) and  below in Sections~\ref{sse:MaxKM} and \ref{sse:3Groups}.

If $\Omega$ is bounded, then  $\mathfrak g_\Omega$ is a lattice in  $\mathfrak g_K$. 

Let $M$ be a $\mathfrak g$-module of highest weight (resp. lowest weight) 
$\Lambda\in X$, then $M$  is the sum of its weight spaces: $M=\oplus_{\lambda\in
X}\;M_\lambda$. We define $M_\Omega=\oplus_{\lambda\in X}\;M_{\lambda,\Omega}$, where 
$M_{\lambda,\Omega}=M_{\lambda,f_{\Omega}(\lambda)}$ and  (in general)
$M_{\alpha,k}=M_\alpha\otimes_{\mathbb C}\{t\in K\mid\omega(t)\geq k\}$. This is a  sub$-\mathfrak g_\Omega-$module of $M\otimes K$, and a lattice when $\Omega$ is bounded. 

 If the module is integrable, then $\Lambda\in X^+$ (resp. $\Lambda\in X^-$) and 
$G(K)$ acts on $M\otimes K$. As we are in equal characteristic $0$, it is 
clear that $U_\Omega$  stabilizes~$M_\Omega$.

\subsection{Maximal Kac-Moody groups}
\label{sse:MaxKM}

\begin{enumerate}
\item[1)] The {\it positively-maximal Kac-Moody algebra} associated to $\mathfrak g$ is 
the Lie  algebra $\hat{\mathfrak g}^p= (\oplus_{\alpha\in\Delta^-}
\;\;\mathfrak g_{\alpha})\oplus\mathfrak h\oplus \hat{\mathfrak n}^+$  where $\hat{\mathfrak n}^+=\prod_{\alpha\in\Delta^+}\;\;\mathfrak g_{\alpha}$ is the completion of 
${\mathfrak n}^+=\oplus_{\alpha\in\Delta^+}\;\;\mathfrak g_{\alpha}$ \cite{Kumar02}.

\item[2)] The {\it positively-maximal Kac-Moody group} $G^{pmax}$ is defined 
in \cite{Kumar02}  (under the name $\mathcal G$); it contains $G$ as a subgroup ($G$ is
denoted by ${\mathcal G}^{\rm min}$ by  Kumar). For any closed subset $\Psi$ of $\Delta^+$,
$G^{pmax}$ contains the pro-unipotent  pro-group $U^{\max}(\Psi)$ with Lie
algebra $\hat{\mathfrak n}(\Psi)=\prod_{\alpha\in\Psi}\;\mathfrak g_{\alpha}$; 
{\it i.e.} $U^{\max}(\Psi)(K)=\prod_{\alpha\in\Psi}\;U_{\alpha}(K)$ where $U_{\alpha}(K)$ is 
isomorphic,  via an isomorphism $x_\alpha$, to $\mathfrak g_{\alpha}\otimes K$ (already
defined when $\alpha$ is real).

One has the {\it Bruhat decomposition}:
$$G^{pmax}(K)=\coprod_{n\in N(K)}\;U^{\max}(\Delta^+)(K)nU^{\max}(\Delta^+)(K),$$ 
and the {\it Birkhoff decomposition}:
$$G^{pmax}(K)=\coprod_{n\in N(K)}\;U^-(K)nU^{\max}(\Delta^+)(K).$$
Moreover,  
$$U^-(K)\cap N(K)U^{\max}(\Delta^+)(K)=N(K)\cap
U^{\max}(\Delta^+)(K)=\{1\}.$$

\begin{enonce*}[remark]{N.B} In all the preceeding or following notations, a sign $^+$ may 
replace $(\Psi)$ when $\Psi=\Delta^+$.
\end{enonce*}

\item[3)] The following subalgebras or subgroups associated to a filter $\Omega$ 
are also defined: 
\begin{itemize}
\item $\hat{\mathfrak g}^p_\Omega= {\mathfrak n}^-_\Omega\oplus\mathfrak h_{\mathcal O}\oplus
\hat{\mathfrak n}^+_\Omega$,  where ${\mathfrak n}^-_\Omega=\oplus_{\alpha\in\Delta^-}\;\;\mathfrak 
g_{\alpha,\Omega}$ and 
$\hat{\mathfrak n}_\Omega(\Psi)=\cup_{S\in\Omega}\;( \prod_{\alpha\in\Psi}\;\mathfrak g_{\alpha,S})$;
\item $U^{\max}_\Omega(\Psi)=\cup_{S\in\Omega}\;(\prod_{\alpha\in\Psi}\;U_{\alpha,S})$, where 
$U_{\alpha,S}=U_{\alpha,f_S(\alpha)}$ is $x_\alpha(\mathfrak g_{\alpha,S})$; as we are in equal characteristic
zero, the Campbell-Hausdorf  formula proves that this is a subgroup of
$U^{\max}(\Psi)(K)$;
\item $U^{pm}_\Omega(\Psi)=G(K)\cap U^{\max}_\Omega(\Psi)$, actually 
$U^{pm}_\Omega(\Psi)=U^+(K)\cap U^{\max}_\Omega(\Psi)$ because by  \cite[7.4.3]{Kumar02},
$U^+(K)=G(K)\cap U^{\max}(\Delta^+)(K)$. We have 
$U^{pm}_\Omega(\Psi)=\cup_{S\in\Omega}\;U^{pm}_S(\Psi)$ and 
$U^{pm}_{\Omega\cup\Omega'}(\Psi)=U^{pm}_{\Omega}(\Psi)\cap U^{pm}_{\Omega'}(\Psi)$.
\end{itemize}

\item[4)] Let $\alpha$ be a simple root, then by \cite[6.1.2, 6.1.3]{Kumar02}, 
 $U^{\max+}(K)=U_\alpha(K)\ltimes U^{\max}(\Delta^+\setminus\{\alpha\})(K)$. Using the same proof, one can show that
$U^{\max+}_\Omega=U_{\alpha,\Omega}\ltimes U^{\max}_\Omega(\Delta^+\setminus\{\alpha\})$ and, 
intersecting with $G(K)$, one gets $U^{pm+}_\Omega=U_{\alpha,\Omega}\ltimes
U^{pm}_\Omega(\Delta^+\setminus\{\alpha\})$.

 The groups $U^{\max}_\Omega(\Delta^+\setminus\{\alpha\})$ and 
$U^{pm}_\Omega(\Delta^+\setminus\{\alpha\})$ above  are normalized by 
$H.U^{(\alpha)}_\Omega$ and $U^{\max}(\Delta^+\setminus\{\alpha\})(K)$ is 
normalized by 
$$G^{(\alpha)}(K)=\langle\, T(K),U_\alpha(K),U_{-\alpha}(K)\,\rangle.$$

\item[5)] One has also to consider the {\it negatively-maximal Kac-Moody 
algebra}  associated to $\mathfrak g$, 
 $\hat{\mathfrak g}^n= (\oplus_{\alpha\in\Delta^+}\;\;\mathfrak g_{\alpha})\oplus\mathfrak h\oplus 
(\prod_{\alpha\in\Delta^-}\;\;\mathfrak g_{\alpha})$ and the associated {\it negatively-maximal 
Kac-Moody group} $G^{nmax}$. More generally, one can change $p$ to $n$ and 
$\pm$ to $\mp$ in 1),2),3), and 4) above in order to obtain similar groups
(with similar properties) in the negative case.
\end{enumerate}

\subsection{The groups $P^{pm}_\Omega$ and $P^{nm}_\Omega$}
\label{sse:3Groups}

\begin{prop} 
\label{pr:3Groups}
Let $\Omega$ be a filter of subsets in $\mathbb A$. We have $3$ subgroups of $G(K)$ 
associated to $\Omega$  and independent of the choice of a set of positive roots in its
$W^v-$conjugacy class:
\begin{enumerate}
\item[1)] The group $U_\Omega$ (generated by all $U_{\alpha,\Omega}$) is equal to 
$U_\Omega=U_\Omega^-.U_\Omega^+.N_\Omega^u=U_\Omega^+.U_\Omega^-.N_\Omega^u$.
\item[2)] The group $U_\Omega^{pm}$ generated by the groups $U_\Omega$ and $U^{pm+}_\Omega$ is 
equal to $U_\Omega^{pm}=U^{pm+}_\Omega.U_\Omega^-.N_\Omega^u$.
\item[3)] Symmetrically, the group $U_\Omega^{nm}$ generated by $U_\Omega$ and  
$U^{nm-}_\Omega$ is equal to $U_\Omega^{nm}=U^{nm-}_\Omega.U_\Omega^+.N_\Omega^u$.
\item[4)] One has: 
\begin{itemize}
\item[i)] $U_\Omega\cap N(K)=N_\Omega^u$
\item[ ii)] $U_\Omega^{pm}\cap N(K)=N_\Omega^u$
\item[iii)] $U_\Omega\cap (N(K).U^\pm (K))=N_\Omega^u.U^\pm _\Omega$ 
\item[iv)] $U_\Omega^{pm}\cap (N(K).U^+(K))=N_\Omega^u.U_\Omega^{pm+}$
\item[v)] $U_\Omega\cap U^\pm (K)=U^\pm _\Omega$
\item[vi)] $U_\Omega^{pm}\cap U^+(K)=U_\Omega^{pm+}$
\end{itemize}
and symmetrically for $U_\Omega^{nm}$.
\end{enumerate}
\end{prop}

\begin{remas} 
\label{rema:3.5}
The group $H = T(\mathcal O)$ normalizes also $U_\Omega^{pm}$ 
and $U_\Omega^{nm}$, moreover, $P^{\rm min}_\Omega$ is contained in $P^{pm}_\Omega =
H.U_\Omega^{pm}$ and in $P^{nm}_\Omega = H.U_\Omega^{nm}$.  The group $U^{++}_\Omega$
generated by the $U_{\alpha,\Omega}$ for 
$\alpha\in\Phi^+$ is included in $U^+_\Omega$, itself included in $U_\Omega^{pm+}$. The 
first  inclusion may be strict even for $\Omega$ reduced to a special point 
and $A$ of affine type. 
The equality $U^{+}_{\Omega}=U_\Omega^{pm+}$ is equivalent to $U_\Omega=U_\Omega^{pm}$, it 
may be false  for $\Omega$ large ({\it e.g.} a negative sector). 
The situation should be better for $\Omega$ narrow. Actually, we shall prove that 
$P_\Omega^{pm} = P_\Omega^{nm}$ when $\Omega$ is reduced to a special point or is a 
spherical face  (\ref{sse:Parahoric}). The problem is then to know if this group is generated
by its intersections
 with the torus and the (real) root groups.

In the classical case of reductive groups, one has 
$G=G^{pmax}=G^{nmax}$ and 
\begin{align*}
U^{++}_\Omega&=U^{+}_\Omega=U^{\max+}_\Omega=U^{pm+}_\Omega,\\ 
U^{--}_\Omega&=U^{-}_\Omega=U^{\max-}_\Omega=U^{nm-}_\Omega;
\end{align*} 
moreover 
$U_\Omega\;(=U^{pm}_\Omega=U^{nm}_\Omega)$ is  the same as the group defined in
\cite[6.4.2, 6.4.9]{BruhatTits72}. The group 
$P^{\rm min}_\Omega$ is called $P_\Omega$ by Bruhat and Tits.
\end{remas}

\begin{proof}{after \cite[6.4.9]{BruhatTits72}}
\begin{itemize}
\item[a)] Let ${\mathcal U}=U_\Omega^{pm}(\Delta^+).U_\Omega^{nm}(\Delta^-).N_\Omega^{u}\subset G(K)$. 
By \ref{sse:MaxKM}.4) and Lemma \ref{le:SL2}, for $\alpha$ simple, one has:
\begin{align*}
\mathcal U & =  U_\Omega^{pm}(\Delta^+\setminus\{\alpha\}).U_\Omega^{nm}(\Delta^-\setminus
\{-\alpha\}). U_{\alpha,\Omega}.U_{-\alpha,\Omega}.N_\Omega^{u}\\
 & =  U_\Omega^{pm}(\Delta^+\setminus\{\alpha\}).U_\Omega^{nm}(\Delta^-\setminus\{-\alpha\}).
U_{-\alpha,\Omega}.U_{\alpha,\Omega}.N_\Omega^{u}\\
&=U_\Omega^{pm}(\Delta^+\setminus\{\alpha\}).U_{-\alpha,\Omega}.U_\Omega^{nm}(\Delta^-\setminus \{-\alpha\}). U_{\alpha,\Omega}.N_\Omega^{u}\\
&=U_\Omega^{pm}(r_\alpha(\Delta^+)).U_\Omega^{nm}(r_\alpha(\Delta^-)).N_\Omega^{u}.
\end{align*}
So ${\mathcal U}$ does not change when $\Delta^+$ is changed by the Weyl group $W^v$.
 
\item[b)] Hence ${\mathcal U}$  is stable by left multiplication by $U_\Omega^{pm+}$ and all 
$U_{\alpha,\Omega}$  for $\alpha\in\Phi$. Moreover, it contains these subgroups, so ${\mathcal U}
\supset U_\Omega^{pm}\supset U_\Omega$.

\item[c)] In $G^{pmax}(K)$, let us prove that  ${\mathcal U}\cap
U^{\max+}(K)=U_\Omega^{pm+}$: if $xyz\in U^{\max+}(K)$ with $x\in U_\Omega^{pm+}$, 
$y\in U_\Omega^{nm-}$ and $z\in N_\Omega^{u}$,  then $yz\in U^{\max+}(K)$ and by the
Birkhoff decomposition (\ref{sse:MaxKM}.2) one has $y=z=1$.

\item[d)] So $\quad U_\Omega^{pm}\cap U^{+}(K)=U_\Omega^{pm}\cap U^{max+}(K)=
U_\Omega^{pm+}$. 

The group 
$U_\Omega(\Delta^+\setminus\{\alpha\}):=U_\Omega\cap U_\Omega^{pm}(\Delta^+\setminus\{\alpha\})
=U_\Omega^{+}\cap U_\Omega^{pm}(\Delta^+\setminus\{\alpha\})$ is normalized by 
$U_{\alpha,\Omega}$ and $U_{-\alpha,\Omega}$. By \ref{sse:MaxKM}.4),  $U^{+}_\Omega=U_{\alpha,\Omega}\ltimes
U_\Omega(\Delta^+\setminus\{\alpha\})$ and  symmetrically for $U^{-}_\Omega$ .

\item[e)] Now we are able to argue as in a), b) above with a new ${\mathcal U}$, where 
$U_\Omega^{pm}(\Delta^+)$ is changed to $U_\Omega^{+}$ and/or  $U_\Omega^{nm}(\Delta^-)$ to
$U_\Omega^{-}$.  This proves 1), 2) and 3).

\item[f)] Concerning 4), v) holds by definition, and vi) was proved in d). We 
prove now iv)  and ii); iii) and i) are similar. Let $n\in N(K)$ and $v\in
U^+(K)$ be such that $nv\in U_\Omega^{pm}$. There exist  $n'\in N_\Omega^{u}$, $u'\in
U_\Omega^{-}$ and $v'\in U_\Omega^{pm+}$ such that $nv=n'u'v'$. Now 
$n'^{-1}n=u'v'v^{-1}$ and by the  Birkhoff decomposition $n=n'\in N_\Omega^{u}$,
$v=u'v'$, so, $u'=1$ and $v=v'\in U_\Omega^{pm+}$.
\end{itemize}
\end{proof}

\subsection{Iwasawa decomposition}
\label{sse:Iwasawa}

\begin{prop}
\label{pr:Iwasawa}
Suppose $\Omega$ narrow, then 
$$
G(K)=U^{+}(K).N(K).U_\Omega.
$$
Suppose, moreover, $\Omega$ almost open. Then the natural map from
$W_Y=N(K)/H$ onto $U^{+}(K)\backslash N(K)/U_\Omega$ is one to one.
\end{prop}

\begin{remas} 
\label{re:Iwasawa}
We also have $G(K)=U^{-}(K).N(K).U_\Omega$ and similarly with the maximal groups $G^{pmax}(K) = U^{pmax+}(K).N(K).U_\Omega$ and $G^{nmax}(K)=U^{nmax-}(K).N(K).U_\Omega$.

As a consequence, when $\Omega$ is narrow, every subgroup $P$ of 
$G(K)$  containing $U_\Omega$ may be written $P=(P\cap (U^{+}(K).N(K))).U_\Omega$.
If, moreover, 
$P\cap (U^{+}(K).N(K))=U^+_P.N_P$ with $U^+_P=P\cap U^{+}(K)$ and  $N_P=P\cap
N(K)$  normalizing $U_\Omega$, then $P=U^+_P.N_P.U^-_\Omega$. We shall use this 
to (almost) identify $U_\Omega^{pm}$ and $U_\Omega^{nm}$ (see Section \ref{sse:Parahoric}).
\end{remas}

The idea of the proof of the Iwasawa decomposition goes back to Steinberg. We follow  \cite[7.3.1]{BruhatTits72},  see also \cite[3.7]{KacPeterson85} and \cite[1.6]{Garland95}. We first need a lemma.

\begin{lemm} 
\label{le:Iwasawa}
Let $\alpha$ be in $\Phi$, then 
$Z_\alpha:=U_\alpha(K).\{1,r_\alpha\}.T(K).U_\Omega^{(\alpha)}$ contains 
$G^{(\alpha)}(K)$.
\end{lemm}

\begin{proof}
By the Bruhat decomposition, 
$$G^{(\alpha)}(K)\subset 
U_\alpha(K).\{1,r_\alpha\}.T(K).U_\alpha(K).$$ 
So it suffices to prove that, for 
$m_\alpha\in N(K)$ such that $\nu^v(m_\alpha)=r_\alpha$ and $u\in U_\alpha(K)$,  $m_\alpha
u\in Z_\alpha$.  If $\varphi_\alpha(u)\geq f_\Omega(\alpha)$, it's clear: $u\in U_{\alpha,\Omega}$.
Otherwise 
$\varphi_\alpha(u)\leq f_\Omega(\alpha)-1\leq -f_\Omega(-\alpha)$ and $u=v'mv"$ with 
$\nu^v(m)=r_\alpha$, 
$v',v''\in U_{-\alpha,-\varphi_\alpha(u)}\subset U_{-\alpha,\Omega}$. So  $m_\alpha u=m_\alpha
v'mv"\in U_\alpha(K).T(K).U_{-\alpha,\Omega}\subset Z_\alpha$, and the lemma is proved.

\begin{enonce*}[remark]{Proof of Proposition \ref{pr:Iwasawa}} 
The set $Z=U^+(K).N(K).U_\Omega$ is stable by left multiplication by 
$U^+(K)$  and $T(K)$. It remains to prove that it is stable by left
multiplication by 
$U_{-\alpha}(K)$ for $\alpha$ a simple root. Let $U(\Phi^+\setminus\{\alpha\})(K)=G(K)
\cap U^{\max}(\Delta^+\setminus\{\alpha\})(K)\subset U^+(K)$, using the Lemma \ref{le:SL2} and discussion in Section \ref{sse:MaxKM}.4), one gets:
\begin{align*}
U_{-\alpha}(K)Z&=U_{-\alpha}(K).U(\Phi^+\setminus\{\alpha\})(K).U_{\alpha}(K).N(K).U_\Omega\\
&\subset U(\Phi^+\setminus\{\alpha\})(K).G^{(\alpha)}(K).N(K).U_\Omega\\
&\subset U(\Phi^+\setminus\{\alpha\})(K).U_\alpha(K).\{1,r_\alpha\}.
T(K).U_\Omega^{(\alpha)}.N(K).U_\Omega\\
&\subset U^+(K).T(K).U_{-\alpha}(K).N(K).U_\Omega\\
&\qquad\cup\;\;U^+(K).T(K).r_\alpha.U_{-\alpha}(K).U_{\alpha}(K).N(K).U_\Omega\\
&\subset U^+(K).T(K).U_{-\alpha}(K).N(K).U_\Omega.
\end{align*}
\end{enonce*}
 It remains to show that $U_{-\alpha}(K).N(K)\subset Z$. But 
$un=n.n^{-1}un\in nU_{\beta}(K)
\subset nU_{-\beta}(K).\{1,r_\beta\}.T(K).U_\Omega^{(\beta)}$ with 
$\beta=-\nu^v(n^{-1})\alpha$. So 
$un\in U_\alpha(K).n.\{1,r_\beta\}.T(K).U_\Omega^{(\beta)}\subset U_{\alpha}(K).N(K).U_\Omega$.

With obvious notation, suppose $n'\in U^+(K)nU_\Omega$. Then, by Lemma~\ref{le:SL2} and Proposition~\ref{pr:3Groups} one has: $n'^{-1}n\in U_\Omega n^{-1}U^+(K)n$. But,
$n^{-1}U^+(K)n = U(n^{-1}\Phi^+)$. Further, \noindent
$U_\Omega.U(n^{-1}\Phi^+)(K)\subset H.U_\Omega(n^{-1}\Phi^-).U(n^{-1}\Phi^+)(K)
\subset U(n^{-1}\Phi^-)(K).H.U(n^{-1}\Phi^+)(K)$. 

Finally, by the Birkhoff decomposition, 
$n'^{-1}n\in H$. 

\end{proof}

\subsection{The group $\widetilde P_\Omega$}
\label{sse:Tilde}

In this section $\Omega$ is asked to be a nonempty set.

Clearly, $U_\Omega^{\max+}$ stabilizes  $\hat{\mathfrak g}^p_\Omega$ and $G(K)$ 
stabilizes $\mathfrak g_K$; so $U_\Omega^{pm+}=G(K)\cap U_\Omega^{\max+}$ stabilizes 
$\mathfrak g_\Omega=\hat{\mathfrak g}^p_\Omega\cap\mathfrak g_K$. Finally, $U_\Omega^{pm}$ and also 
$U_\Omega^{nm}$ (or $H$) stabilize $\mathfrak g_\Omega$. If $M$ is a highest weight integrable $\mathfrak g-$module, then $U_\Omega^{\max+}$ 
stabilizes $M_\Omega$. The group $U_\Omega^{\max-}$ stabilizes $\hat{M}_\Omega=
\prod_{\lambda\in X}\,M_{\lambda,\Omega}$ and $G(K)$ stabilizes $M\otimes K$. Finally, 
$U_\Omega^{pm}$ and also $U_\Omega^{nm}$ (or $H$) stabilize $M_\Omega$ for every 
highest (or lowest)  weight integrable $\mathfrak g-$module $M$.

\begin{defi} 
\label{de:Tilde}
The group $\widetilde P_\Omega$ is the subgroup of all
elements in $G(K)$ stabilizing $\mathfrak g_\Omega$ and $M_\Omega$ for every highest (or 
lowest) weight integrable $\mathfrak g-$module $M$.
\end{defi}

 Hence, $\widetilde P_\Omega$ contains $U_\Omega^{pm}$, $U_\Omega^{nm}$, $U_\Omega$ 
and $H$. When $\Omega$ is narrow,  we have $\widetilde P_\Omega=({\widetilde P_\Omega}\cap U^+(K).N(K)).U_\Omega= ({\widetilde P_\Omega}\cap U^-(K).N(K)).U_\Omega$.

\begin{lemm} 
\label{le:Tilde}
Let $\widetilde N_\Omega = \widetilde P_\Omega\cap N(K)$, 
then $\widetilde P_\Omega\cap (U^+(K).N(K))=U_\Omega^{pm+}.\widetilde N_\Omega$ and 
$\widetilde P_\Omega\cap (U^-(K).N(K))=U_\Omega^{nm-}.\widetilde N_\Omega$.
 Moreover, $\widetilde N_\Omega$ normalizes $U_\Omega$ and is the stabilizer (in 
$N(K)$ for the action $\nu$ on $\mathbb A$) of the $\mathcal P-$enclosure $cl_{\mathcal P}(\Omega)$ 
of $\Omega$, where $\mathcal P\subset X$ is the
 union of $\Delta$ and the set of all weights of $\mathfrak h$ in all the modules $M$ above. 
\end{lemm}

\begin{proof}
\begin{itemize}
\item[a)] Let $n\in N(K)$ and $u\in U^+(K)$ be such that 
$un\in \widetilde P_\Omega$ and $w=\,\nu^v(n)$. For $\mathcal M=M_\Omega$ or $\mathfrak g_\Omega$, 
$g\in \widetilde P_\Omega$ and $\mu,\mu'\in \mathfrak h^*$, we define  $_{\mu'}\vert
g\vert_\mu$  as the restriction of $g$ to $\mathcal M_\mu$ followed by the projection
onto $\mathcal M_{\mu'}$  (with kernel $\oplus_{\mu''\neq\mu'}\; \mathcal M_{\mu''})$. Now for all  $\mu$, 
$_{w\mu}\vert un\vert_\mu=\;_{w\mu}\vert n\vert_\mu$ and 
$n=\oplus_\mu\;_{w\mu}\vert n\vert_\mu$ 
(in an obvious sense); so $n\in\widetilde N_\Omega$. We have 
$\widetilde P_\Omega\cap (U^+(K).N(K))=(\widetilde P_\Omega\cap U^+(K)). \widetilde
N_\Omega$ and  it remains to determine $\widetilde N_\Omega$ and $\widetilde
P_\Omega\cap U^+(K)$ (or 
$\widetilde P_\Omega\cap U^-(K)$).

\item[b)]$\widetilde P_\Omega\cap U^+(K)=U_\Omega^{pm+}$: the inclusion $\supset$ 
is already  proved in the discussion before Definition~\ref{de:Tilde}. So, consider
$u=\prod_{\alpha\in\Delta^+}\;u_\alpha\in  U^{\max+}(K)$ such that $u$ stabilizes ${\mathfrak 
g}_\Omega$ (the order  on the $u_\alpha\in U_\alpha(K)$ is chosen such that the height
of $\alpha$ is  increasing from right to left). We shall prove by induction that
each $u_\alpha$ is in 
$U_{\alpha,\Omega}$. We may suppose $u_{\alpha'}\in U_{\alpha',\Omega}$ for $u_{\alpha'}$ on the
 right of $u_\alpha$; moreover, as $U_{\alpha',\Omega}$ stabilizes ${\mathfrak g}_\Omega$, 
we may suppose these $u_{\alpha'}$ equal to $1$. So 
$u=(\prod_{\beta\neq\alpha}\;u_\beta).u_\alpha$
 where the $u_\beta$ are in $U_\beta(K)$ and $ht(\beta)\geq ht(\alpha)$. But 
$_{\alpha}\vert u\vert_0=\,_{\alpha}\vert u_\alpha\vert_0$ sends $\mathfrak h_{\mathcal O}$ into 
$\mathfrak g_{\alpha,\Omega}$, so $u_\alpha\in U_{\alpha,\Omega}$.
 Now if $u\in\widetilde P_\Omega\cap U^+(K)$, it is in $U^{\max+}(K)\cap G(K)$
 and stabilizes ${\mathfrak g}_\Omega$; by the above argument, 
$u\in U^{\max+}_\Omega\cap G(K)=U_\Omega^{pm+}$.

\item[c)] Let $n=n_0t$, $n_0\in N(\mathbb C)$, $\nu^v(n)=w$ and $t\in T(K)$, then 
$nM_{\lambda,k}=M_{w\lambda,k+\omega(\lambda(t))}$. Consider now the action on $\mathbb A$:  
$nD(\lambda,k)=n_0D(\lambda,k+\omega(\lambda(t)))=D(w\lambda,k+\omega(\lambda(t)))$. But, $\mathfrak g_\Omega$ is 
generated by $\mathfrak h_{\mathcal O}$ and the $\mathfrak g_{\alpha,\Omega}$ for $\alpha\in\Delta$, so $n$ 
is in $\widetilde P_\Omega$ if and only if, for all $\lambda\in\mathcal P$ 
$f_\Omega(\lambda)+\omega(\lambda(t))=f_\Omega(w\lambda)$ if and only if, for all $\lambda\in\mathcal P$,
$nD(\lambda,f_\Omega(\lambda))=D(w\lambda,f_\Omega(w\lambda))$. This is equivalent to 
the fact that $n$ stabilizes the set  $cl_{\mathcal P}(\Omega)$. Moreover, as $\widetilde N_\Omega$ stabilizes $\mathfrak g_\Omega$,  it
normalizes $U_\Omega$.
\end{itemize}
\end{proof}

We know that  $N^{\rm min}_\Omega=H.N^u_\Omega\subset\widetilde 
N_\Omega$.  So, to determine $\widetilde N_\Omega$, we only have to determine the
subgroup $\widetilde W_\Omega=\widetilde N_\Omega/H=\nu(\widetilde N_\Omega)$ of 
$\nu( N(K))$; it contains $W_\Omega^{\rm min}=\nu(N_\Omega^{\rm min})$.

\begin{exems}
\label{ex:Bounded}
\begin{enumerate}
\item[1)] {\it Let us now assume that $\Omega$ is bounded.} As $\mathcal P\supset\Delta\cup
X^+\cup X^-$, it is easy to prove that each $\chi\in X$ is a positive linear
combination of some $\lambda\in\mathcal P$. 
Hence, the intersection $cl_{\mathcal P}(\Omega)$ of all $D(\lambda,f_\Omega(\lambda))$'s 
(for $\lambda\in\mathcal P$) is a nonempty convex compact set.
 But $\widetilde N_\Omega$ stabilizes $cl_{\mathcal P}(\Omega)$ and, as it acts 
affinely, it fixes a point $x_\Omega$ in $cl_{\mathcal P}(\Omega)$.
\item[2)] {\it Suppose now $\Omega$ narrow}, then 
$$\widetilde P_\Omega=
U_\Omega^{pm+}.\widetilde N_\Omega.U_\Omega=U_\Omega^{pm+}.U_\Omega.\widetilde N_\Omega
=U_\Omega^{pm+}.U_\Omega^-.\widetilde N_\Omega$$ and 
$$\widetilde P_\Omega=
U_\Omega^{nm-}.\widetilde N_\Omega.U_\Omega=U_\Omega^{nm-}.U_\Omega.\widetilde N_\Omega
=U_\Omega^{nm-}.U_\Omega^+.\widetilde N_\Omega.$$ 
In particular $\widetilde P_\Omega$, which 
contains always $P_\Omega^{pm}$ and $P_\Omega^{nm}$, is not much greater than them in
this case.
\end{enumerate}
\end{exems}

\subsection{The (parahoric) group $P_\Omega$ and the ``fixator'' $\hat P_\Omega$}
\label{sse:Parahoric}

{\sl From now on, we suppose $\mathfrak g$  symmetrizable.}

As $\mathfrak g_\Omega$ is generated by $\mathfrak h_{\mathcal O}$ and the  $\mathfrak 
g_{\alpha,\Omega}$ for $\alpha\in\Delta$, the derived algebra of $\mathfrak g_\Omega$ is $\mathfrak g'_\Omega= 
(\sum_{\alpha\in\Delta}\;[\mathfrak g_{\alpha,\Omega},\mathfrak g_{-\alpha,\Omega}])\oplus
(\oplus_{\alpha\in\Delta}\;\mathfrak g_{\alpha,\Omega})$. Consider the quotient algebra 
$\overline{\mathfrak g}_\Omega=\mathfrak g_\Omega/\varpi\mathfrak g_\Omega=(\mathfrak h\otimes\kappa)\oplus 
(\oplus_{\alpha\in\Delta}\;\mathfrak g_{\alpha,\Omega}/\varpi\mathfrak g_{\alpha,\Omega})$. As $\mathfrak g$
is  symmetrizable, $[\mathfrak g_{\alpha},\mathfrak g_{-\alpha}]=\mathbb C\alpha^\vee$ for all $\alpha\in\Delta$, so the derived 
algebra of $\overline{\mathfrak g}_\Omega$  is $\overline{\mathfrak g}'_\Omega= (\sum_{\alpha\in\Delta_\Omega}\;\kappa\alpha^\vee)\oplus
(\oplus_{\alpha\in\Delta}\;\mathfrak g_{\alpha}\otimes\kappa)$,
 where $\Delta_\Omega=\{\alpha\in\Delta\mid f_\Omega(\alpha)+f_\Omega(-\alpha)=0\}$ is the set of 
 $\alpha\in\Delta$ such that $\alpha(\Omega)$ is reduced to a point in $\mathbb Z$.

\medskip
As $\mathfrak g$ is  symmetrizable,  the orthogonal 
$(\overline{\mathfrak g}'_\Omega)^\perp$ of $\overline{\mathfrak g}'_\Omega$ in  $\overline{\mathfrak 
g}_\Omega$    is $\{x\in\mathfrak h\otimes\kappa\mid\alpha(x)=0,\forall\alpha\in\Delta_\Omega\}$.  
If $\Omega$ is a set, the action of $\widetilde P_\Omega$  (by inner automorphisms) 
is compatible with  the invariant bilinear form; so $\widetilde P_\Omega$
stabilizes $\mathfrak g_\Omega$, $\overline{\mathfrak g}_\Omega$, $\overline{\mathfrak g}'_\Omega$ and 
$(\overline{\mathfrak g}'_\Omega)^\perp$. Let $\widetilde P'_\Omega$ be the fixator of
$(\overline{\mathfrak g}'_\Omega)^\perp$ for this action of $\widetilde P_\Omega$.

\begin{defi} 
\label{de:Parahoric}
For a filter $\Omega$, $P_\Omega=\cup_{S\in\Omega}\;(\cap_{S'\subset S}\;\widetilde P'_{S'})$
 is a "parahoric" group associated to $\Omega$.
\end{defi}
 
An element $g$ of $U_\Omega$ (or $U_\Omega^{pm+}$, $U_\Omega^{nm-}$) is in some
$U_S$ (or $U_S^{pm+}$, $U_S^{nm-}$) for $S\in\Omega$, hence in $U_{S'}$ (or
$U_{S'}^{pm+}$, $U_{S'}^{nm-}$) for any $S'\subset S$. But $U_{S'}^{max+}$ and $H$
induce the identity on $\mathfrak h_{\mathcal O}$, so $\widetilde P'_{S'}$ contains
$H$, $U_{S'}^{pm+}$ and also $U_{S'}^{nm-}$. Moreover $U_{S'}$ is generated
by elements in $U_{S'}^{pm+}$ or $U_{S'}^{nm-}$. Finally $P_\Omega$ contains $U_\Omega$,
$U_\Omega^{pm+}$, $U_\Omega^{nm-}$ and $H$. 

The group $N_\Omega=P_\Omega\cap N(K)$ contains $N_\Omega^{\rm min}$ (and is often  
equal to it, as we shall see). The quotient group $W_\Omega=N_\Omega/H$ contains 
$W_\Omega^{\rm min}$ and is included  in $\widetilde W_\Omega:=\cup_{S\in\Omega}\;\widetilde
W_S$. Actually, $W_\Omega=\cup_{S\in\Omega}\;(\cap_{S'\subset S}\;\widetilde W'_{S'})
$ with $\widetilde W'_{S'}=(N(K)\cap \widetilde P'_{S'})/H$. If $S'$ is a non empty bounded set ,
then, by Example~\ref{ex:Bounded}, $\widetilde W_{S'}$ fixes a point $x_{S'}$ in
$cl_{\mathcal P}(S')\subset\{x\in\mathbb A\mid\alpha(x)=\alpha(S')\quad\forall\alpha\in\Delta_{S'}\}\subset
\{x\in\mathbb A\mid\alpha(x)=\alpha(S)\quad\forall\alpha\in\Delta_{S}\}$ 
(if $S'\subset S$); in particular, it
is isomorphic to its image in $W^v$. But, by definition, the image in $W^v$ of
$\widetilde W'_{S'}$ is the fixator (in the image of $\widetilde
W_{S'}$) of the direction of the affine space $\{x\in\mathbb A\mid\alpha(x)=\alpha(S')\quad\forall\alpha\in\Delta_{S'}\}\supset cl_{\mathcal P}(S')$.
Hence, by
Lemma~\ref{le:Tilde}, $\widetilde W'_{S'}$ is the fixator in $W_Y$ of
$\{x\in\mathbb A\mid\alpha(x)=\alpha(S')\quad\forall\alpha\in\Delta_{S'}\}$. It follows that $W_\Omega$ is always the fixator in $W_Y$ of
$\{x\in\mathbb A\mid\alpha(x)=\alpha(\Omega)\quad\forall\alpha\in\Delta_{\Omega}\}$. In particular $N_\Omega$
normalizes $U_\Omega$.

 When there exists $x,y\in\{x\in\mathbb A\mid\alpha(x)=\alpha(\Omega)\quad\forall\alpha\in\Delta_{\Omega}\}$
  such that $y-x$ is in the open-Tits-cone  (in particular when $\Omega$ is a
spherical face), it is known that $W_\Omega=W_\Omega^{\rm min}$.

 When $\Omega$ is narrow, $P_\Omega=U_\Omega^{pm+}.U_\Omega^{-}.N_\Omega
=U_\Omega^{nm-}.U_\Omega^{+}.N_\Omega=U_\Omega^{pm}.N_\Omega=U_\Omega^{nm}.N_\Omega$.

 In particular, when $\Omega$ is a spherical face (or a special point), 
$N_\Omega=N_\Omega^{\rm min}$ and $P_\Omega=P_\Omega^{pm}=P_\Omega^{nm}$ is called the 
{\it parahoric subgroup }associated to $\Omega$.

\begin{defi} 
\label{de:Fixator}
The "fixator" $\hat P_\Omega$ associated to $\Omega$ is the group generated by $P_\Omega$ and the fixator $\hat N_\Omega$ (in $N(K)$ for the action $\nu$) of  $\Omega$.
\end{defi}

Actually, $\hat N_\Omega$ is also the fixator of the {\it support} of $\Omega$: the 
smallest affine subspace of $\mathbb A$ generated by a set in $\Omega$. Clearly 
${\rm supp}(\Omega)\subset\{x\in\mathbb A\mid\alpha(x)=\alpha(\Omega)\quad\forall\alpha\in\Delta_{\Omega}\}$, 
so $\hat N_\Omega\supset N_\Omega$. As $\hat N_\Omega$ 
normalizes $P_\Omega$ (and all the groups previously defined), we have
 $\hat P_\Omega=P_\Omega.\hat N_\Omega$. Clearly,  $\hat P_\Omega\cap
N(K)=\hat N_\Omega$  and $\hat P_\Omega\supset U_\Omega^{pm+}$, $U_\Omega^{nm-}$.

 When $\Omega$ is narrow, $\hat P_\Omega=U_\Omega^{pm+}.U_\Omega^{-}.\hat N_\Omega
=U_\Omega^{nm-}.U_\Omega^{+}.\hat N_\Omega$.

 This group should be the fixator of $\Omega$ for the action of $G(K)$ 
on the ``ugly-building'' we shall build now. But this will be proved only for 
some $\Omega$, see \ref{sse:ExGood} below.

\begin{exems}
\label{ex:G(O)}
An explicit computation: Suppose $\Omega$ reduced to the special point $0$, 
the origin of $V=\mathbb A$ chosen as in Remark \ref{re:ActionNK}.2). Then $f_0(\alpha)=0,\forall\alpha$,  $\mathfrak 
g_0=\mathfrak g\otimes_{\mathbb C}\mathcal O$ and $M_{0}=M\otimes_{\mathbb C}\mathcal O$. Hence, the
definition of the ind-group structure of $G$ \cite[7.4.6 and 7.4.7]{Kumar02} tells
us that $\widetilde P_0\subset G(\mathcal O)$. Moreover, $Lie(G)=\mathfrak g$ and the
highest or lowest weight modules are defined by morphisms of ind-varieties
[loc. cit.; 7.4.E(6) and 7.4.13] so $\widetilde P_0=G(\mathcal O)$. Now
$\overline{\mathfrak g}_0=\mathfrak g\otimes_{\mathbb C}\kappa$ and (as $0$ is special)  $\hat
N_0=N_0^{\rm min}$, $(\overline{\mathfrak g}'_0)^\perp=\mathfrak  c\otimes_{\mathbb C}\kappa$ where  $\mathfrak 
c$ is the center of $\mathfrak g$; so
$G(\mathcal O)=\widetilde P_0=P_0=\hat P_0$\quad ($=G_0$ with the notation of~\ref{sse:Good}).

In the classical case of reductive groups, $W_\Omega$ is always equal to 
$W_\Omega^{\rm min}$. If $\Omega$ is narrow ({\it i.e.} included in a closed-face), 
$P_\Omega=P_\Omega^{\rm min}$ and $\hat P_\Omega$ are as defined by Bruhat and Tits 
(cf. Remark \ref{rema:3.5}). In particular, $\hat P_x$ is the same as in Bruhat-Tits 
and the  following definition gives the (pretty) Bruhat-Tits building.
\end{exems}

\subsection{The hovel and its apartments}

\begin{defi} 
\label{de:Hovel}
The {\it hovel} $\mathcal I=\mathcal I(G,K)$ of $G$ over $K$ is the 
quotient of the set $G(K)\times\mathbb A$ by the relation:

\quad $(g,x) \sim (h,y)\quad \Leftrightarrow\quad \exists n\in N$ such that
 $y = \nu(n)x$ and $g^{-1}hn \in \hat{P}_x$.
\end{defi}

One proves easily \cite[7.4.1]{BruhatTits72} that $\sim$ is an equivalence 
relation.  Moreover, $\hat{P}_x\cap N(K) = \hat{N}_x$. 
So, the map $x \mapsto cl(1,x)$ identifies $\mathbb A$ with its image $A_f=A(T,K)$, 
the {\it apartment of $T$ in $\mathcal I(G,K)$}. 

 The left action of $G(K)$ on $G(K)\times\mathbb A$ descends to an 
action on $\mathcal I$. The {\it apartments} of 
$\mathcal I$ are the $g.A_f$ for $g\in G(K)$. The action of $N(K)$ on $\mathbb A = A_f$ 
is through $\nu$; in  particular, $H$ fixes (pointwise) $A_f$. By construction,
the fixator of $x \in \mathbb A$ is $\hat{P}_x$ and, for $g \in G(K)$, one has 
$gx \in \mathbb A \Leftrightarrow g \in N(K)\hat{P}_x$.

 From the definition of the groups $\hat{P}_x$, it is clear that, for $\alpha
 \in \Phi$ and $u \in K$, $x_\alpha(u)$ fixes $D(\alpha,\omega(u))$. Hence,
 for $k\in\mathbb Z$, the group $H.U_{\alpha,k}$ fixes $D(\alpha,k)$.

\section{The hovel, first properties}
\label{se:HovelProperties}

 First, we define the notion of good fixator for a filter $\Omega$ of $\mathbb A$. 
It formalizes the fact that the fixator $G_\Omega$ of $\Omega$ for the action of
$G(K)$ on $\mathcal I$ has a nice decomposition and the fact that $G_\Omega$ acts
transitively on the apartments containing $\Omega$ (\ref{sse:Good}). Thanks to a technical
proposition (Proposition~\ref{pr:Good}), we can show, in particular, that faces, sectors,
sector-germs, walls and half-apartments in $\mathbb A$ do have a good fixator (\ref{sse:ExGood}).
This, in turn, gives a lot of applications (\ref{sse:Applications}), like the retraction
associated to a sector-germ. We finish this section with the structure of the
residue buildings (\ref{sse:Residue}).

\subsection{Good fixators}
\label{sse:Good}

When $\Omega\subset\Omega'\subset\mathbb A$, then 
$\hat{P}_{\Omega}\supset\hat{P}_{\Omega'}$.  As $\hat{P}_{x}$ is the fixator of
$x\in\mathbb A$, $\hat{P}_{\Omega}$ is included in  the fixator $G_\Omega$ of $\Omega$ (for
the action of $G(K)$ on $\mathcal I\supset\mathbb A$).  Actually, when $\Omega$ is a set
$G_\Omega=\bigcap_{x\in\Omega}\;\hat{P}_{x}$, and when $\Omega$  is a filter
$G_\Omega=\bigcup_{S\in\Omega}\;G_S$.

 We have $G_\Omega\cap N(K)=\hat{N}_{\Omega}$ and 
$G_\Omega\supset\hat{P}_{\Omega}$ which  contains $U_\Omega^{pm+}.U_\Omega^{nm-}.\hat
N_\Omega$ and 
$U_\Omega^{nm-}.U_\Omega^{pm+}.\hat N_\Omega$.

\medskip
For $\Omega$ a filter of subsets in $\mathbb A$, the subset of $G(K)$ 
consisting  of the $g\in G(K)$ such that $g.\Omega\subset\mathbb A$ is: $G(\Omega\subset\mathbb A)= 
\bigcup_{S\in\Omega}\; (\bigcap_{x\in S}\;N(K).\hat{P}_{x})$. Indeed
 $g.\Omega\subset\mathbb A\;\Leftrightarrow\;\exists S\in\Omega,\;g.S\subset\mathbb A
\;\Leftrightarrow\;\exists S\in\Omega,\;\forall x\in S,\;gx\in\mathbb A
\;\Leftrightarrow\;$ (by \ref{de:Hovel})  $\exists S\in\Omega,\;\forall x\in S,\;g\in
N(K).\hat{P}_{x}$.

\begin{defi} 
\label{de:Good}
Consider the following properties:

\par\qquad (GF+) $G_\Omega=\hat{P}_{\Omega}=U_\Omega^{pm+}.U_\Omega^{nm-}.\hat N_\Omega$,

\par\qquad (GF$-$) $G_\Omega=\hat{P}_{\Omega}=U_\Omega^{nm-}.U_\Omega^{pm+}.\hat N_\Omega$,

\par$\qquad\;$ (TF)$\;$ $G(\Omega\subset\mathbb A)=N(K).G_{\Omega}$.

\par We say that $\Omega$ in $\mathbb A$ {\it has a good fixator} if it satisfies these three
properties.

\par We say that $\Omega$ in $\mathbb A$ {\it has an half-good fixator} if it satisfies (TF)
and (GF+) or (GF$-$).

\par We say that $\Omega$ in $\mathbb A$ {\it has a transitive fixator} if it satisfies (TF).

\end{defi}

 By point a) in the proof of Proposition \ref{sse:3Groups}, this definition doesn't depend on the choice of $\Delta^+$ in its 
$ W^v-$conjugacy class and $N(K)$ permutes the filters with good fixators 
and the corresponding fixators. By \ref{de:Fixator} and \ref{de:Hovel}, a point has a good fixator.

 In the classical case of reductive groups, every $\Omega$ has a good fixator 
and $\hat{P}_{\Omega}$ is as defined by Bruhat and Tits  \cite[7.1.8, 7.1.11, 7.4.8]{BruhatTits72}. 

\begin{rema}
\label{ex:Independence}
If $\Omega$ in ${\mathbb A}$ has a transitive fixator.  Then $G_\Omega$ is transitive
on the apartments containing $\Omega$:  if $g\in G(\Omega\subset\mathbb A)$, there exists
$n\in N(K)$ such that $g\vert_\Omega=n\vert_\Omega$; moreover  if $g.{\mathbb A}\supset\Omega$, then 
$g^{-1}\Omega\subset\mathbb A$ and $g^{-1}=np\in N(K).G_{\Omega}$, so $g.{\mathbb A}=
p^{-1}n^{-1}.{\mathbb A}=p^{-1}.{\mathbb A}$. In particular $G_\Omega$ and all invariant 
subgroups of $G_{\Omega}$ do not depend of the particular choice of the apartment 
$A$ containing $\Omega$.
\end{rema}

\begin{prop} 
\label{pr:Good}
\begin{enumerate}
\item[1)] Suppose $\Omega\subset\Omega'\subset
cl(\Omega)$. If $\Omega$ in $\mathbb A$ has a good (or half-good) fixator,  then this also holds
for $\Omega'$  and $G_{\Omega}=\hat{N}_{\Omega}.G_{\Omega'}$, 
$N(K).G_{\Omega}=N(K).G_{\Omega'}$. In particular, any apartment containing 
$\Omega$ contains its enclosure $cl(\Omega)$.

 Conversely, if $supp(\Omega)=\mathbb A$ (or $supp(\Omega')=supp(\Omega)$, hence
$\hat{N}_{\Omega'}=\hat{N}_{\Omega}$), $\Omega$ has an half-good fixator and $\Omega'$ has 
a good fixator,  then $\Omega$ has a good fixator.

\item[2)] If a filter $\Omega$ in $\mathbb A$ is generated by a family $\mathcal F$ of filters 
with  good (or half-good) fixators, then $\Omega$ has a good (or half-good) fixator
 $G_{\Omega}=\bigcup_{F\in\mathcal F}\;G_{F}$.

\item[3)] Suppose that the filter $\Omega$ in $\mathbb A$ is the union of an increasing 
sequence $(F_i)_{i\in\mathbb N}$  of filters with good (or half-good) fixators and that, for  some
$i$, the space ${\rm supp}(F_i)$ has a finite fixator $W_0$ in $W_Y$, then 
$\Omega$ has a good (or half-good) fixator $G_{\Omega}=\bigcap_{i\in\mathbb N}\;G_{F_i}$.

\item[4)] Let $\Omega$ and $\Omega'$ be two filters in $\mathbb A$. Suppose $\Omega'$ satisfies
(GF+) (resp. (GF+) and (TF)) and that there exist a finite number of
positive, closed, vectorial chambers 
$\overline {C^v_1},\cdots$, $\overline {C^v_n}$ such that: 
$\Omega\subset\cup_{i=1,n}\;\Omega'+\overline {C^v_i}$. Then $\Omega\cup\Omega'$ satisfies
(GF+) (resp. (GF+) and (TF)) and $G_{\Omega\cup\Omega'}=G_{\Omega}\cap G_{\Omega'}$.
\end{enumerate}
\end{prop}

\begin{rema}
\label{re:Good}
In 4) above, the same results are true when changing $+$ to $-$.

 If $\Omega'$ has a good fixator, $\Omega\subset\cup_{i=1,n}\;\Omega'+\overline
{C^v_i}$ and $\Omega\subset\cup_{i=1,n}\;\Omega'-\overline {C^v_i}$, then $\Omega\cup\Omega'$
has a good fixator.

 If $\Omega$ satisfies (GF$-$), $\Omega'$ satisfies (GF+), $\Omega$ or $\Omega'$ satisfies
(TF), $\Omega\subset\cup_{i=1,n}\;\Omega'+\overline {C^v_i}$ and
$\Omega'\subset\cup_{i=1,n}\;\Omega-\overline {C^v_i}$, then $\Omega\cup\Omega'$ has a good
fixator.

\end{rema}

\begin{proof}
\begin{enumerate}
\item[1)] When $\Omega\subset\Omega'\subset cl(\Omega)$, we always have
$U^{pm+}_{\Omega}=U^{pm+}_{\Omega'}= U^{pm+}_{cl(\Omega)}$,
$U^{nm-}_{\Omega}=U^{nm-}_{\Omega'}=U^{nm-}_{cl(\Omega)}$, $G_{\Omega'}\subset G_\Omega$,
$G(\Omega'\subset\mathbb A)\subset G(\Omega\subset\mathbb A)$ and 
$\hat{N}_{\Omega'}= N(K)\cap G_{\Omega'}\subset \hat{N}_{\Omega}$ (with equality when
$supp(\Omega')=supp(\Omega)$);  so the first assertion of 1) is clear. 
 The second assertion is a consequence of Remark~\ref{ex:Independence}.
 
 For the last assertion we know that
$G_{\Omega'}=U_\Omega^{pm+}.U_\Omega^{nm-}.\hat N_{\Omega'}=U_\Omega^{nm-}.U_\Omega^{pm+}.\hat
N_{\Omega'}$, $\hat N_{\Omega'}=\hat N_{\Omega}$ and $G_{\Omega'}=G_{\Omega}$; so the
fixator $G_{\Omega}$ is good.

\item[2)] If $\Omega$ is generated by the family $\mathcal F$ of filters, we have 
\begin{gather*}
G_{\Omega}=\bigcup_{F\in\mathcal F}\;G_{F},\ 
U^{pm+}_{\Omega}=\bigcup_{F\in\mathcal F}\;U^{pm+}_{F},\\
U^{nm-}_{\Omega}=\bigcup_{F\in\mathcal F}\;U^{nm-}_{F}, \ 
\hat{N}_{\Omega}=\bigcup_{F\in\mathcal F}\;\hat{N}_{F}
\end{gather*} and $G(\Omega\subset
\mathbb A)=\bigcup_{F\in\mathcal F}\;G(F\subset\mathbb A)$; so 2) is clear.

\item[3)] If $\Omega$ in $\mathbb A$ is the union of an increasing 
sequence $(F_i)_{i\in\mathbb N}$  of filters, we have
 $G_{\Omega}=\bigcap_{i\in\mathbb N}\;G_{F_i}$, 
$U^{pm+}_{\Omega}=\bigcap_{i\in\mathbb N}\;U^{pm+}_{F_i}$, 
$U^{nm-}_{\Omega}=\bigcap_{i\in\mathbb N}\;U^{nm-}_{F_i}$ and
$G(\Omega\subset \mathbb A)=\bigcap_{i\in\mathbb N}\;G(F_i\subset \mathbb A)$. By hypothesis we may suppose that all 
${\rm supp}(F_i)$ have the same finite fixator $W_0$, so,
$\hat{N}_{\Omega}=\hat{N}_{F_i}=W_0.H$.

 If  $g\in\bigcap_{i\in\mathbb N}\;G_{F_i}=\bigcap_{i\in\mathbb N}\;U^{pm+}_{F_i}.
U^{nm-}_{F_i}.H.W_0$,  by extracting a subsequence, there exists $n_0\in N(K)$
such that $gn_0^{-1}\in \bigcap_{i\in\mathbb N}$ $U^{pm+}_{F_i}.U^{nm-}_{F_i}.H$, and, 
because $U^\pm (K)\cap B^\mp (K)=\{1\}$ (\ref{ssse:Borels}), this intersection is equal to
$U^{pm+}_{\Omega}.U^{nm-}_{\Omega}.H$. So,
$G_\Omega =U^{pm+}_{\Omega}.U^{nm-}_{\Omega}.{\hat N}_\Omega$. 

 If $g\in G(\Omega\subset
\mathbb A)=\bigcap_{i\in\mathbb N}\;N(K)G_{F_i}$, then, for all $i$, $g\in w_iG_{F_i}$ for some
$w_i\in\hat W$, unique modulo $W_0$ as $G_{F_i}\cap N(K)=\hat N_{F_i}=W_0.H$.
Extracting a subsequence, we may suppose that $w_i$ is independent on $i$, so $g\in
w_i.(\bigcap_j\;G_{F_j})=w_i.G_\Omega$ and $G(\Omega\subset\mathbb A)\subset N(K).G_\Omega$.

\item[4)] By induction, we may suppose
$\Omega\subset\Omega'+\overline{C^v_1}$. We may also assume that 
$\overline {C^v_1}$ is the closed positive fundamental chamber  $\overline {C^v_f}$.

Suppose (GF+) and (TF) for $\Omega'$. Let $u\in U^{pm+}_{\Omega'}$, 
$v\in U^{nm-}_{\Omega'}$ and $n\in N(K)$ be such that $uvn\in (G_\Omega.N(K))\cap
(G_{\Omega'}.N(K))$; we now replace $\Omega$ and $\Omega'$ by 
appropriate  sets in these filters. Clearly,
 $U^{pm+}_{\Omega'}= U^{pm+}_{\Omega'+\overline C^v_f}\subset 
U^{pm+}_{\Omega\cup\Omega'}\subset U^{pm+}_{\Omega}$. So, for all $x\in\Omega$, $vn\in G_{x}.N(K)$ and, as a point
has a good fixator, we may write $vn=v'_xu'_xn'_x$ with 
$v'_x\in U^{nm-}_{x}$, $u'_x\in U^{pm+}_{x}$ and $n'_x\in N(K)$. Hence 
$n'_xn^{-1}=(u'_x)^{-1}(v'_x)^{-1}v$ and, by Birkhoff (2.1.7), $n'_x=n$, $u'_x=1$, 
$v'_x=v\in U^{nm-}_{\Omega'}\cap U^{nm-}_{x}=U^{nm-}_{\Omega'\cup\{x\}}$.
Now, we have $u\in U^{pm+}_{\Omega\cup\Omega'}$  and $v\in
\cap_{x\in\Omega}\;U^{nm-}_{\Omega'\cup\{x\}}=U^{nm-}_{\Omega\cup\Omega'}$, so
$uvn\in U^{pm+}_{\Omega\cup\Omega'}.U^{nm-}_{\Omega\cup\Omega'}.N(K)\subset (G_\Omega\cap 
G_{\Omega'}).N(K)$.

Suppose (GF+) for $\Omega'$. Let $uvn$ as above be in $ G_\Omega\cap G_{\Omega'}$, we
have still the same results, but moreover $n\in \hat{N}_{\Omega'}$  and $n'_x\in
\hat{N}_{x}$. So $n=n'_x\in \hat{N}_{\Omega'\cup\{x\}}$, 
$\forall x\in\Omega$, hence $n\in \hat{N}_{\Omega\cup\Omega'}$. Therefore we get
$uvn\in U^{pm+}_{\Omega\cup\Omega'}.U^{nm-}_{\Omega\cup\Omega'}.\hat{N}_{\Omega\cup\Omega'}$.
\end{enumerate}
\end{proof}

\subsection{Examples of filters with good fixators}
\label{sse:ExGood}

\begin{enumerate}
\item[1)] If $x\leq y$ in $\mathbb A$, then $\{x,y\}$, $[x,y]$ and $cl(\{x,y\})$ have 
good fixators and $G_{\{x,y\}}=G_{[x,y]}$. Moreover, if $x\not=y$,
$]x,y]=[x,y]\setminus\{x\}$ has a good fixator: it satisfies (GF-) and (TF) by
Proposition ~\ref{pr:Good}  4) and, as $]x,y]\subset[x,y]\subset cl(]x,y])$, it has a good
fixator by Proposition ~\ref{pr:Good} 1).

\item[2)] A local face in $\mathbb A$ has a good fixator: ${\rm germ}_x(x+F^v)$ is generated by
the sets $F_n=(x+F^v)\cap(y_n-{\overline{F^v}})$ for $y_n=x+{1\over n}\xi$, $\xi\in F^v$ and
$n\in\mathbb N$; moreover (for $F^v\not=\{0\}$) $]x,y_n]\subset F_n\subset cl(]x,y_n])$, so by
1) above and Proposition ~\ref{pr:Good} (1) and 2)) $F_n$ and the local face have a good
fixator. Now ${\rm germ}_x(x+F^v)\subset F(x,F^v)\subset \overline{F}(x,F^v)=
cl({\rm germ}_x(x+F^v))$; so, by Proposition ~\ref{pr:Good} 1), any face or closed face has a good
fixator.

\item[3)] A sector in $\mathbb A$ has a good fixator: $x+C^v$ is the 
increasing union of the sets $F_n=(x+C^v)\cap cl(]x,y_n])$ where $y_n=x+n\xi$,
$\xi\in C^v$ and $n\in\mathbb N$. Moreover these $F_n$ have $\mathbb A$ as support and 
$]x,y_n]\subset F_n\subset cl(]x,y_n])$, so $F_n$ and $x+C^v$ have good fixators.

\item[4)] A sector-germ has a good fixator. The fixator of $\mathfrak S_{\pm \infty}$ 
is $H.U^\pm (K)$, since every element in 
$U^\pm (K)$ is a finite product of elements in groups $U_\alpha(K)$ for 
$\alpha\in\Phi^\pm $. 

 On the contrary, $U^{\max+}(K)$ is not the union of the $U^{\max+}_\Omega$ for 
$\Omega\in\mathfrak S_{+\infty}$.

\item[5)] The apartment $\mathbb A$ itself has a good fixator  $G_{\mathbb A}=H$: 
$\mathbb A$  is the increasing union of $cl(\{-n\xi,n\xi\})$ for $\xi\in C^v_f$.

\item[6)] For the same reasons, a  wall $M(\alpha,k)$  has a good fixator which is 
$$U_{\alpha,k}.U_{-\alpha,-k}.\{1,r_{\alpha,k}\}.H.$$

\item[7)] Exercise:
An half-apartment $D(\alpha,k)$  has a good fixator
$HU_{\alpha,k}$.
 If $x_+-x_-\in\mathcal T^o$, then $cl(\{F(x_-,F^v_-),F(x_+,F^v_+)\})$ 
has a good fixator for all vectorial faces $F^v_\varepsilon$ (where $\varepsilon = \pm$).
\end{enumerate}

\subsection{Applications}
\label{sse:Applications}
\begin{enumerate}
\item[1)] By \ref{sse:ExGood}.5), the fixator (resp. stabilizator) of the apartment $\mathbb A=A_f$ 
is $H$  (resp. N(K)). In particular, the maps $g\mapsto g.\mathbb A$ and $g\mapsto
g.T.g^{-1}$  give bijections $\{$apartments of $\mathcal I(G,K)\}\leftrightarrow
G(K)/N(K)\leftrightarrow
\{$maximal split tori of $G(K)\}$.

 Moreover, the action of $N(K)$ on $\mathbb A$ preserves the affine structure of 
$\mathbb A$,  its lattice of cocharacters $Y$, $\mathcal T$ and $\mathcal T^o$. So, any apartment 
$A$ in $\mathcal I(G,K)$ is endowed with a canonical structure of real affine space, 
an affine  action of a Weyl group $W(A)$, a lattice $Y(A)$ of cocharacter
points, Tits cones and a  preorder relation. More generally, all structures in
$\mathbb A$ invariant under $N(K)$ are  transferred to any apartment by the
$G(K)-$action: in an apartment, the notions  of (spherical)
face, special point, cocharacter point,  wall, sector, sector-germ or filter
with good fixator are well defined  (independently of the apartment containing
them, as they all have good fixators).

 When we speak of an isomorphism between apartments, we mean an affine 
isomorphism  exchanging the walls and the Tits cones.

\item[2)] Let $A_1,A_2$ be two apartments and $x,y$ be two points in  $A_1\cap
A_2$. If $x\leq y$ in $A_1$, then, by Remark \ref{ex:Independence} and \ref{sse:ExGood}.1), there exists $g\in \hat{P}_{cl(x,y)}$,  such that $A_2=g.A_1$, hence $A_1\cap A_2\supset
cl(x,y)$   and $x\leq y$ in $A_2$. In particular, the relation $\leq$ is
defined  on the whole hovel $\mathcal I(G,K)$ (note that $x\leq y$ implies by
definition  that $x$ and $y$ are in a same apartment). We shall
 see below (\ref{sse:Preorder}) that this relation is transitive, so it is a preorder-relation 
(reflexive,  transitive, perhaps not antisymmetric).

 The intersection of two apartments $A_1,A_2$ is {\it order-convex}: if 
$x,y\in A_1\cap A_2$ and $x\leq y$, then the segment $[x,y]$ of $A_1$ is in 
$A_1\cap A_2$  and equal to the corresponding segment in $A_2$. In particular,
any affine subspace of $A_1$ whose direction meets the open Tits cone $\mathcal T^\circ(A_1)$ and which is contained in
$A_1\cap A_2$ is also an affine subspace of $A_2$.

\item[3)] For any face (or any narrow filter) $F$ and any sector germ $\mathfrak S$ in 
$\mathcal I(G,K)$, there exists an apartment $A$ containing $F$ and $\mathfrak S$: Using 
the $G(K)-$action one may suppose $\mathfrak S=\mathfrak S_{\pm \infty}$. Now $F=g.F'$ 
with $F'$ a face in $\mathbb A$. By the Iwasawa decomposition, $g=unv$ with 
$u\in U^\pm (K)\subset G_{\mathfrak S}$, $n\in N(K)$ and $v\in U_{F'}\subset G_{F'}$. 
So $F=un.F'\subset un.\mathbb A=u.\mathbb A$ and $\mathfrak S\subset u.\mathbb A$. 

 By order-convexity, any apartment containing $F$ and $\mathfrak S={\rm germ}$ $(y+
{\overline{C^v}})$ contains $F+C^v$ (and even $cl(F+C^v)\supset
F+{\overline{C^v}}$, when $F$ has a good fixator, by 1) and 4) of Proposition~\ref{pr:Good}).  In
particular, any apartment containing $x$ and $\mathfrak S$ contains the sector $\mathfrak s$
of  direction $\mathfrak S$ and base point $x$.  By \ref{sse:ExGood}.3) and Remark \ref{ex:Independence}, any two such apartments are  conjugated by $G_{\mathfrak s}$.

\item[4)] If $\Omega_1=F(x,F_1^v)$ is a face of base point $x$ and $\Omega_2$ a 
narrow filter  containing $x$, there exists an apartment $A$ containing both of
them: in  an apartment $A_1$ containing $\Omega_1$ we choose a vectorial chamber
$C^v$ such that $\overline{C^v}\supset F_1^v$; now an apartment $A$ containing 
$\Omega_2$ and the germ of the  sector $x+C^v$ contains $\Omega_1$ and $\Omega_2$. If
moreover $\Omega_2$ is also a face, then 
$G_{\Omega_1\cup\Omega_2}$ acts transitively on the apartments containing $\Omega_1$  and
$\Omega_2$ by  \ref{sse:ExGood}.2, Proposition~\ref{pr:Good} 4) and Remark \ref{ex:Independence}. Actually, one can prove that
$\Omega_1\cup\Omega_2$ has a good  fixator when the faces $\Omega_1$ and $\Omega_2$ are of
opposite signs or if one of them is spherical.

 If $C=F(x,C^v)$ is a chamber (in $\mathbb A$) and $M(\alpha,k)$ one of its walls 
(with $C\subset D(\alpha,k)$), then $U_{\alpha,k}=U_{\alpha,C}$ acts transitively on  the
chambers $C'$ adjacent  to $C$ along $M(\alpha,k)$: this is a consequence of \ref{sse:ExGood}.2
and Remark \ref{ex:Independence} as $G_C$ may be written 
$$U_{\alpha,C}.U^{pm}_{C}(\Phi^+\setminus\{\alpha\}).
U^{nm-}_{C}.H$$
and, in this  decomposition, all factors but $U_{\alpha,C}$ fix the
chamber $C_0'$ in $\mathbb A$ adjacent  to $C$ along $M(\alpha,k)$. In particular, any
such chamber $C'$ and $D(\alpha,k)$ are contained  in a same apartment.
\end{enumerate}

\kern40pt

\subsection{Retraction with respect to a sector-germ}
\label{sse:Retraction}

  Let $\mathfrak S$ be a sector-germ in an apartment $A$ of $\mathcal I(G,K)$. For  $x\in
\mathcal I(G,K)$, choose an apartment $A'$ containing $\mathfrak S$ and $x$.  As $\mathfrak S$ has a
good fixator, there exists a $g$ in $\hat{P}_{\mathfrak S}$ such  that $A=g.A'$.
If $g$ and $g'$ are two such elements, then $g^{-1}g'$ induces  an automorphism
of $A'$ fixing the sector germ $\mathfrak S$, hence this  automorphism is the identity:
the map $A'\rightarrow A$, $y\mapsto g.y$ is  unique. Moreover, $\mathfrak S\cup\{x\}$
has a good fixator (see Remark~\ref{re:Good}), so 
$\hat{P}_{{\mathfrak S}\cup\{x\}}$ is transitive on the possible apartments $A'$: 
the point $g.x$ does not depend on the choice of $A'$. So, one may define 
$\rho_{A,{\mathfrak S}}(x)=g.x$.

\begin{defi} 
\label{de:Retraction}
The map $\rho=\rho_{A,{\mathfrak S}}$: $\mathcal I\rightarrow A$, $x\mapsto \rho_{A,{\mathfrak S}}(x)$ is the {\it retraction of $\mathcal I$ onto $A$ with center $\mathfrak S$}. It depends only on $A$ and $\mathfrak S$.
\end{defi}

The restriction of $\rho$ to $A$ is the identity.  It is clear that, up to canonical 
isomorphisms, $\rho_{A,{\mathfrak S}}$ depends  only on $\mathfrak S$. We set $\rho_{\pm
\infty}=\rho_{\mathbb A,{\mathfrak S}_{\pm \infty}}$.

A segment-germ $[x,y)$ for $x\neq y$ in an apartment $A$ (cf. \ref{ssse:Faces})
 is  a narrow filter. When $x\leq y$ (resp. $y\leq x$), its enclosure is a 
 closed-face and $[x,y)$ has a good fixator (\ref{sse:ExGood}.1) and 
3) of Proposition~\ref{pr:Good}); we say that $[x,y)$  is {\it positive} (resp. {\it negative}) and that $[x,y]$ and $[x,y)$ are {\it generic}.

 For any sector-germ $\mathfrak S$ and any segment-germ $[x,y)$, there exists an 
apartment  containing $\mathfrak S$ and $[x,y)$, {\it i.e.} containing $[x,z]$ for some
$z\in[x,y]\setminus\{x\}$. 

 A segment $[x,y]$ in an apartment is compact and, for 
$z\in[x,y]$, the set $[z,x)\cup[z,y)$ is an open neighbourhood of $z$. So, if 
$\mathfrak S$ is a sector-germ, there exist an integer $n$, points
$x_0=x,x_1,\dots,x_n=y\in[x,y]$  and apartments $A_1,A_2,\dots,A_n$ such that
$A_i$ contains $\mathfrak S$ and $[x_{i-1},x_i]$. As  a consequence, for all apartments
$A'$ containing $\mathfrak S$, $\rho_{A',{\mathfrak S}}([x,y])$  is the piecewise linear path
 $[\rho x_0,\rho x_1]\cup[\rho x_1,\rho x_2]\cup\dots\cup[\rho x_{n-1},\rho x_n]$.

 We shall give a better description of this piecewise linear path when 
$x\leq y$ in the last section.

\begin{rema}
\label{re:Micro}
The fixator of some spherical sector-face-germ $\mathfrak F={\rm germ}(x+F^v)$ 
contains clearly the group $P(F^\mu)$ associated in \cite{Rousseau06} to a
microaffine face $F^\mu=F^v\times{\mathcal F}$ (for some $\mathcal F$ containing $x$); 
and it was proved in [loc. cit.; 3.5] that $G=P(F^\mu).N.P(E^\mu)$ for any microaffine
faces $F^\mu$ and $E^\mu$ of the same sign. Actually, using Proposition~\ref{pr:Iwasawa} above, the proof of this result is still valid if only one among $E^v$ and $F^v$ is spherical and
the signs of $E^v$ and $F^v$ may be opposite. So, any two sector-face-germs in ${\mathcal I}$ 
are contained in a same apartment, if at least one of these sector-face-germs
is spherical. For an abstract definition of affine hovels, this property is used in \cite{Rousseau08a} as a good substitute to axioms (A3) and (A4) in Tits' definition of affine
buildings \cite{Tits86}, see also \cite[Appendix 3]{Ronan89}. 
\end{rema}

\subsection{Residue buildings} 
\label{sse:Residue}

Let us denote the set of all positive (resp. negative) segment-germs 
$[x,y)$  with $x<y$ (resp. $y<x$) by $\mathcal I_x^+$ (resp. $\mathcal I_x^-$). The set $\mathcal I_x^+$  (resp. $\mathcal I_x^-$) can be given two structures of a building. An
apartment $\mathfrak a_+$ in $\mathcal I_x^+$ (resp. $\mathfrak  a_-$ in $\mathcal I_x^-$) is the 
intersection $A\cap\mathcal I_x^+$  (resp. $A\cap\mathcal I_x^-$) of an apartment $A$ of 
$\mathcal I$ containing $x$ with $\mathcal I_x^+$   (resp. $\mathcal I_x^-$) (or, more precisely,
the set of all $[x,y)$ for $y\in A$ and $x<y$ (resp. $x>y$)). Now, on any 
apartment $\mathfrak a_\pm$, one can put two structures of a  Coxeter complex:
\begin{itemize}
\item the restricted one, modelled on $(W_x^{\rm min}, S_x)$, where $W_x^{\rm min}$ 
is the subgroup  of $W$ generated by all reflections with respect to true walls
passing through $x$,  and where $$S_x = \{s_H \mid H \hbox{ is a wall of } F(x,
C_f^v)\hbox{ containing }x\}$$ may be infinite.  We restrain the action of $W_x^{\rm min}$ on $\mathfrak
a_\pm$. The faces of this structure are the faces 
$F(x,F^v)$ with $F^v$ positive (resp. negative) (or, more precisely, the set 
of all segment-germs $[x,y)$ contained in $F(x,F^v)$).
 \item the unrestricted one, modelled on $(W^v,S^v)$, where we force $x$ to 
be a special point and consider the faces ${\rm germ}_x(x + F^v)$ (local-face in $x$) 
with $F^v$  positive (resp. negative) (note that ${\rm germ}_x(x + F^v)\subset 
F(x,F^v)$).  So we add new (ghost) walls $M(\alpha,k)$ for $\alpha\in\Phi$ and
$k\in\mathbb R\setminus\mathbb Z$,
$\alpha(x)+k=0$.
\end{itemize}

\begin{prop}  
\label{pr:Residue}
The set $\mathcal I_x^+$ or $\mathcal I_x^-$, endowed with its apartments with their restricted (resp. unrestricted) structures of Coxeter complex, is a building.
\end{prop}

\begin{proof}
We have to verify the last two axioms of a building 
(as in \cite[IV.1]{Brown89} or \cite[2.4.1]{Remy02}). 
We focus on the positive case, the negative one is obtained in the same way. 
In both Coxeter structures, \ref{sse:Applications}.4) shows that, given two faces $F_1$ and $F_2$, 
there exists an apartment containing $x$ and both of them. Further, the group 
$G_{F_1\cup F_2}$ acts transitively on the apartments 
containing $F_1\cup F_2$. Hence, for any two such apartments $A$ and $A'$, there 
exists  an element $g\in G_{F_1\cup F_2}$ such that $A' = g\cdot A$ which also
gives  an isomorphism $\mathfrak a_+\simeq \mathfrak a_+'$ fixing $F_1\cup F_2$.
\end{proof}

Note that the unrestricted building structure can be thick  only when it 
coincides  with the restricted one, {\it i.e.} when $x$ is special (thick means that
any panel is a face of at least three chambers). The buildings $\mathcal I_x^\pm$ may
be spherical for the restricted structures, as $W_x^{\rm min}$ may be finite when
$x$ is not special.

Now, we consider  $\mathcal I_x^\pm$ endowed with the unrestricted structure. 
On the set of chambers $Ch(\mathcal I_x^\varepsilon)$ of $\mathcal I_x^\varepsilon$ (for a sign 
$\epsilon\in \{+,-\}$),  we have a distance $d_\varepsilon :  Ch(\mathcal I_x^\pm)\times Ch(\mathcal I_x^\pm) \to W^v$ defined as follows, cf. \cite[2.2 to 2.4]{Remy02}.

 If  $c,c'\in Ch(\mathcal I_x^\varepsilon)$, choose an apartment $A$ containing $c,c'$ and 
a chamber $c_0={\rm germ}_x(x+\varepsilon C_0)$ in $\mathfrak  a_\varepsilon=A\cap \mathcal I_x^\varepsilon$; identify
$(A,C_0)$  to $(A_f, C^v_f)$; this enables us to identify 
$W^v(A)$ to $W^v=W^v(A_f)$. Now, if $c=wc_0$ and $c'=w'c_0$ for some $w,w'\in 
W^v$, then $d_\varepsilon (c,c') = w^{-1}w'$. Note that if we choose $c_0=c$, then 
$c'=d_\varepsilon(c,c').c$. 

 Now, we define a codistance 
$$
d_x^* : \big (Ch(\mathcal I_x^+)\times Ch(\mathcal I_x^-)\big )\cup \big (Ch(\mathcal I_x^-)\times 
Ch(\mathcal I_x^+)\big )\to W^v
$$ in the following way. If $(c,e)\in Ch(\mathcal I_x^\varepsilon)\times Ch(\mathcal I_x^{-\varepsilon})$,
by \ref{sse:Applications}.4), there exists an apartment $A$ containing $x$, $c$ and $e$, unique up 
to isomorphism. If $c' = {\rm germ}_x(x + C')$ is a chamber in $A$, we denote the 
chamber  opposite $c'$ in $A$ by $-c'$, {\it i.e.}  $-c' = {\rm germ}_x(x-C')$.
Choose a chamber
$c_0={\rm germ}_x(x\pm\varepsilon C_0)$ in $\mathfrak  a_{\pm\varepsilon}=A\cap \mathcal I_x^{\pm\varepsilon}$ and identify 
$(A,C_0)$  to $(A_f, C^v_f)$. If $c=\pm w_1.c_0$ and 
$e=\mp w_2.c_0$, the codistance between $c$ and $e$ is then  $d_x^* (c,e) =
d_{-\varepsilon}(-c,e) =  d_{\varepsilon}(c,-e)=w_1^{-1}w_2$. It does not depend on the
choices. 

\begin{prop} 
\label{pr:Twin}
The two buildings $\mathcal I_x^+$ and $\mathcal I_x^-$, endowed with their unrestricted structures of buildings and the codistance $d_x^* $, form a twinned pair of buildings.  
\end{prop}

\begin{enonce*}[remark]{N.B} With analogue arguments, one shows that this still holds if  
$\mathcal I_x^+$ and $\mathcal I_x^-$ are endowed with the restricted structures.
\end{enonce*}

\begin{proof}
We have to check the axioms of  twinning as given in \cite[2.2]{Tits92}, see also \cite[2.5.1]{Remy02}.

 Indeed, the first axiom (Tw1) is fulfilled:  $d_x^*(e,c) =w_2^{-1}w_1
=d_x^*(c,e)^{-1}$. Let now $c\in Ch(\mathcal I^\varepsilon_x)$ and $e, e'\in  Ch(\mathcal I^{-\varepsilon}_x)$ be chambers  such that $d_x^*(c,e) = w$ and $d_{-\epsilon}(e, e')
= s \in S^v$ with  $\ell (ws) = \ell (w)-1$.
 Let $A$ be an apartment containing $x$, $c$ and $e$ and choose $c_0=-c$; since 
$\ell (ws) = \ell (w) -1$, the wall~$H$ generated by the panel of $e$  of 
type $\{s\}$ separates the latter from~$-c$. In other words, $c$ 
and $e$ are on the same side of $H$. Therefore, by \ref{sse:Applications}.4) there exists an 
apartment $A'$ containing $c$, $e$ and $e'$. In this apartment, $e=w.(-c)$ and 
$e'=(wsw^{-1}).e$, so, $e'=ws.(-c)$ and $d_x^*(c,e') = ws$. This is the  second
axiom (Tw2). 

 To check the third axiom (Tw3), let again $c$ and $e$ be
two chambers such that $d_x^*(c,e) = w$, and let $s\in S^v$. In an apartment $A$ 
containing $c$ and $e$, the chamber $h$ adjacent to $e$ along the panel of type 
$\{s\}$ satisfies $d_{-\epsilon}(e, h) = s $ and $d_x^*(c, h) = ws$. 
\end{proof}

An apartment $A$ of $\mathcal I$ containing $x$ gives a twin apartment  $\mathfrak  a=\mathfrak 
a_+\cup\mathfrak  a_-$,  where $\mathfrak  a_\pm=A\cap\mathcal I_x^\pm$. If $c_0$ is a chamber in $\mathfrak 
a$, there is (as in any twin  building) a retraction $\rho$ of center $c_0$ of
$\mathcal I_x^+$ onto $\mathfrak  a_+$ and of $\mathcal I_x^-$  onto $\mathfrak  a_-$; it preserves the
distances or codistances to $c_0$.

\section{Littelmann paths}
\label{se:Littelmann}

In this paragraph, we give a brief account and
some new results on Littelmann's theory of paths \cite{Littelmann94},  \cite{Littelmann95}, \cite{Littelmann98}. First,  we recall the definitions of $\lambda-$paths, billiard paths, LS paths, and Hecke paths; we then compare the last two notions (see Section \ref{sse:LambdaPaths}). In analogy with \cite{GaussentLittelmann05} where the dimension of a gallery is defined, we introduce some statistics on paths used to characterize the LS paths (see Section \ref{sse:Statistics} and \ref{sse:NewLS}). Note that the symmetrizability of $\mathfrak g$ assumed since \ref{sse:Parahoric} is useless in this section.

\subsection{$\lambda-$paths}
\label{sse:LambdaPaths} 

We consider piecewise linear continuous paths 
$\pi:[0,1]\rightarrow \mathbb A$ such that each (existing) tangent vector $\pi'(t)$
 is in an orbit $W^v.\lambda$ of some $\lambda\in {\overline{C^v_f}}$ under the 
vectorial Weyl group $W^v$. Such a path is called a {\it $\lambda-$path}; it is 
increasing with respect to the preorder relation of \ref{ssse:AffineWeyl}. If $\pi(0)$, $\pi(1)$  and~
$\lambda$  are in $Y$, we say that $\pi$ is ``in $Y$''.

 For any $t\neq 0$ (resp. $t \neq1$), we let 
$\pi'_-(t)$ (resp. $\pi'_+(t)$) denote the derivative of $\pi$ at $t$ from the left 
(resp. from the right). Further, we define $w_\pm(t)\in W^v$ to be the smallest 
element  in its 
$W^v_\lambda-$class such that $\pi'_\pm(t)=w_\pm(t).\lambda$ (where $W^v_\lambda$ is the fixator 
in
$W^v$ of $\lambda$). Moreover, we denote by $\pi_-(t)=\pi(t)-[0,1)\pi_-'(t)=
[\pi(t),\pi(t-\varepsilon)\,)$ (resp. $\pi_+(t)=\pi(t)+[0,1)\pi_+'(t)=
[\pi(t),\pi(t+\varepsilon)\,)$ (for $\varepsilon>0$ small) the positive (resp. negative) 
segment-germ of $\pi$ at $t$ (cf. \ref{ssse:Faces}).

 The reverse path $\overline\pi$ defined by  
${\overline\pi}=\pi(1-t)$ has symmetric properties, it is a $(-\lambda)-$path.

 If, for all $t$, $w_+(t)\in W^v_{\pi(t)}.w_-(t)$, we shall say that $\pi$ is
a {\it billiard path}. This corresponds to what is stated in \cite[Lemma 4.4]{KapovichMillson05},  but seems stronger than the definition given in [loc. cit.; 2.5] which looks
more like our definition of $\lambda-$path.

 For any choices of $\lambda\in {\overline{C^v_f}}$, $\pi_0\in \mathbb A$, 
$r\in\mathbb N\setminus\{0\}$ and sequences 
${\underline\tau}=(\tau_1,\tau_2,\dots,\tau_r)$ of elements in $W^v/W^v_\lambda$ and 
${\underline a}=(a_0=0<a_1<a_2<\dots<a_r=1)$ of elements in $\mathbb R$, we define 
a $\lambda-$path $\pi=\pi(\lambda,\pi_0,{\underline\tau},{\underline a})$ 
by the formula:
$$
\pi(t)=\pi_0+\sum_{i=1}^{j-1}\;(a_i-a_{i-1})\tau_i(\lambda)+(t-a_{j-1})\tau_j(\lambda)
\quad \hbox{ for } \quad a_{j-1}\leq t\leq a_j.
$$

\noindent Any $\lambda-$path may be defined in this way. We shall always assume 
$\tau_j\neq\tau_{j+1}$.

\subsubsection{LS and Hecke paths}  
\label{ssse:LSandHecke}

We consider now more specific paths.

\begin{defi}
\label{de:LS}
\cite{Littelmann94} 
A {\it Lakshmibai-Seshadri path} (or {\it LS path}) of shape 
$\lambda\in Y^+$ is a $\lambda-$path $\pi=\pi(\lambda,\pi_0,{\underline\tau},{\underline a})$ 
starting in $\pi_0\in Y$ and such that:
for all $j=1,\dots,r-1$, there exists an {\it $a_j-$chain} from $\tau_j$ to 
$\tau_{j+1}$ {\it i.e.} a sequence of cosets in $W^v/W^v_\lambda$: 
$$\sigma_{j,0}=\tau_{j}\;,\;\sigma_{j,1}=r_{\beta_{j,1}}\tau_j\;,\dots,\;
\sigma_{j,s_j}=r_{\beta_{j,s_j}}\dots r_{\beta_{j,1}}\tau_j=\tau_{j+1}$$
where $\beta_{j,1},\dots,\beta_{j,s_j}$ are positive real roots such that, for all 
$i=1,...,s_j$: 
\begin{itemize} 
\item[i)] $\sigma_{j,i}<\sigma_{j,i-1}$,\quad for the Bruhat-Chevalley order
 on $W^v/W^v_\lambda$,
\item[ii)] $a_j\beta_{j,i}(\sigma_{j,i}(\lambda))\in\mathbb Z$,
\item[iii)] $\ell_\lambda(\sigma_{j,i})=\ell_\lambda(\sigma_{j,i-1})-1$, here 
$\ell_\lambda(-)$  is the length in $W^v/W^v_\lambda$.
\begin{enonce*}[remark]{N.B} Actually Littelmann requires the following additional condition 
\end{enonce*}
\item[iv)] $\pi$ is normalized {\it i.e.} $\pi_0=0$.
\end{itemize}
\end{defi}

\begin{defi}
\label{de:Hecke}
\cite[3.27]{KapovichMillson05} 
A {\it Hecke path} of shape $\lambda$ is a 
$\lambda-$path such that, for all $t\in [0,1]\setminus\{0,1\}$, 
$\pi_+'(t)\leq_{W^v_{\pi(t)}}\pi_-'(t)$, which means  that there exists a
${W^v_{\pi(t)}}-$chain from $\pi_-'(t)$ to $\pi_+'(t)$, {\it i.e.} finite sequences 
$(\xi_0=\pi_-'(t),\xi_1,\dots,\xi_s=\pi_+'(t))$ of vectors in $V$ and 
$(\beta_1,\dots,\beta_s)$  of positive real roots such that, for all $i=1,\dots,s$:
\begin{itemize}
\item [v)] $r_{\beta_i}(\xi_{i-1})=\xi_i$,
\item[vi)] $\beta_i(\xi_{i-1})<0$,
\item[vii)] $r_{\beta_i}\in{W^v_{\pi(t)}}$ {\it i.e.} $\beta_i(\pi(t))\in\mathbb Z$: 
$\pi(t)$ is in a wall of direction Ker$(\beta_i)$.
\end{itemize}
\end{defi}

\begin{remas} 
Conditions v) and vii) tell us that $\pi$ is a billiard path. More precisely, the path is folded at $\pi(t)$ by applying successive reflections along the walls $M(\beta_i,-\beta_i(\pi(t))\,)$. Moreover condition vi) 
tells us that the path is ``positively folded'' (cf. \cite{GaussentLittelmann05}).

The definition of affine paths in [Littelmann-98] is a little bit different; 
in particular, it is stable by concatenation.
\end{remas}

\subsubsection{LS versus Hecke} 
\label{ssse:LSvHecke}

Let $\pi=\pi(\lambda,\pi_0,{\underline\tau}, {\underline
a})$ be a $\lambda-$path.  The conditions in Definition \ref{de:Hecke} are trivially satisfied for $t\neq
a_1,\dots,a_{r-1}$.  So, we compare conditions i), ii), iii) to conditions v), vi),
vii) at $t=a_j$, $1\leq j\leq r-1$, for $s=s_j$ and $\beta_i = \beta_{j,i}$. As
$\pi_-'(t)=\tau_j(\lambda)$, the  condition v) tells us that $\xi_i=\sigma_{j,i}(\lambda)$.

\begin{lemm} 
\label{le:Equiv1}
Conditions i) and vi) are equivalent. If they 
are satisfied (for all $i=1,\dots,s_j$), then $w_+(t)< w_-(t)$ in the 
Bruhat-Chevalley order of $W^v/W^v_\lambda$.
\end{lemm}

\begin{proof}
This is clear as $\sigma_{j,i}=r_{\beta_{j,i}}\sigma_{j,i-1}$ and 
$\xi_{i-1}=\sigma_{j,i-1}(\lambda)$.
\end{proof}

\begin{rema} 
When $\lambda$ is in $Y$, the conditions i) and ii) tell us that 
$a_j\in\mathbb Q$ (as required by Littelmann for LS paths).
\end{rema}

\begin{lemm} 
\label{le:Equiv2}
Suppose that $\pi_0\in Y$, $\lambda\in Y^+$ and that 
conditions ii) are satisfied for $1\leq j'<j$ 
and $1\leq i\leq s_{j'}$. Then the set of conditions ii) for $1\leq i\leq s_{j}$ 
(and this $j$) is equivalent to the set of conditions vii) for $1\leq i\leq s_{j}$ 
(and $t=a_j$).
 If $\pi_0\in Y$, $\lambda\in Y^+$ and conditions ii) (or vii)) are satisfied 
for all  $1\leq j\leq r-1$ and 
$1\leq i\leq s_j$, then $\pi(1)\in Y$, hence $\pi$ is in~$Y$.
\end{lemm}

\begin{proof}
From the definition, one has 
$$\pi(a_j)=\pi_0+\sum_{i=1}^{j}\;(a_i-a_{i-1})\tau_i(\lambda)
=\pi_0+a_j\tau_j(\lambda) + \sum_{i=1}^{j-1}\;a_i(\tau_i(\lambda)-\tau_{i+1}(\lambda))$$ 
and (with the $\sigma_{j,i}$ as in Definition~\ref{de:LS}):
$$a_j(\tau_{j+1}(\lambda)-\tau_j(\lambda))=\sum_{i=1}^{s_j}\;a_j(\sigma_{j,i}(\lambda)-\sigma_{j,i-1}(\lambda))
=\sum_{i=1}^{s_j}\;a_j \beta_{j,i} (\sigma_{j,i}(\lambda))\beta^\vee{_{j,i}}.$$ 
Hence, the conditions (ii) for $1\leq i\leq s_j$ imply that 
$a_j(\tau_{j+1}(\lambda)-\tau_j(\lambda))\in Q^\vee\subset Y$. In particular, 
conditions ii) for all $i,j$ imply that $\pi(1)\in~Y$. 

One has $ a_j\beta_{j,i}(\sigma_{j,i}(\lambda)) = r_{\beta_{j,1}}\cdots r_{\beta_{j,i}}(\beta_{j,i})(a_j\tau_j(\lambda))$, so the condition ii) above for $i=1,\dots,s_j$ may be written:
$$ r_{\beta_{j,1}}\cdots r_{\beta_{j,i}}(\beta_{j,i})(a_j\tau_j(\lambda))\in\mathbb Z.$$

\noindent It is easy to verify that these conditions, for all $i=1,\dots,s_j$, mean 
that the roots 
$\beta_{j,i}$ satisfy $ \beta_{j,i}(a_j\tau_j(\lambda))\in\mathbb Z$. If we assume 
ii)  for each $j'<j$ and $i\leq s_{j'}$, this is equivalent 
to $\beta_{j,i}(\pi(a_j))\in\mathbb Z$. 
\end{proof}

Any LS path $\pi$ is a Hecke path in $Y$. The  reverse
path $\overline\pi$ has symmetric properties. The reverse
path of a Hecke path in $Y$ has  also symmetric properties.

 Conversely, any Hecke path $\pi$ of shape $\lambda$ in $Y$ is not far from  being
a LS path. Condition iii)  only is missing. Actually, by condition i) one has 
$s_j\leq \ell_\lambda(\tau_j) - \ell_\lambda(\tau_{j+1})$; so condition iii) is equivalent to 
$s_j= \ell_\lambda(\tau_j) - \ell_\lambda(\tau_{j+1})$. Hence $\pi$ is 
a LS path if and only if the $W^v_{\pi(t)}-$chains are of maximal lengths. See \cite[Proposition 3.24]{KapovichMillson05} for a more precise statement.

\subsection{Statistics on paths}
\label{sse:Statistics}

We define two statistics on $\lambda-$paths and compare them with the one one would have defined inspired by \cite{GaussentLittelmann05}.

\subsubsection{Dual and co-dimension}
\label{ssse:DdimCodim}

\begin{defi}
\label{de:DdimCodim}
The {\it dual 
dimension} of a $\lambda-$path $\pi$, denoted by ${\rm d\!\dim} (\pi)$ and the  {\it
codimension} of  $\pi$,  denoted  by ${\rm codim}(\pi)$, are the non negative
integers:
$${\rm d\!\dim} (\pi) = \sum_{t>0} \ell_{\pi(t)}(w_-(t)),\qquad
{\rm codim} (\pi) = \sum_{t<1} \ell_{\pi(t)}(w_+(t)),$$ 
where $\ell_{\pi(t)}(\ )$ is the relative length function on the Coxeter group 
$W^v$ with respect to $W^v_{\pi(t)}$ defined as follows:  $\ell_{\pi(t)}(w)$ is 
the  number of walls $M(\alpha)$ for $\alpha\in\Phi^+_{\pi(t)}$ separating the
fundamental chamber 
$C^v_f$ from $wC^v_f$; it coincides with the usual length on $W^v_{\pi(t)}$.
\end{defi}

It seems that the sums are infinite, but, actually, there are only a finite 
number of  possible $w_-(t)$ or $\pi'_-(t)=w_-(t)\lambda$ (resp. $w_+(t)$ or
$\pi'_+(t)=w_+(t)\lambda$). Moreover, for any $t$, $\ell_{\pi(t)}(w_-(t))$ (resp.
$\ell_{\pi(t)}(w_+(t))$) is the number of roots $\beta\in\Phi^+_{\pi(t)}$ such that
$\beta(\pi'_-(t))<0$ (resp. $\beta(\pi'_+(t))<0$). Hence ${\rm d\!\dim} (\pi)$ (resp.
${\rm codim} (\pi)$) is the number of pairs $(t,M(\beta,k))$ consisting of a $t>0$
(resp. $t<1$) and a wall associated to $\beta\in\Phi^+$ such that 
$\pi(t)={\overline \pi}(1-t)\in M(\beta,k)$ and  $\pi(t-\varepsilon)={\overline
\pi}(1-t+\varepsilon)\in D^\circ(\beta,k)=\mathbb A\setminus D(-\beta,-k)$ (resp.
$\pi(t+\varepsilon)\notin D(\beta,k)$), for all small $\varepsilon>0$; this number is clearly finite.

 To be short, ${\rm d\!\dim} (\pi)$ is the number (with multiplicities) of  all
walls  positively leaved by the reverse path $\overline\pi$ (load-bearing  walls
for $\overline\pi$  as in \cite{GaussentLittelmann05}); and ${\rm codim} (\pi)$ is
the number (with multiplicities) of all walls negatively leaved by $\pi$.

 In the following, for $\beta\in\Phi^+$ and $\pi$ a $\lambda-$path, we define 
$pos_\beta(\pi)$  (resp. $neg_\beta(\pi)$) as the number (with multiplicities) of
walls of direction Ker$(\beta)$  leaved positively (resp. negatively) by $\pi$.
Hence:
$${\rm d\!\dim} (\pi)=\sum_{\beta>0}pos_\beta({\overline\pi}) {\rm\quad and\quad }
{\rm codim}(\pi)=\sum_{\beta>0}neg_\beta({\pi})$$

\subsubsection{Classical case} 
\label{ssse:Classical}

Let $\pi$ be a $\lambda-$path in $Y$ and set
$\nu=\pi(1)-\pi(0)$. If $\Phi$ is finite, \cite{GaussentLittelmann05} suggests us to define the 
dimension  of $\pi$ as: $\dim(\pi)=\sum_{\beta>0}pos_\beta({\pi})$ (so, ${\rm
d\!\dim}(\pi)=\dim({\overline\pi})$) and to prove (for Hecke paths) that
$\dim(\pi)\leq\rho(\lambda+\nu)$
 where $\rho_{\Phi^+}=\rho$ is defined by $2\rho=\sum_{\beta>0}\ \beta$.

Actually,
$\beta(\nu)=pos_\beta(\pi)-pos_\beta({\overline\pi})=neg_\beta({\overline\pi})-neg_\beta(\pi)$. 
So, $\dim(\pi)\leq\rho(\lambda+\nu)$ if and only if $\sum_{\beta>0}pos_\beta({\pi})\leq
\rho(\lambda-\nu)+\sum_{\beta>0}\;\beta(\nu)=\rho(\lambda-\nu)+\sum_{\beta>0}pos_\beta(\pi)-
\sum_{\beta>0}pos_\beta({\overline\pi})$ if and only if ${\rm d\!\dim}(\pi)\leq\rho(\lambda-\nu)$.

First, one has
$pos_\beta(\pi)+neg_\beta(\pi)=neg_\beta({\overline\pi})+pos_\beta({\overline\pi})$. 
Further, ${\rm \dim}(\pi)+{\rm codim}(\pi)= \sum_{\beta>0}
(pos_\beta(\pi)+neg_\beta(\pi)) =
\sum_{\beta>0} (neg_\beta({\overline\pi})+pos_\beta({\overline\pi}))$ 
is the number of pairs $(t,M(\beta,k))$ consisting of
a $t<1$ (resp. $t>0$)  such that $\pi(t)\in M(\beta,k)$ and 
$\pi_+(t)\not\subset M(\beta,k)$ (resp. $\pi_-(t)\not\subset M(\beta,k)$). This 
number is invariant if we replace $\pi$ by $\pi_1$ defined by: $\pi_1(t)=\pi(t)$
for $t\leq t_1$ and
$\pi_1(t)=w\pi(t)$ for $t\geq t_1$, for some $t_1\in[0,1]$ and  $w\in
W_{\pi(t_1)}^{\rm min}$. In addition, any billiard path of shape $\lambda$ is obtained
by a sequence of such transformations starting from the straight $\lambda-$path
$\pi_\lambda$ ($\pi_\lambda(t)=t\lambda$). So, 
${\rm \dim}(\pi)+{\rm codim}(\pi)={\rm \dim}(\pi_\lambda)+{\rm codim}(\pi_\lambda) =
{\rm \dim}(\pi_\lambda)=\sum_{\beta>0}\;\beta(\lambda)=\rho(2\lambda)$. Therefore, for any 
billiard path $\pi$ in $Y$, $\dim(\pi)\leq\rho(\lambda+\nu)$ if and only if ${\rm
codim}(\pi)\geq\rho(\lambda-\nu)$.

\subsection{A new characterization of LS paths}
\label{sse:NewLS}

The goals of this section are to prove, in case $\Phi$ is infinite, the inequalities 
${\rm codim}(\pi)\geq\rho(\lambda-\nu)\geq {\rm d\!\dim}(\pi)$  for Hecke
paths in $Y$, and to obtain a new characterization of LS paths.  We choose $\rho_{\Phi^+}=\rho\in X$ such 
that $\rho(\alpha^\vee) =1$ for all simple roots
$\alpha$. It  is clear that $\lambda-\nu$ is a linear combination of coroots; so
$\rho(\lambda-\nu)$  does not depend on the choice of $\rho$.

\subsubsection{The characterization}
\label{ssse:CharacLS}

\begin{prop}
\label{pr:Inequalities}
Let $\pi$ be a Hecke path of shape $\lambda$ in 
$Y$ and $\nu=\pi(1)-\pi(0)$. Then  
$${\rm d\!\dim} (\pi) \leq \rho(\lambda - \nu)\leq{\rm codim} (\pi) 
\quad{\rm and}\quad {\rm d\!\dim}(\pi)+{\rm codim}(\pi)=2\rho(\lambda - \nu) $$ 
with equality if 
and only if $\pi$ is a LS path.
\end{prop}

The proof of Proposition \ref{pr:Inequalities} follows the same strategy as the proof of Proposition 4 in 
\cite{GaussentLittelmann05} and occupies the next three subsections.

\begin{coro}
\label{co:FiniteHecke}
Let $y_0,\,y_1\in Y$ and $\lambda\in Y^+$. Then the
number of Hecke paths $\pi$ of shape $\lambda$ starting in $y_0$ and ending in 
$y_1$ is finite. 
\end{coro}

\begin{enonce*}[remark]{N.B}
Using Littelmann's path model, it was already  clear that  the
number of LS paths satisfying the same conditions is finite, but our proof is
purely combinatorial.
\end{enonce*}

\begin{proof}
By Proposition \ref{pr:Inequalities} and the definition of ${\rm codim}$,  $\ell_{\pi(0)}(w_+(0))\leq{\rm
codim}(\pi)\leq2\rho(\lambda-\nu)$, with $\nu=y_1-y_0$. As $\pi(0)$ is a special
point, this means that there is a finite number of possible $w_+(0)$. So there
is a finite number of possible $w_{\pm}(t)$, or
$\sigma_{j,i}\,$, or $\beta_{j,i}\,$, or $a_j$ satisfying conditions i) and ii)  of
Definition \ref{de:LS}. In conclusion, the number of Hecke paths
$\pi=\pi(\lambda,y_0,{\underline\tau},{\underline a})$ is finite (perhaps zero).
\end{proof}

\subsubsection{The operator $\tilde e_\alpha$}
\label{ssse:Operator}

\begin{defi} 
\label{de:Operator}
Let $\pi$ be a $\lambda-$path and $\alpha$ a simple root. Set $Q = \min \{\alpha(\pi([0,1]))\cap\mathbb Z\}$, the minimal  integral value attained by the 
function
$\alpha(\pi(\ ))$ and let $q$ be the  greatest number in $[0,1]$ such that
$\alpha(\pi([0,q]))\geq Q$. If $q<1$ ({\it i.e.} if 
$Q>\min \{\alpha(\pi([0,1]))\}$), let $\theta > q$ be such that 
$$ \alpha(\pi(q))=\alpha(\pi(\theta))=Q \quad \hbox{ and } \quad
\alpha(\pi(t))<Q \hbox{ for } q<t<\theta.$$
 We cut the path $\pi$ into three parts in the following way. Let  $\pi_1,
\pi_2$ and $\pi_3$ be the paths defined by
$$
\pi_1(t) = \pi(tq);\quad \pi_2(t) = \pi(q+t(\theta - q)) - \pi(q);  \quad
\pi_3(t) = \pi(\theta+t(1-\theta)) - \pi(\theta)
$$ for $t\in [0,1]$. Then, by definition, $\pi = \pi_1 * \pi_2 * \pi_3$,  where
$*$  means the concatenation of paths as defined in \cite[1.1]{Littelmann94}. The
path $\tilde e_\alpha \pi$ is equal to $\pi_1 * r_\alpha (\pi_2) * \pi_3$. 
After a suitable reparametrization $\tilde e_\alpha \pi$ is a $\lambda-$path in $Y$.
\end{defi}

We use also the operators $e_\alpha$ and 
$f_\alpha$ ($\alpha$ simple) defined by Littelmann in \cite[1.2 and 1.3]{Littelmann94}. We do not recall the complete definition here, but note that when they
exist, 
$e_\alpha\pi = \pi_1 * r_\alpha(\pi_2) * \pi_3$ (resp. 
$f_\alpha\pi = \pi_1 * r_\alpha(\pi_2) * \pi_3$), where the path $\pi$ is  cut
into well-defined parts
$\pi =  \pi_1 * \pi_2 * \pi_3$. Further,   $e_\alpha\pi(1) = \pi(1) +
\alpha^\vee$ and 
$f_\alpha\pi(1) = \pi(1) - \alpha^\vee$. After a suitable reparametrization, 
$e_\alpha\pi$ and $f_\alpha\pi$ are $\lambda-$paths in $Y$. More importantly,  
Littelmann obtains a  characterization of LS paths by using these operators. He
proves \cite[5.6]{Littelmann94} that  a $\lambda-$path $\pi$ with $\pi(0) = 0$ is a
LS path if, and only if, there exist  some simple roots
$\alpha_{i_1},...,\alpha_{i_s}$ such that 
$$
e_{\alpha_{i_1}}\circ\cdots\circ e_{\alpha_{i_s}} (\pi) = \pi_\lambda,
$$ where for all $t\in [0,1]$, $\pi_\lambda(t) = t\lambda$.

\begin{lemm}
\label{le:Operator}

i) If $\pi$ is a Hecke path in $Y$ and 
$e_\alpha\pi$ (resp. 
$\tilde e_\alpha\pi$) is defined, then  ${\rm d\!\dim}(e_\alpha\pi)={\rm d\!\dim}(\pi)
-1$ and ${\rm codim}(e_\alpha\pi)={\rm codim}(\pi) -1$ (resp. ${\rm d\!\dim}(\tilde
e_\alpha\pi)  ={\rm d\!\dim}(\pi) +1$ and ${\rm codim}(\tilde e_\alpha\pi) ={\rm
codim}(\pi)-1$), and, similarly, if $f_\alpha\pi$  is defined, then ${\rm
d\!\dim}(f_\alpha\pi)={\rm d\!\dim}(\pi)+1$ and ${\rm codim}(f_\alpha\pi)={\rm
codim}(\pi)+1$.

ii) If $\pi$ is a Hecke path in $Y$ such that $\tilde e_\alpha\pi$ is defined, 
then $\tilde e_\alpha\pi$ is again a Hecke path in $Y$.

iii) If $\pi$ is a Hecke path in $Y$ such that $\tilde e_\alpha\pi$ is not 
defined  but $e_\alpha\pi$ (resp. $f_\alpha\pi$) is, then $e_\alpha\pi$ (resp.
$f_\alpha\pi$)  is again a Hecke path in $Y$.
\end{lemm}

We prove the Lemma in Section \ref{ssse:ProofOperator}.

\subsubsection{Proof of Proposition \ref{pr:Inequalities}} 
\label{ssse:ProofInequalities}

By translation, we may (and shall often) suppose $\pi$ 
normalized, {\it i.e.} $\pi(0)=0$. It is clear that  ${\rm
d\!\dim}(\pi_\lambda)={\rm codim}(\pi_\lambda)=0$. As a corollary of i) and the 
characterization  of LS paths, if $\pi$ is a LS $\lambda-$path then ${\rm
d\!\dim}(\pi)={\rm codim}(\pi)=\rho(\lambda-\nu)$. The other implication is
obtained by induction on $\rho(\lambda -\nu)$. We suppose $\pi(0)=0$.  There is
only one $\lambda-$path $\pi$ such that $\pi(1) =
\lambda$; it  is $\pi_\lambda$. And in this case, ${\rm d\!\dim}(\pi_\lambda) = 0$. 

 If $\nu \ne \lambda$, then $w_+(0)\ne {\rm id}$ and there exists a simple root 
$\alpha$  such that $e_\alpha\pi$ or $\tilde e_\alpha \pi$ is defined. If, for
all 
$\beta$ simple, $\tilde e_\beta \pi$ is not defined, then the claim follows 
immediately by induction and by Lemma~\ref{le:Operator}. Otherwise, we apply all possible 
operators $\tilde e_\beta$ to $\pi$ to end up with a Hecke path $\eta$ such 
that $\eta(1) = \pi(1) = \nu$, ${\rm d\!\dim}(\eta) = {\rm d\!\dim}(\pi) + k$,  ${\rm
codim}(\eta) = {\rm codim}(\pi) - k$ ($k> 0$)  and there still exists $\alpha$
such that $e_\alpha\eta$ is defined. But then, by induction, 
${\rm d\!\dim}(\eta) -1 = {\rm d\!\dim}(e_\alpha\eta) \leq  \rho(\lambda -
e_\alpha\eta(1)) = \rho(\lambda - \nu) - 1$, which implies that  ${\rm
d\!\dim}(\pi)< \rho(\lambda - \nu)$. Moreover, ${\rm codim}(\pi)+{\rm
d\!\dim}(\pi)={\rm codim}(\eta)+{\rm d\!\dim}(\eta) ={\rm codim}(e_\alpha\eta)+{\rm
d\!\dim}(e_\alpha\eta)+2= 2\rho(\lambda - e_\alpha\eta(1))+2= 2\rho(\lambda - \nu)$
(by induction).

Suppose now that ${\rm d\!\dim}(\pi) = \rho(\lambda - \nu)> 0$, then for  dimension
reasons, $\tilde e_\alpha\pi$ is never defined. But $e_\alpha\pi$ is and  ${\rm
d\!\dim}(e_\alpha\pi) = \rho(\lambda - \nu) -1$. Repeating the same argument 
leads to a sequence of  simple roots $\alpha_{i_1},...,\alpha_{i_s}$ such that 
$e_{\alpha_{i_1}}\circ\cdots\circ e_{\alpha_{i_s}} (\pi) = \pi_\lambda$, 
in other words, $\pi$ is a LS path. This proves the proposition.\qed

\begin{rema} 
This proof implies also that a Hecke path in $Y$ is LS  if
and only if, for all simple roots $\alpha_j,\alpha_{i_1},...,\alpha_{i_s}$, the
minimum of 
$\alpha_j(e_{\alpha_{i_1}}\circ\cdots\circ e_{\alpha_{i_s}} (\pi) )$ is in $\mathbb Z$, 
cf. \cite[4.5]{Littelmann95}.
\end{rema}

\subsubsection{Proof of Lemma \ref{le:Operator}}
\label{ssse:ProofOperator}

We suppose $\pi(0)=0$. Let us start  with
the  operator $e_\alpha$ and dual dimensions. The paths $\pi$ and $e_\alpha\pi$
are cut into three parts, meaning that 
$$
\begin{cases}
\pi(t) =\pi_1(t), e_\alpha\pi(t) = \pi_1(t)& \hbox{ if } 0\leq t\leq 1/3\\
\pi(t) =\pi_2(t)+\pi(1/3), e_\alpha\pi(t) = r_\alpha(\pi_2)(t)+\pi_1(1/3) & \hbox{ if } 1/3\leq t\leq 2/3\\
\pi(t) =\pi_3(t)+\pi(2/3), e_\alpha\pi(t) = \pi_3(t)+e_\alpha\pi(2/3) & \hbox{ if } 2/3\leq t\leq 1.
\end{cases}
$$

For the first part of $\pi$, that is for $t\leq 1/3$, there is nothing  to
prove. Because $\alpha$  is a simple root, if $1/3< t \leq 2/3$, the relative
position of $\pi_-(t)$ with respect  to a wall $M(\beta,k)$ (with $\beta
\ne\alpha$) is the same as the relative position of 
$e_\alpha\pi_-(t)$ with respect to $r_\alpha M(\beta,k)$. Further, if  
$2/3<t\leq 1$, 
$e_\alpha\pi (t) = \pi(t) + \alpha^\vee$. So, again, up to translation  the
relative positions are the  same. It remains to check the positions relatively
to the walls $M(\alpha,-Q),  M(\alpha,-Q-2)$ at $t=2/3$. But $2/3$ is the
smallest real number $t$ such that 
$\pi(t)\in M(\alpha,-Q)$, therefore $\pi_-(2/3) \not\subset D(-\alpha,Q)$,
$e_\alpha\pi(2/3)\in M(-\alpha,Q+2)$
 and $e_\alpha\pi_-(2/3)\subset D(-\alpha,Q+2)$. 
Therefore, ${\rm d\!\dim} (e_\alpha\pi) = {\rm d\!\dim}(\pi) - 1$.

For the formulas ${\rm d\!\dim} (\tilde e_\alpha\pi) = {\rm d\!\dim} (\pi) +1$ and 
${\rm d\!\dim} (f_\alpha\pi) = {\rm d\!\dim} (\pi) +1$, similar arguments show that it 
suffices to look at the case  $t=2/3$ in the corresponding cuts of the  path
$\pi$. For  the operator $\tilde e_\alpha$, one has $\pi(2/3)=\tilde
e_\alpha\pi(2/3)\in M(-\alpha,Q)$ and $\pi_-(2/3)\subset D(-\alpha,Q)$ whereas  
$\tilde e_\alpha\pi_-(2/3)\not\subset D(-\alpha,Q)$. This proves the formula 
for the  operator $\tilde e_\alpha$. And for the operator $f_\alpha$ one has:
$\pi(2/3)\in M(-\alpha,Q+1)$,  $\pi_-(2/3)\subset D(-\alpha,Q+1)$ whereas
$f_\alpha\pi(2/3)\in D(-\alpha,Q-1)$
 and $f_\alpha\pi_-(2/3)\not\subset D(-\alpha,Q-1)$. The proof of i) for the 
dual dimensions is then complete.

For the codimensions, similar arguments show that it suffices to look at the case
$t=1/3$ and the root $\alpha$. For the operator $e_\alpha$, 
$e_\alpha\pi(1/3)=\pi(1/3)\in M(\alpha,-Q-1)$, $\pi_+(1/3)\not\subset D(\alpha,-Q-1)$, 
$e_\alpha\pi_+(1/3)\not\subset D(-\alpha,Q+1)$; therefore  ${\rm
codim}(e_\alpha\pi) = {\rm codim}(\pi) - 1$. For the operator $f_\alpha$ (resp.
$\tilde e_\alpha$), $f_\alpha\pi(1/3)=\pi(1/3)\in M(\alpha,-Q)$ (resp.  $\tilde
e_\alpha\pi(1/3)=\pi(1/3)\in M(\alpha,-Q)$ ), $\pi_+(1/3)\not\subset D(-\alpha,Q)$ and
$f_\alpha\pi_+(1/3)\not\subset D(\alpha,-Q)$ (resp. $\pi_+(1/3)\not\subset D(\alpha,-Q)$
and $f_\alpha\pi_+(1/3)\not\subset D(-\alpha,Q)$ ), therefore ${\rm
codim}(f_\alpha\pi) = {\rm codim}(\pi) + 1$ (resp. ${\rm codim}(\tilde
e_\alpha\pi) = {\rm codim}(\pi) - 1$ ).

Concerning ii), using the same arguments again, one has to take only care of the 
places $t=1/3$ and $t=2/3$ in the path $\pi$.  For $t=1/3$, the 
$W^v_{\pi(1/3)}-$chain for 
 $\tilde e_\alpha\pi$ is obtained from the one for $\pi$ just by adding 
 $\xi_{s+1}=\tilde e_\alpha\pi'_+(1/3)=r_\alpha(\pi'_+(1/3))$ and $\beta_{s+1}=\alpha$; 
as  $\alpha(\pi'_+(1/3))<0$ the conditions are satisfied. For $t=2/3$, the
$W^v_{\pi(2/3)}-$chain for 
 $\tilde e_\alpha\pi$ is obtained from the one for $\pi$ just by adding 
 $\xi_{-1}=\tilde e_\alpha\pi'_-(2/3)=r_\alpha(\pi'_-(2/3))$ and $\beta_{0}=\alpha$; as 
 $\alpha(\pi'_-(2/3))>0$ the conditions are satisfied (after a shift of the  indices
of the chain).  Therefore, $\tilde e_\alpha\pi$ is a Hecke path and ii) is
proved.

It remains to prove iii). Let us start with $e_\alpha$. Once again, it suffices  to
check the values $t=1/3$ and
$t=2/3$.  The situation around the point $\pi(1/3)$ is the same as above. 
Because $\tilde e_\alpha \pi$ is not defined, $\alpha(\pi'_+(2/3))\geq~0$. Let 
$(\xi_0,\dots,\xi_s)$, $(\beta_1,\dots,\beta_s)$ be the $W^v_{\pi(2/3)}-$chain 
from $\pi'_-(2/3)$ to $\pi'_+(2/3)$. If $\alpha=\beta_u$, $1\leq u\leq s$ (and $u$ 
is minimal for this property), then
\begin{gather*}
(r_\alpha\xi_0,r_\alpha\xi_1,\dots,r_\alpha\xi_{u-1}=\xi_u,\xi_{u+1},\dots,\xi_s),\\
(r_\alpha\beta_1,\dots,r_\alpha\beta_{u-1}, \beta_{u+1},\dots,\beta_s)
\end{gather*}
is the 
$W^v_{e_\alpha\pi(2/3)}-$chain  from $e_\alpha\pi'_-(2/3)$ to $e_\alpha\pi'_+(2/3)$.  If
no such $u$ exists and $\alpha(\pi'_+(2/3))> 0$ (resp. $\alpha(\pi'_+(2/3))= 0$), 
then this chain is $(r_\alpha\xi_0, $ $\dots, $ $r_\alpha\xi_s, \xi_{s+1} = \xi_s)$, $(r_\alpha\beta_1,\dots,r_\alpha\beta_s,\beta_{s+1}=\alpha)$  
(resp.  $(r_\alpha\xi_0,\dots,r_\alpha\xi_s=\xi_s)$, 
$(r_\alpha\beta_1,\dots,r_\alpha\beta_s)$ ). This proves that $e_\alpha\pi$ is still  a
Hecke path. The proof for $f_\alpha\pi$  follows similar lines
and is left to the reader! \qed

\section{Segments in the hovel} 
\label{se:Segments}

 This section contains the most important application of the definition of 
the hovel $\mathcal I$. We first prove that the retraction of any segment $[x,y]$ (with
$x\leq y$) in $\mathcal I$ is a Hecke path in $\mathbb A$ (see Theorem~\ref{th:RetractingSeg}). Then, we give a parametrization of all segments retracting on a given Hecke path sharing the
same end (Theorem~\ref{th:RetractingOnHecke} and Corollary~\ref{co:NumberPara}). The algebraic structure of the set of parameters is studied in \ref{sse:Algebraic} and allows us to define a generalization of Mirkovi\'c-Vilonen cycles. Then, we state another characterization of LS paths
in terms of a new statistic, but depending on extra data and not solely on
the path (\ref{sse:AnotherCharLS}). To finish, we prove a result on the structure of $\mathcal I$ (Theorem~\ref{th:Preorder}).

The field $K$ is as in Section \ref{sse:ActionNK}. Note however that, in the classical case
where $G$ is a split reductive group, all what follows holds for any field  $K$
endowed with a discrete valuation; we just have to use the Bruhat-Tits building
instead of the hovel constructed in Section \ref{se:Hovel}.

\subsection{Retracting segments}
\label{sse:RetractingSeg}

We consider a negative sector germ $\mathfrak S$ and 
denote by $\rho$ the retraction of center 
$\mathfrak S$ without specifying on which apartment $A$ (containing $\mathfrak S$) $\rho$ maps~
$\mathcal I$, as $\rho$  does not depend on $A$ up to canonical isomorphisms. Actually
we identify any pair $(A,\mathfrak S)$  of an apartment $A$ containing $\mathfrak S$ with the
fundamental pair $(A_f=\mathbb A,\mathfrak S_{-\infty})$,  this is well determined up to
translation.

 We consider two points $x,y$ in the hovel with $x\leq y$. The segment 
$[x,y]$  is the image of the path $\pi:[0,1]\longrightarrow \mathcal I$ defined by
$\pi(t)=x+t(y-x)$ in  any apartment containing $x$ and $y$ (\ref{sse:Applications}.2). As each
segment in $[x,y]$ has a  good fixator, the derivative $\pi'(t)$ is independent
of the apartment containing a  neighbourhood of $\pi(t)$ in $[x,\pi(t)]$ or
$[\pi(t),y]$, up to the Weyl group $W^v$.

 We saw in \ref{sse:Retraction} that the image $\rho\pi$ is a piecewise linear continuous 
path  in~$A$. By the previous paragraph, there exists a $\lambda$ in the 
fundamental closed-chamber $\overline{C^v_f}$ such that $\rho\pi'(t)=w_t.\lambda$ 
for each $t\in[0,1]$ (different from $\pi^{-1}(x_i)$ for $x_i$ as in \ref{sse:Retraction}) 
and some $w_t\in W^v$ (chosen minimal with this  property). Hence $\rho\pi$ is a
$\lambda-$path (in particular the map $\rho$ is increasing with respect to the
``preorder'' of \ref{sse:Applications}.2)  and may be described as 
$\pi(\lambda,\rho\pi(0),\underline \tau,\underline a)$. We shall prove that $\rho\pi$ 
is a  Hecke path and often a LS path.

We choose some $t\in]0,1[=[0,1]\setminus\{0,1\}$ and we set 
$z = \rho\pi(t)$.  We denote by 
$\rho\pi'_-$ (resp. $\rho\pi'_+$) the left (resp. right) derivative of $\rho\pi$ 
in $t$  and $w_-$ (resp. $w_+$) the minimal element in $W^v$ such that
$\rho\pi'_-=w_-\lambda$  (resp. $\rho\pi'_+=w_+\lambda$).

\begin{prop}
\label{pr:RetractingSeg}
We have $\rho\pi_+'\leq_{W^v_z}\rho\pi_-'$ (cf. Definition \ref{de:Hecke})
and $w_+\leq w_-$ in the Bruhat-Chevalley order of $W^v/W^v_\lambda$. More 
precisely, there exist $s\in\mathbb N$ and a sequence $\beta_1,\dots,\beta_s$
 of positive real roots such that:

 -- for $1\leq i\leq s$, there exists a wall of direction ${\rm Ker}(\beta_i)$ 
containing $z=\rho\pi(t)$,

 -- if one defines $\xi_0=\rho\pi'_-$, $\xi_1=r_{\beta_1}.\xi_0$, \dots, 
$\xi_s=r_{\beta_s}.\dots.r_{\beta_1}.\xi_0$, one has $\xi_s=\rho\pi'_+$ and 
$\beta_k(\xi_{k-1})<0$ for $1\leq k\leq s$,

 -- if one defines $\sigma_0=w_-$, $\sigma_1=r_{\beta_1}.w_-$, \dots, 
$\sigma_s=r_{\beta_s}\cdots r_{\beta_1}.w_-$, then $\rho\pi'_+=\sigma_s\lambda$ and, for 
$1\leq k\leq s$,  one has 
$\sigma_k<\sigma_{k-1}$ in the Bruhat-Chevalley order of $W^v/W^v_\lambda$,

 -- there exists in $\mathfrak  a_+=A\cap\mathcal I^+_z$ an (unrestricted) gallery
$\delta=(c_0,c_1',\dots,c_n')$ from
$c_0={\rm germ}_z(z+C^v_f)$  to $c_n'={\rm germ}_z(z+w_+C^v_f)\supset z+[0,1)\rho\pi_+'$, 
the type of which is associated  to a (given) reduced decomposition of $w_-$.
The panels along which this gallery is folded  (actually, positively folded:
see the proof) are successively the walls $z+{\rm Ker}(\beta_1)$, \dots ,  
$z+{\rm Ker}(\beta_s)$.
\end{prop}

\begin{proof}
Let $A^0$ be an apartment containing $[x,y]$; set  
$\pi_-=[\pi(t),x)$  and $\pi_+=[\pi(t),y)$. By \ref{sse:Applications}.3 and \ref{sse:Retraction} there exist
apartments $A^+$ and 
$A^-$ containing the sector $\mathfrak s$ (of direction $\mathfrak S$ and base point 
$\pi(t)$) and  respectively $\pi_+$ and $\pi_-$. We choose $A^-$ for the image
$A$ of $\rho$, so 
$\pi(t)=\rho\pi(t)=z$ and $\pi_-=\rho\pi_-=z-[0,1)\pi'_-$.

 As $A^0$ and $A^-$ contain $\pi_-$, there exists $g\in \hat P_{\pi_-}$
 such that $A^0=g.A^-$. In the decomposition 
$\hat P_{\pi_-}=U^{nm-}_{\pi_-}.U^{pm+}_{\pi_-}.\hat N_{\pi_-}$, the 
group 
$\hat N_{\pi_-}$ fixes $\pi_-$ and stabilizes $A^-$, so, one has 
$A^0=u^-u^+A^-$ with $u^-\in U^{nm-}_{\pi_-}$ and $u^+\in U^{pm+}_{\pi_-}$.
 Let us consider the apartment $A^1=(u^-)^{-1}A^0=u^+A^-$; it contains 
 $\pi_-=(u^-)^{-1}(\pi_-)$ and $\pi_+^1=(u^-)^{-1}(\pi_+)$, which are 
 opposite segment germs; moreover $\rho(\pi_+^1)=\rho(\pi_+)$. On the 
 other side, $A^1$ contains the chamber $C_0=F(z,C^v_f)$ in $A^-$, 
 which is opposite $\mathfrak s$. We replace in the following $\pi_+$ by 
 $\pi^1_+$ and $A^0$ by $A^1$.

 In $A^-$, $\pi_-\in w_-\mathfrak s$, hence, the germ opposite $\pi_-$ is in 
$w_-C_0$. So, in $A^1$, $\pi_+^1\in w_-^1C_0$, where  $w_-^1$ corresponds to 
$w_-$ in the identification of $W^v(A^-)$ and  $W^v(A^1)$ via $u^+$.

 We choose in $\mathfrak  a^1_+=A^1\cap \mathcal I_z^+$ (an apartment of $\mathcal I_z^+$ 
endowed  with the unrestricted building structure) a minimal gallery $m = (c_0,
c_1,..., c_n)$
 between  $c_0 = {\rm germ}_z(z+C_f^v)$ and $c_n \supset \pi_+^1=z+[0,1)w_-^1\lambda$ 
of type $\tau = (i_1,...,i_n)$, $i_j\in I$; hence $w_-^1=r_{i_1}\cdots 
r_{i_n}$ is a  reduced decomposition.  The restriction of $\rho$ to the residue 
twin building $(\mathcal I_z^+,\mathcal I_z^-, d_z^*)$ (with the unrestricted structure) 
preserves the codistance to $\overline{\mathfrak s}={\rm germ}_z(\mathfrak s)$, 
which is a chamber in $\mathcal I_z^-$. Therefore, this restriction is the retraction 
$\rho_z : \mathcal I_z^+\to \mathfrak  a^-_+=A^-\cap \mathcal I_z^+$ of centre $\overline{\mathfrak s}$. 
We have $\rho\pi_+^1 = \rho_z\pi_+^1$.

 The retracted gallery $\delta=\rho_z(m) = (c_0,c_1'=\rho_z(c_1),...,c_n'=
\rho_z(c_n))$ in $A=A^-$ is a 
positively folded gallery, meaning that $\rho_z(c_j) = \rho_z(c_{j+1})$ implies 
that $\rho_z(c_j)$ is on the positive side of the wall $H_j$ spanned by the panel 
of type $\{i_j\}$ of $\rho_z(c_j)$ (note that $H_j$ is a wall for the 
unrestricted structure).  Otherwise, suppose that $\rho_z(c_j) = \rho_z(c_{j+1})$
is on the negative side of $H_j$. Then,  because $\overline {\mathfrak s}$ is the
opposite fundamental chamber in $z$, it is always on  the negative side of
$H_j$. Further, let $\mathfrak  a$ be a twin apartment containing 
$\overline {\mathfrak s}$ and $c_j$, as the retraction preserves the codistance to 
$\overline {\mathfrak s}$, we also have  that $\overline {\mathfrak s}$ and $c_j$ are on the
same side of the wall spanned by the panel  of type $\{i_j\}$ of $c_j$ in $\mathfrak 
a$. Therefore, \ref{sse:Applications}.4) implies that, modifying the latter  if needed, we can
assume that $c_{j+1}$ is still in $\mathfrak  a$. But, on one side, 
$c_j \not = c_{j+1}$ then, computing in $\mathfrak  a$, 
$\ell(d_z^*( \overline {\mathfrak s}, c_{j+1})) =  \ell(d_z^*( \overline {\mathfrak s}, 
c_j)) -1$;  on the other side, $\ell(d_z^*( \overline {\mathfrak s}, c_{j+1})) = 
\ell(d_z^*( \overline {\mathfrak s}, \rho_zc_{j+1})) = \ell(d_z^*( \overline {\mathfrak s}, 
\rho_zc_{j}))  = \ell(d_z^*( \overline {\mathfrak s}, c_{j}))$. Contradiction!

 If the wall $H^1_j$ separating $c_j$ from $c_{j+1}$ in $\mathfrak  a^1_+$ is a 
ghost   wall {\it i.e.} not a true wall (for the restricted structure), then the
enclosure of $c_j$  in the hovel contains $c_{j+1}$ and there is an apartment
of $\mathcal I$ containing $\mathfrak s$, $c_j$  and $c_{j+1}$, so $\rho_z(c_j) \neq
\rho_z(c_{j+1})$.

 Let us now denote by $j_1,...,j_s$ the indices such that 
$c'_j=\rho_z(c_{j}) = \rho_z(c_{j+1})=c'_{j+1}$.  For any $k\in\{1,...,s\}$, 
$H_{j_k}$ is a true wall spanned by the panel $\rho(H^1_{j_k}\cap c_j)$ and we
denote the positive real root associated  with $H_{j_k}$ by $\beta_k$ (i.e.
$H_{j_k}$ is of direction ${\rm Ker}(\beta_k)$). Actually,  the gallery $\delta$ is
obtained from the minimal gallery $\delta^0 = (c_0^0=c_0,c_1^0,\dots,c_n^0)$  of
type $(i_1,...,i_n)$ beginning in $c_0$, ending in
$c_n^0=w_- (c_0)$  and staying inside $A^-$ by applying successive (positive) 
foldings along the walls associated to the indices $\{j_1,...,j_s\}$, starting 
with $H_{j_1}$,  then folding along $H_{j_2}$... At each step, one gets a
positively folded gallery 
$\delta^k = (c_0^k=c_0,c_1^k,\dots,c_n^k)$ ending closer and closer to the 
chamber $c_0$. So, this proves the last assertion of the proposition.

 Let us denote $\xi_0=\pi_-'=\rho\pi_-'=w_-\lambda$ (in $A_-$) and 
$\xi_k=r_{\beta_k}\cdots r_{\beta_1}\xi_0$. As $\rho_z(c_n)\supset 
\rho\pi_+=z+[0,1)\rho\pi_+'$  and $\rho\pi_+'\in W^v\rho\pi_-'$ one has
$\xi_s=\rho\pi_+'$, and more generally, 
$z+[0,1)\xi_k\subset c_n^k$. As $\delta^0$ is a minimal gallery from $c_0$ to 
$z+[0,1)\pi'_-$, $c^0_{j+1},\dots,c^0_{n}$ and $z+[0,1)\pi'_-$ are on the  same
side of any wall separating $c^0_0$ from $c^0_{j+1}$; in particular, 
$(c^k_{j_k+1},\dots,c^k_{n})$ is a minimal gallery, entirely on the same side 
of $H_{j_k}$ and $z+[0,1)\xi_k\not\subset H_{j_k}$. But 
$c^k_{j_k}=\rho_z(c_{j_k})=\rho_z(c_{j_k+1})=c^k_{j_k+1}$ and we saw  that this
chamber is on the positive side of the wall $H_{j_k}$ (of direction 
Ker$(\beta_{k})$). So,
$c^k_{j_k+1},\dots,c^k_{n}$ are on the positive side of $H_{j_k}$; this means
that $\beta_k(\xi_k)>0$ {\it i.e.} $\beta_k(\xi_{k-1})<0$. Hence, the sequences 
$(\xi_0,\xi_1,\dots,\xi_s)$ and $(\beta_1,\dots,\beta_s)$ give a $W^v_z-$chain from 
$\pi_-'=\rho\pi_-'$ to $\rho\pi_+'$. This proves the proposition, in view of 
Lemma~\ref{le:Equiv1}.
\end{proof}

\begin{theo} 
\label{th:RetractingSeg}

Let $\pi=[x,y]$ be a segment in an apartment 
$A'$ with $x<y$, and $\rho$ the retraction of $\mathcal I$ with center the fundamental 
sector-germ $\mathfrak S_{-\infty}$ onto an apartment $A$. Then the retracted segment
$\rho\pi$ is a Hecke path in $A$.

 If moreover $x$ and $y$ are cocharacter points ({\it i.e.} $x,y\in Y(A')$),  then
$\rho\pi$ is a Hecke path in $Y(A)$. 
\end{theo}

\begin{proof}
The path $\rho\pi$ is Hecke by Proposition~\ref{pr:RetractingSeg} and Definition~\ref{de:Hecke}. If $x,y\in Y(A')$, computing in $A'$, $\lambda=W^v(y-x)\cap\overline{C^v_f}$
is in $Y^+(A')$. Moreover, by \ref{sse:Applications}.1, $\rho(Y(A'))\subset Y(A)$, so $\rho\pi$ is 
in $Y(A)$. 
\end{proof}
 
\subsection{Segments retracting on a given Hecke path}
\label{sse:RetractingOnHecke}

 If we consider all segments $\pi=[x,y]$ from some $x\in\mathbb A$ to some 
$y\in\mathcal I$  ({\it i.e.} $\pi(t)=x+t(y-x)$) whose retraction $\rho\pi$ is a given Hecke
path $\pi_1$ in $A$ starting at $x$, then, there are too many of them. For
example, take for $\pi_1$ the path $t\mapsto t\lambda$ (with $\lambda$ dominant), then,
already at $0$, one has infinitely many choices to define a segment starting at
$0$ and retracting onto $\pi_1$. Therefore, in this subsection, we fix
$y=\pi(1)$.

 More precisely, let $y\in\mathcal I$ and $\pi_1$ a Hecke path in $A$ with 
$\pi_1(1)=\rho(y)$,  we define $S(\pi_1,y)$ as the set of all segments 
$\pi=[x,y]$ in $\mathcal I$ such that $\pi_1=\rho\pi$.
 
\begin{theo}
\label{th:RetractingOnHecke}

The set $S(\pi_1,y)$ is nonempty and is parametrized by 
exactly $N={\rm d\!\dim}(\pi_1)$ parameters in the residue field $\kappa$. More precisely, the set $P(\pi_1,y)$ 
of parameters is a finite union of subsets of $\kappa^N$, each being a product of
$N$ factors either equal to $\kappa$ or to $\kappa^*$.
\end{theo}

\begin{rema} 
\label{re:QuasiAffine}
In particular, $P(\pi_1,y)$ is a Zariski open subset of 
$\kappa^N$ stable under the natural action of the torus $(\kappa^*)^N$, in other
words, a quasi-affine toric variety.
\end{rema}

\begin{proof}
 We shall prove that, for $t\in[0,1]$, the segments 
$\pi^t:[t,1]\rightarrow\mathcal I$ with $\pi^t(1)=y$, retracting onto 
$\pi_1^t=\pi_1{}_{\vert_{[t,1]}}$, are parametrized by exactly 
$\sum_{t'>t}\ell_{\pi_1(t)}(w_-(t))$ parameters. This is clear for $t=1$.
Suppose the result true for some $t$. So, $\pi(t)$ is  given and we shall prove
now that the number of parameters for the choice of the  segment-germ
$\pi_-(t)$ of origin $\pi(t)$ is $\ell_{\pi_1(t)}(w_-(t))$. This result and
arguments after the Definition~\ref{de:DdimCodim} imply that $\pi_-(t)$ determines $\pi{}_{\vert_{[t',t]}}$
where $t'\;(<t)$  is $0$ or the next number in $[0,1]$ such that 
$\overline{\pi_1}$ leaves positively a wall in $t'$.

 Now, $\pi(t)$ is given and we want to find out how many parameters govern 
the choice of $\pi_-(t)$. We choose $A$ so that it contains $\mathfrak S$ and
$\pi_+(t)$; so $\pi_+(t)=\rho\pi_+(t)$  and we set $z=\pi(t)$. There is an
apartment $A^-$
 containing $\mathfrak S$ and $\pi_-(t)$, thus, $\rho$ is an isomorphism from $A^-$ to 
$A$. The  (unrestricted) chamber $c_0^-={\rm germ}_z(z-C^v_f)$ (cf. \ref{sse:Residue}) is in $A\cap
A^-$. We  choose a reduced decomposition $w_-(t)=r_{i_1}.\dots.r_{i_n}$ in
$W^v$. The associated  minimal (unrestricted) gallery of type $(i_1,\dots,i_n)$ 
from $c_0^-$ to $\pi_-(t)$ is denoted by $m^-=(c_0^-,c_1^-,\dots,c_n^-)$. 
Clearly, $\pi_-(t)$ is entirely determined by the gallery 
$m^-$, so it seems to depend on $n=\ell(w_-(t))$ parameters in $\kappa$. But 
actually, if  the wall separating $c^-_{j-1}$ from $c^-_{j}$ (or $\rho
c^-_{j-1}$ from $\rho c^-_{j}$)  is a ghost wall {\it i.e.} not a true wall (see Section~\ref{sse:Residue}), the
chamber $c^-_{j}$ is  determined by $c^-_{j-1}$; whereas, if this wall is true
(some $M(\alpha,k)$ for some 
$\alpha\in\Phi^-$), then the choice of $c^-_{j}$ depends on an element in 
$U^\times_{\alpha,k} = U_{\alpha,k}/U_{\alpha,>k}\simeq \kappa$ (cf. \ref{sse:Applications}.4). 
Hence, the true  number of parameters is $\ell_{\pi_1(t)}(w_-(t))$.

 But, we forgot to check that $\pi_-(t)$ is opposite $\pi_+(t)$. 
Actually, as we shall see now, removing at most one value for each parameter,
this condition is fulfilled. This proves the first part of the theorem.

 As $\pi_1$ is a Hecke path, conditions i), v), vi) and vii) of Definitions~\ref{de:LS} and \ref{de:Hecke} are satisfied for some roots
$\beta_1$,\dots,$\beta_s$ and we can use the results of Section~\ref{ssse:LSvHecke}. Let $\delta^0=(c_0,c_1^0,\dots,c_n^0)$  be the minimal
gallery of type $(i_1,\dots,i_n)$ in $\mathfrak  a_+=A\cap\mathcal I^+_x$ starting from
$c_0={\rm germ}_z(z+C^v_f)$;  its end $c_n^0=w_-c_0$ contains $z+[0,1)\pi_-'(t)$. We 
shall fold this gallery stepwise.
 Since $r_{\beta_1}w_-<w_-$, the wall $z+{\rm Ker}(\beta_1)$ separates $c_0$ from 
$c_n^0$: it is  the wall between some adjacent chambers $c^0_{j-1}$ and
$c^0_{j}$. We define 
$\delta^1=(c_0,c_1^1=c_1^0,\dots,c^1_{j-1}=c^0_{j-1},
c^1_{j}=r_{\beta_1}c^0_{j},\dots,c^1_{n}=r_{\beta_1}c^0_{n})$, so
$c^1_{j-1}=c^1_{j}$ and $c^1_{n}=r_{\beta_1}w_-c^0$.
 But $r_{\beta_2}r_{\beta_1}w_-<r_{\beta_1}w_-$, so the wall $z+{\rm Ker}(\beta_2)$ 
separates $c_0$ from $c_n^1$: it is the wall between some adjacent chambers 
$c^1_{k-1}$ and $c^1_{k}$.
 We define $\delta^2=(c_0,c_1^2=c_1^1,\dots,c^2_{k-1}=c^1_{k-1},
c^2_{k}=r_{\beta_2}c^1_{k},\dots,c^2_{n}=r_{\beta_2}c^1_{n})$. At the end of this 
procedure, we  get a gallery $\delta^s=(c_0,c_1^s,\dots,c_n^s)$ of type
$(i_1,\dots,i_n)$ in $A$  starting from $c_0$ and ending in
$c_n^s=w_+(t)c_n^0\supset\pi_+(t)$. Moreover, this  gallery is positively
folded along true walls.

 As $\pi_+'(t)\in W^v\pi_-'(t)$, to prove that $\pi_-(t)=z-[0,1)\pi_-'(t)$ 
and $\pi_+(t)=z+[0,1)\pi_+'(t)$ are opposite segment-germs, it suffices to 
prove that $c_n^s$ and $c_n^-$ are opposite chambers. For this, we prove that, 
except perhaps for one  choice of each parameter, $c_j^s$ and $c_j^-$ are
opposite for $0\leq j\leq n$. This is true   for $j=0$. Suppose $c_{j-1}^-$
opposite $c_{j-1}^s$. Then $c_{j}^-$ (resp.
$c_{j}^s$)  is adjacent to $c_{j-1}^-$ (resp. $c_{j-1}^s$) along an 
(unrestricted) panel of type 
$i_j$. If the wall containing these two panels is not true ({\it i.e.} restricted), 
then $c_{j}^-$ and $c_{j}^s$ are automatically opposite. Now, if this wall is 
true, by \ref{sse:Residue} and the general properties of twin buildings (see \cite[2.5.1]{Remy02}) among the  chambers adjacent (or equal) to $c_{j-1}^-$ along the panel
of type $i_j$, there is a unique chamber not opposite $c_{j}^s$. Hence, all 
but (perhaps) one choice for $c_{j}^-$ is opposite $c_{j}^s$; and the 
corresponding parameter has to be chosen in $\kappa$ or in $\kappa^*$. Therefore, the
set $S(\pi_1,y)$ is nonempty.

 Let us have a closer look at the set of
parameters. Choose $\pi\in S(\pi_1,y)$ and $t\in[0,1]$. We show now that 
$\pi_-(t)$ is obtained with the above procedure. We have the gallery
$(c_0^-,\dots,c_n^-)$ as above in
$A^-$. We choose the apartment $A$ containing $\mathfrak S$ and a chamber
$c_n\supset\pi_+(t)$ opposite $c_n^-$. Using the same properties of twin 
buildings, we find a gallery $\delta=(c_0,c_1,\dots,c_n)$ of type
$(i_1,\dots,i_n)$ in $A$, folded only along true walls, and such that, for all
$j$, $c_j$ and $c_j^-$ are opposite. In particular, $c_0$ is as defined above.
So, using $\delta$ instead of $\delta^s$, $\pi_-(t)$ is defined as before. Moreover,
the number of possibilities for $\delta$ is finite. Hence, the set of parameters
for $S(\pi_1,y)$ is a finite union of subsets of
$\kappa^N$, each being a product of $N$ factors either $\kappa$ or $\kappa^*$.
\end{proof}

\begin{coro}
\label{co:NumberPara}
Suppose $\pi_1$ is a Hecke path in $Y(A)$. 
Then the number ${\rm d\!\dim}(\pi_1)$ of  parameters for $S(\pi_1,y)$ is at most
$\rho_{\Phi^+}(\lambda-\pi_1(1)+\pi_1(0))$, with equality  if and only if $\pi_1$ is 
a LS path.
\end{coro}
 
\begin{proof}
This is a simple consequence of Proposition~\ref{pr:Inequalities} and  Theorem~\ref{th:RetractingOnHecke}.
\end{proof}

\subsection{Algebraic structure of $S(\pi_1,y)$ and Mirkovi\'c-Vilonen cycles}
\label{sse:Algebraic}

 To simplify notation, we suppose that $y=0$ in $A$ and (as before) 
$\lambda\in Y^+$. Moreover, we suppose that $K=\mathbb C(\!(\varpi)\!)$.
\begin{enumerate}
\item[1)] The set $\mathcal G_\lambda$ of segments $\pi$ in $\mathcal I$ of shape $\lambda$ and 
ending in $0$ may be identified with the set of its starting points $\pi(0)$
{\it i.e.} with $G_0.(-\lambda) = G(\mathcal O).(-\lambda)$. For $\nu\in Y$, let us define $\mathcal G_{\lambda,\nu}$ as the subset of $\mathcal G_\lambda$ consisting of the segments $\pi$ with
$\rho(\pi(0))=-\nu$. Thus, $\mathcal G_{\lambda,\nu}$ is identified with
$U^-(K).(-\nu)\cap G_0.(-\lambda)$. As $-\lambda\in G(K).0$, we can see $\mathcal G_{\lambda,\nu}$ as a subset of the affine grassmannian $\mathcal G=G(K)/G(\mathcal O)$, cf. Example~\ref{ex:G(O)}. We shall  view the algebraic structure of $\mathcal G_{\lambda,\nu}$ as inherited from
$U^-(K)$.

By Theorem~\ref{th:RetractingSeg} and Corollary~\ref{co:FiniteHecke}, $\mathcal G_{\lambda,\nu}$ is the finite (disjoint) union of the
subsets $S(\pi_1,0)$ for $\pi_1$ a Hecke path of shape $\lambda$ in $A$ from 
$-\nu$ to $0$.

\item[2)] Now, we better describe the parameters for $S(\pi_1,0)$ found in Theorem~\ref{th:RetractingOnHecke}. Let $0<t_1<\dots<t_m\leq 1$ be the values of $t$ such that
$n_i=\ell_{\pi_1(t_i)}(w_-(t_i))>0$ and $t_0=0$, $t_{m+1}=1$. For  $1\leq i\leq
m$ given, there exist negative roots 
$\alpha_{i,j}$ and integers $k_{i,j}$, $1\leq j\leq n_i$, such that
$M(\alpha_{i,n_i},k_{i,n_i})$, $\dots,$ $M(\alpha_{i,1},k_{i,1})$ are the true walls 
successively crossed by a minimal gallery from
$c_0^-={\rm germ}_{\pi_1(t_i)}(\pi_1(t_i)-C^v_f)$ to
$\pi_{1-}(t_i)$. Further, for any $a\in\mathbb C$, let us set 
$x_{i,j}(a)=x_{\alpha_{i,j}}(a\varpi^{k_{i,j}})\in U_{\alpha_{i,j},k_{i,j}}^\times$.
Moreover, let $\pi\in S(\pi_1,0)$ and $g\in U^-(K)$ such that
$\pi(t_i)=g\pi_1(t_i)$; since $g^{-1}\pi_-(t_i)$ is also the end of a minimal 
gallery from $c_0^-$ of the same type, for any $t\in [t_{i-1},t_i[$,
$$g^{-1}\pi(t) = x_{i,n_i}(a_{i,n_i})\cdots x_{i,1}(a_{i,1}) \pi_1(t)$$ 
for some
parameters $a_{i,n_i},...,a_{i,1}$ that have to be chosen in $\mathbb C$ or $\mathbb C^*$
according to the proof of Theorem~\ref{th:RetractingOnHecke}.

Iterating this procedure, one obtains that if $\pi\in S(\pi_1,0)$, then 
there exists some $(a_{i,j})\in P(\pi_1,0)\subset\mathbb C^N$ such that
$$\pi(0)=(\prod_{i\geq1;j\leq n_i}\;x_{i,j}(a_{i,j})).(-\nu)$$ 
where
the product is taken in lexicographical order from right to left. More
generally, for $t_{i_0-1}\leq t\leq t_{i_0}$, 
$$\pi(t)=(\prod_{i\geq i_0;j\leq n_i}\;x_{i,j}(a_{i,j})).(\pi_1(t)).$$

Thus, we define a map 
$$\mu:\mathbb C^N\supset P(\pi_1,0)\rightarrow U^-(K), 
(a_{i,j})_{i\leq m;j\leq n_i}\mapsto\prod_{i\leq m;j\leq n_i}\;x_{i,j}(a_{i,j})$$ 
such that the composition 
$$\overline{\mu}=\rm{proj}\circ\mu:\mathbb C^N\supset
P(\pi_1,0)\rightarrow U^-(K)/U^-(K)_{-\nu}$$ 
is injective. But, 
$$U^-(K)\subset U^{nmax-}(K)=\prod_{\alpha\in\Delta^-}U_\alpha(K)=\prod_{\alpha\in\Delta^-}\mathbb C(\!(\varpi)\!)$$ 
and, as $\mu$ involves finitely many groups $U_{\alpha,k}$ with $\alpha\in\Phi^-$, there
exists $y$ in $Y$ such that the image of $\mu$ is contained in $U^-(K)_y\subset
U^{nmax-}(K)_y = \prod_{\alpha\in\Delta^-}U_{\alpha,-\alpha(y)} = \prod_{\alpha\in\Delta^-}\mathbb C[[\varpi]]$. This last group has the structure of a pro-group in the sense of \cite{Kumar02} and
the map $\mu$ is clearly a morphism for this algebraic structure.

\item[3)] Hence, $\mathcal G_{\lambda,\nu}$ is a finite (disjoint!) union of sets 
$S(\pi_1,0)$ each in bijection with a quasi-affine irreducible variety
$P(\pi_1,0)$ and these sets are indexed by the Hecke paths $\pi_1$ of shape $\lambda$
from $-\nu$ to $0$ in $A$. The maximal dimension of these varieties is
$\rho_{\Phi^+}(\lambda-\nu)$, and the varieties of maximal dimension correspond to LS
paths from $-\nu$ to $0$ in $A$. A Mirkovi\'c-Vilonen cycle inside $\mathcal G_{\lambda,\nu}$ should be the closure of a set $S(\pi_1,0)$ (for $\pi_1$ a LS path)
and $P(\pi_1,0)$ should be isomorphic to a dense open subvariety of this cycle.

This holds in the classical case of reductive groups. These cycles in $\mathcal G_{\lambda,\nu}$ are dense in the Mirkovi\'c-Vilonen cycles corresponding to $-\lambda$ 
and $-\nu$ and described by using the reverses of the paths above, cf.
\cite{GaussentLittelmann05}. 
\end{enumerate}

\subsection{Another characterization of LS paths} 
\label{sse:AnotherCharLS}

 Suppose $\pi_1$ is a Hecke path of shape $\lambda$ in the apartment $A$. For each $t$,
$0<t<1$, let $w_-(t)$ (resp. $w_+(t)$ ) be the minimal element in $W^v$ such that
$\pi'_{1-}(t)=w_-(t)\lambda$ (resp. $\pi'_{1+}(t)=w_+(t)\lambda$ ). By Proposition~\ref{pr:RetractingSeg} and Theorem~\ref{th:RetractingOnHecke}, there 
exists an unrestricted gallery $\delta_t=(d_0,\dots,d_n)$ in $\mathfrak  a_+=A\cap\mathcal I^+_{\pi_1(t)}$, of type $(i_1,\dots,i_n)$ associated to a (given) reduced 
decomposition of $w_-(t)$,  starting from 
$d_0=c_0={\rm germ}_{\pi_1(t)}(\pi_1(t)+C^v_f)$ and ending in 
$d_n\supset\pi_1(t)+[0,1)\pi'_{1+}(t)$. Moreover, this gallery may be taken
positively folded along true walls. For each $t$, we choose such a gallery and we
set $\tilde\pi_1=(\pi_1,(\delta_t)_{0<t<1})$.

 The gallery $\delta_t$ is minimal for almost all $t$ (when 
$\pi'_{1-}(t)=\pi'_{1+}(t)$ ). Let us define $neg(\delta_t)$ as the number of all
unrestricted walls $H_j$ (containing the panel of type $i_j$ in $d_j$ or
$d_{j-1}$) which are true walls and separate $d_j$  from $d_0$. Actually, as
$\delta_t$ is positively folded, such an $H_j$ separates $d_j$ from
$d_{j-1}$ {\it i.e.} $d_j\neq d_{j-1}$.

\begin{defi} 
\label{de:CodimTilde}
The {\it codimension} of $\tilde\pi_1$ is:
$$ {\rm codim}(\tilde\pi_1)=\ell_{\pi(0)}(w_+(0))+\sum_{0<t<1}neg(\delta_t).$$
\end{defi}

 By the same arguments as for Definition~\ref{de:DdimCodim}, ${\rm codim}(\tilde\pi_1)$ is a nonnegative
integer; actually,\quad ${\rm codim}(\tilde\pi_1)\leq
\ell_{\pi(0)}(w_+(0))+{\rm d\!\dim}(\pi_1)-\ell_{\pi(1)}(w_-(1))$.

\begin{prop} 
\label{pr:InequalitiesTilde}
Let $\pi_1$ be a Hecke path in $Y$. For each 
choice of $\tilde\pi_1$, ${\rm codim}(\tilde\pi_1)\geq{\rm codim}(\pi_1)$. 
Further, $\pi_1$ is a LS path if and only if there is equality for (at least) one
choice of $\tilde\pi_1$.
\end{prop}
 
\begin{rema} 
\label{re:InequalitiesTilde}
Therefore, ${\rm codim}(\tilde\pi_1)\geq{\rm codim}(\pi_1)
\geq\rho_{\Phi^+}(\lambda-\nu)\geq{\rm d\!\dim}(\pi_1)$, with equalities if and only if 
$\pi_1$ is a LS-path (for good choices of~$\tilde\pi_1$).
\end{rema}

\begin{proof}
 It is clear that any true wall $H$ separating $d_0$ from
$\pi_{1+}(t)=\pi_1(t)+[0,1)\pi'_{1+}(t)$ is among the walls $H_j$, and, if $j$ 
is chosen maximal for this property, $H$ separates $d_0$ from $d_j$. So
$\ell_{\pi_1(t)}(w_+(t))\leq neg(\delta_t)$ for $0<t<1$ and ${\rm
codim}(\tilde\pi_1)\geq{\rm codim}(\pi_1)$. 

 Suppose $\ell_{\pi_1(t)}(w_+(t))=neg(\delta_t)$, then every true wall $H$ 
separating $d_0$ from $\pi_{1+}(t)$ is leaved negatively once and only once by the gallery 
$\delta_t$; in particular $\delta_t$ cannot be negatively folded along such a wall  and
cannot cross it positively. Moreover $\delta_t$ cannot leave negatively any other
true wall. As, by hypothesis, $\delta_t$ may only be folded along a true wall, this
gallery remains inside the (unrestricted) enclosure of $d_0$ and $\pi_{1+}(t)$.
The number of foldings of $\delta_t$ is 
$s=\ell_{\pi_1(t)}(w_-(t))-\ell_{\pi_1(t)}(w_+(t))$ and $\delta_t$ is positively
folded. One can now argue as at the end of the proof of Proposition~\ref{pr:RetractingSeg}. One
obtains positive roots
$\beta_1,\dots,\beta_s$ such that conditions (i) and (ii) of Definition~\ref{de:LS} are fulfilled; 
condition (iii) is then a consequence of
$s=\ell_{\pi_1(t)}(w_-(t))-\ell_{\pi_1(t)}(w_+(t))$ and
$\pi_1$ is a LS path.

Conversely, if $\pi_1$ is a LS path, the construction of $\delta_t$ as in Theorem~\ref{th:RetractingOnHecke} may be performed by using a set $(\beta_1,\dots,\beta_s)$ of positive roots with
$s = \ell_{\pi_1(t)}(w_-(t)) - \ell_{\pi_1(t)}(w_+(t))$. This gallery $\delta_t$ is 
folded exactly $s$ times (positively and along true walls), its length is 
$n=\ell_{\pi_1(t)}(w_-(t))$; so, once we get rid of the stutterings, we get a
minimal gallery $\delta_t^{ns}$ from $c_0$ to
$\pi_{1+}(t)=\pi_1(t)+[0,1)w_+(t)\lambda$. Hence, as the foldings were positive,
$$neg(\delta_t)=neg(\delta_t^{ns})=\ell_{\pi_1(t)}(w_+(t)),$$ 
so  
${\rm codim}(\tilde\pi_1)={\rm codim}(\pi_1)$.
\end{proof}

\subsection{Preorder relation on the hovel} 
\label{sse:Preorder}

\begin{theo} 
\label{th:Preorder}
On the hovel $\mathcal I$, the relation $\leq$ 
(defined in \ref{sse:Applications}.2) is a preorder relation. More precisely, if three different 
points $x$, $y$  and $z$ in $\mathcal I$ are such that $x\leq y$ and $y\leq z$ then
$x\leq z$ and, in particular, $x$ and 
$z$ are in a same apartment.
\end{theo}

\begin{rema} 
\label{re:Preorder}
This result precises the 
structure of the hovel $\mathcal I$. It is a generalization of Lemme 7.3.6 in
\cite{BruhatTits72}.  It may also be seen as a generalization of the Cartan 
decomposition proved by Garland for  $p-$adic loop groups \cite{Garland95}, even if
it is weaker than this decomposition in the affine  case. As Garland asserts, the
Cartan decomposition holds only after some twisting; this is more or less
equivalent to the fact that not any two points in $\mathcal I$ are in a same apartment.

More precisely, let us look at the simplest affine Kac-Moody group  $G = \SL_2^{(1)}$. If $K = \mathbb C(\!(\varpi)\!)$ and $\mathcal O = \mathbb C[\![\varpi]\!]$, then, up to the center (which is in $T(K)$), the group $G(K) = \SL_2^{(1)}(\mathbb C(\!(\varpi)\!))$ is a semidirect product $G(K) = K^* \ltimes \SL_2(K[u,u^{-1}])$, with $K^*\subset T(K)$. We saw in Example~\ref{ex:G(O)} that $G_0 = \hat P_0 = G(\mathcal O) = \mathcal O^*\ltimes \SL_2(\mathcal O[u,u^{-1}])$ (up to the center). The Cartan decomposition would tell that $G(K) = G(\mathcal O) T(K) G(\mathcal O)$, hence $\SL_2(K[u,u^{-1}]) = \SL_2(\mathcal O[u,u^{-1}]) T_1(K)  \SL_2(\mathcal O[u,u^{-1}])$, where $T_1(K)$ is the torus of diagonal matrices in $\SL_2(K)$. The four coefficients of a matrix in the right hand side span the same sub$-\mathcal O[u,u^{-1}]-$module of $K[u,u^{-1}]$ as a matrix in $T_1(K)$. So, this module is generated by $h$ and $h^{-1}$, i.e. by a single element of $K^*$. But, the $\mathcal O[u,u^{-1}]-$module spanned by the coefficients of $g = \begin{pmatrix}1& \varpi^{-1}(1+u)\\ 0 & 1\end{pmatrix}$ is not generated by a single element. Hence, Cartan decomposition fails and the points $0$ and $g.0$ are not in a same apartment. If they were in an apartment $A'$, then $A' = h_1.A_f$ and $g^{-1}.A' = h_2.A_f$ for $h_1,h_2\in G_0$ (by Remark~\ref{ex:Independence}), hence $h_1^{-1}gh_2\in N(K)$ (by \ref{sse:Applications}.1), contradiction!
\end{rema}

\begin{lemm} 
\label{le:Preorder}
In the situation of Theorem \ref{th:Preorder}, there exists an 
apartment $A$ containing $x$ and $[y,z)$. Moreover, this apartment is unique up
to isomorphism.
\end{lemm}

\begin{proof}
In an apartment $A_1$ containing $x$ and $y$, there exists a 
vectorial chamber $C^v$ such that $x\in y+\overline{C^v}$.
Moreover, there exists an apartment~$A$ containing $\mathfrak S = {\rm germ} (y + 
\overline{C^v})$ and $[y,z)$; this  apartment also contains $y+ \overline
{C^v}\ni x$ (cf. \ref{sse:Applications}.3). The uniqueness is a consequence of \ref{sse:ExGood}.2, \ref{sse:ExGood}.4, Proposition~\ref{pr:Good} 4) and Remark~\ref{ex:Independence}.
\end{proof}

\begin{enonce*}[remark]{Proof of the theorem}\ 

\begin{enumerate}
\item[1)] For $z'\in[y,z[$ such that $x\leq z'$, we choose an apartment $A$ 
containing $[z',x)$ and $[z',z)$ (\ref{sse:Applications}.4); this apartment has an associated 
system of real roots $\Phi(A)$ and we define the finite set $\Phi(z')$ of the 
roots $\alpha\in\Phi(A)$ such that $\alpha(z')>\alpha(x_1)$ and $\alpha(z')>\alpha(z_1)$ for 
some $x_1\in[x,z']\cap A$ and some $z_1\in[z,z']\cap A$. As $[z',x)$ and 
$[z',z)$ are generic, \ref{sse:Applications}.4) shows that $\Phi(z')$  depends, up to isomorphism,
only on $[z',x)$ and $[z',z)$ but not on $A$.

By Lemma \ref{le:Preorder} (mutatis mutandis), there exists an apartment $A_{z'}$ 
containing $z$ and $[z',x)$ and this apartment is unique up to isomorphism. We 
define $N_{z'}$  as the finite number of walls (in $A_{z'}$) of direction
Ker$(\alpha)$ for some $\alpha\in\Phi(z')$ and  separating $z'$ from $z$ (in a strict
sense). We shall argue by induction on 
$(\vert\Phi(y)\vert,N_y)$ (with lexicographical order).

\item[2)] By Lemma \ref{le:Preorder}, there is an apartment $A_1$ containing $x$ and $[y,z)$. 
We choose a  vectorial chamber $C^v$ in $A_1$ such that its associated system of
positive roots $\Phi^+(C^v)$ contains the roots $\alpha\in\Phi(A_1)$ such that 
$\alpha(y)>\alpha(x)$ or $\alpha(y)=\alpha(x)$ and $\alpha(z_1)>\alpha(y)$ (for some 
$z_1\in[y,z]\cap A_1$); in particular $[x,y]\subset y-\overline{C^v}$. Now if 
$\alpha\in\Phi^+(C^v)$ is such that $\alpha(z_1)<\alpha(y)$  (for some $z_1\in[y,z]\cap
A_1$) then $\alpha(y)>\alpha(x)$; hence $\Phi(y)$ (computed in $A_1$)  is the set of
roots $\alpha\in\Phi^+(C^v)$ such that $\alpha(z_1)<\alpha(y)$  (for some $z_1\in[y,z]\cap
A_1$).

Let $\mathfrak S$ be the sector-germ associated to $-C^v$ in $A_1$ and $\rho$ the 
retraction of center $\mathfrak S$ onto $A_1$.

\item[3)] Suppose $z_1\in]y,z]$ is such that no wall (in $A_y$ or any apartment 
containing~$z_1$  and $[y,x)$) of direction Ker$\alpha$ for some $\alpha\in\Phi(y)$
separates $y$ from~$z_1$. We shall  prove that the enclosure of $\mathfrak S$ and $z_1$
contains $y$ and $x$, so, there is an apartment  containing $x$, $y$, $z_1$ and~
$\mathfrak S$. Hence, the theorem is true if $z_1=z$ and this gives the first step of the
induction when $\Phi(y)=\emptyset$ or $N_y=0$.

As in \ref{sse:Retraction}) we get a sequence $y_0=y,\,y_1,\dots,\,y_n=z_1\in[y,z_1]$ and 
apartments $A_1,A_2,\dots,A_n$ such that $A_i$ contains $\mathfrak S$ and $[y_{i-1},y_i]$. 
The characterization of $\Phi(y)$ in 2) above and the hypothesis on the walls 
prove that~$y$ is in the enclosure of $y_1$ and $\mathfrak S$, then $x$ is also in this 
enclosure. So $A_2$, which contains $y_1$ and $\mathfrak S$  contains also $y$ and~$x$. We
can replace $y_1$ by $y_2$ and $A_1$ by $A_2$; by induction  on $n$ we obtain the
result of 3).

\item[4)] We choose for $z_1\in]y,z]$ the point satisfying the hypothesis of 3) 
which is the nearest to $z$, it exists as $\Phi(y)$ and $N_y$ are finite. We may
(and do) suppose $z_1\neq z$.  We choose for $A_1$ the apartment containing $x$,
$y$, $z_1$ and the $\mathfrak S$ explained in 3).  An apartment $A_2$ containing $\mathfrak S$ and
$[z_1,z)$ is sent isomorphically by $\rho$ onto
$A_1$.  This enables us to identify $\Phi(z_1)$ with the set $\Phi'(z_1)$ of the 
roots $\alpha\in\Phi(A_1)$ such that $\alpha(z_1)>\alpha(x)$ (hence $\alpha\in\Phi^+(C^v)$) and 
$\alpha(z_1)>\alpha(\rho z_2)$ (for some $z_2\in [z_1,z]$ near $z_1$). By Proposition~\ref{pr:RetractingSeg}, $\rho([z_1,z))=z_1+[0,1)w^+\lambda\,$, $[y,z_1)=y+[0,1)w^-\lambda$ for some 
$\lambda\in{\overline{C^v}}$  and $w^+,w^-\in W^v$ such that $w^+\leq w^-$. In
particular, for $\alpha\in\Phi^+(C^v)$, 
$\alpha(z_1)>\alpha(\rho z_2)$ means $\alpha(w^+\lambda)<0$, so 
$\Phi'(z_1)\subset\{\alpha\in\Phi^+(C^v)\mid\alpha(w^+\lambda)<0\}$ and (as $w^+$ is chosen 
minimal)  this set is of cardinal $\ell(w^+)$. Now we saw in 2) that 
$\Phi(y)=\{\alpha\in\Phi^+(C^v)\mid\alpha(w^-\lambda)<0\}$. Hence, as $w^+\leq w^-$, 
$\vert\Phi'(z_1)\vert\leq\ell(w^+)\leq\ell(w^-)\leq\vert \Phi(y)\vert$.

 If $\vert\Phi'(z_1)\vert<\vert\Phi(y)\vert$ the theorem is true by induction. 
Otherwise, the four numbers above are equal; in particular, as  
$w^+\leq w^-$, one has $w^+=w^-$ and $\Phi'(z_1)=\Phi(y)$. 

 We consider the segment $[y,z]$ as a linear path $\pi:[0,1]\rightarrow[y,z]$, 
$\pi(0)=y$, $\pi(1)=z$, $z_1=\pi(t_1)$ and $\pi'(t)\in W^v\lambda$, $\forall t$. The 
number $N_{z_1}$ is  calculated in an apartment $A_{z_1}$ containing $z$ and
$[z_1,x)$ using $\Phi(z_1)$ and 
$[z_1,z]$. We may suppose $A_{z_1}$ also containing ${\rm germ}_{z_1}(z_1-C^v)$; then 
there is  an isomorphism from $(A_{z_1},\Phi(z_1),[z_1,z])$ to
$(A_1,\Phi'(z_1),[z_1,Z])$ where 
$Z=z_1+(1-t_1)w^+\lambda$, so $N_{z_1}$ may be computed with this last triple. Arguing 
the same way, we see that $N_y$ may be calculated with the triple 
$(A_1,\Phi(y),[y,Z'])$  with $Z'=y+w^-\lambda$. Actually, $z_1=y+t_1.w^-\lambda$ and we
saw that $\Phi'(z_1)=\Phi(y)$ and $w^+\lambda=w^-\lambda$, so $Z'=Z$. Moreover, by the 
choice of $z_1$, there is a wall of  direction Ker$(\alpha)$ for some $\alpha\in\Phi(y)$
containing $z_1$. Hence, $N_{z_1}<N_y$, and this proves the theorem by induction.
\end{enumerate}
  
\end{enonce*}

\nocite{*}
\bibliography{GausRous}
\bibliographystyle{cdraifplain}

\end{document}